%
%
%
%
%
%
%
%
\magnification=\magstephalf      
%
%
\vsize=7.5truein                 
\hsize=5.2truein                 
\newskip\stdskip                 
\stdskip=6pt plus3pt minus3pt    
\medskipamount=\stdskip          
\parindent=0pt                   
\parskip=\stdskip                
\abovedisplayskip=\stdskip       
\belowdisplayskip=\stdskip       
\mathsurround=0.75pt             
\overfullrule=0pt                
%
%
\def\ppar{\par\goodbreak\vskip 8pt plus 4pt minus 4pt}     
%
%
\def\stdspace{\hskip 0.75em plus 0.15em\ignorespaces}
\let\qua\stdspace 
%
%
%
%
%
%
%
\def\hexnumber#1{\ifcase#1 0\or 1\or 2\or 3\or 4\or 5\or 6\or 7\or 8\or
 9\or A\or B\or C\or D\or E\or F\fi}
%
%
\font\thirtnmsa=msam10 scaled 1315    
\font\tenmsa=msam10          \font\ninemsa=msam9
\font\sevenmsa=msam7         \font\sixmsa=msam6
\font\fivemsa=msam5
%
%
\newfam\msafam                  \textfont\msafam=\tenmsa
\scriptfont\msafam=\sevenmsa    \scriptscriptfont\msafam=\fivemsa
\edef\hexa{\hexnumber\msafam}        
\def\msa{\fam\msafam\tenmsa}         
%
%
\font\thirtnmsb=msbm10 scaled 1315   
\font\tenmsb=msbm10      \font\ninemsb=msbm9
\font\sevenmsb=msbm7     \font\sixmsb=msbm6
\font\fivemsb=msbm5
%
\newfam\msbfam                   \textfont\msbfam=\tenmsb       
\scriptfont\msbfam=\sevenmsb     \scriptscriptfont\msbfam=\fivemsb
\edef\hexb{\hexnumber\msbfam}    
\def\msb{\fam\msbfam\tenmsb}     
%
%
\font\thirtneufm=eufm10 scaled 1315   
\font\teneufm=eufm10                 \font\nineeufm=eufm9
\font\seveneufm=eufm7                \font\sixeufm=eufm6
\font\fiveeufm=eufm5
%
\newfam\eufmfam                    \textfont\eufmfam=\teneufm
\scriptfont\eufmfam=\seveneufm     \scriptscriptfont\eufmfam=\fiveeufm
\edef\hexf{\hexnumber\eufmfam}      
\def\frak{\fam\eufmfam\teneufm}     
%
%
%
\font\thirtnrm=cmr10 scaled 1315    
\font\ninerm=cmr9                   \font\sixrm=cmr6   
%
\font\thirtni=cmmi10 scaled 1315    
\font\ninei=cmmi9                   \font\sixi=cmmi6  
%
\font\thirtnsy=cmsy10 scaled 1315   
\font\ninesy=cmsy9                  \font\sixsy=cmsy6  
%
\font\thirtnbf=cmbx10 scaled 1315   
\font\ninebf=cmbx9                  \font\sixbf=cmbx6  
%
%
\font\thirtnex=cmex10 scaled 1315   
\font\nineex=cmex9                  
%
%
\font\thirtnit=cmti10 scaled 1315  
\font\nineit=cmti9                  
%
\font\thirtnsl=cmsl10 scaled 1315  
\font\ninesl=cmsl9                  
%
\font\thirtntt=cmtt10 scaled 1315  
\font\ninett=cmtt9                  
%
%
%
%
\def\small{%
%
%
\textfont0=\ninerm \scriptfont0=\sixrm \scriptscriptfont0=\fiverm
\def\rm{\fam0\ninerm}
%
%
\textfont1=\ninei \scriptfont1=\sixi \scriptscriptfont1=\fivei
%
%
\textfont2=\ninesy \scriptfont2=\sixsy \scriptscriptfont2=\fivesy
%
%
\textfont3=\nineex \scriptfont3=\nineex \scriptscriptfont3=\nineex
%
%
\textfont\bffam=\ninebf \scriptfont\bffam=\sixbf
\scriptscriptfont\bffam=\fivebf \def\bf{\fam\bffam\ninebf}%
%
%
\textfont\itfam=\nineit \def\it{\fam\itfam\nineit}%
\textfont\slfam=\ninesl \def\sl{\fam\slfam\ninesl}%
\textfont\ttfam=\ninett \def\tt{\fam\ttfam\ninett}%
%
%
%
\textfont\msafam=\ninemsa \scriptfont\msafam=\sixmsa
\scriptscriptfont\msafam=\fivemsa \def\msa{\fam\msafam\ninemsa}%
%
%
\textfont\msbfam=\ninemsb \scriptfont\msbfam=\sixmsb
\scriptscriptfont\msbfam=\fivemsb \def\msb{\fam\msbfam\ninemsb}%
%
%
\textfont\eufmfam=\nineeufm  \scriptfont\eufmfam=\sixeufm
\scriptscriptfont\eufmfam=\fiveeufm \def\frak{\fam\eufmfam\nineeufm}%
%
%
%
\normalbaselineskip=11pt%
\setbox\strutbox=\hbox{\vrule height8pt depth3pt width0pt}%
%
%
\normalbaselines\rm
%
%
\stdskip=4pt plus2pt minus2pt    
\medskipamount=\stdskip          
\parskip=\stdskip                
\abovedisplayskip=\stdskip       
\belowdisplayskip=\stdskip       
\def\ppar{\par\goodbreak\vskip 6pt plus 3pt minus 3pt}%
%
%
\def\section##1{\global\advance\sectionnumber by 1
\vskip-\lastskip\penalty-800\vskip 20pt plus10pt minus5pt 
\egroup{\bf\number\sectionnumber\quad##1}\bgroup\small         
\vskip 6pt plus3pt minus3pt
\nobreak\resultnumber=1}
}    
%
\def\beginsmall{\bgroup\small}
\let\endsmall\egroup
%
%
%
%
\def\large{%
\textfont0=\thirtnrm \scriptfont0=\ninerm \scriptscriptfont0=\sevenrm
\def\rm{\fam0\thirtnrm}%
\textfont1=\thirtni \scriptfont1=\ninei \scriptscriptfont1=\seveni
\textfont2=\thirtnsy \scriptfont2=\ninesy \scriptscriptfont2=\sevensy
\textfont3=\thirtnex \scriptfont3=\thirtnex \scriptscriptfont3=\thirtnex
\textfont\bffam=\thirtnbf \scriptfont\bffam=\ninebf
\scriptscriptfont\bffam=\sevenbf \def\bf{\fam\bffam\thirtnbf}%
\textfont\itfam=\thirtnit \def\it{\fam\itfam\thirtnit}%
\textfont\slfam=\thirtnsl \def\sl{\fam\slfam\thirtnsl}%
\textfont\ttfam=\thirtntt \def\tt{\fam\ttfam\thirtntt}%
\textfont\msafam=\thirtnmsa \scriptfont\msafam=\ninemsa
\scriptscriptfont\msafam=\sevenmsa \def\msa{\fam\msafam\thirtnmsa}%
\textfont\msbfam=\thirtnmsb \scriptfont\msbfam=\ninemsb
\scriptscriptfont\msbfam=\sevenmsb \def\msb{\fam\msbfam\thirtnmsb}%
\textfont\eufmfam=\thirtneufm  \scriptfont\eufmfam=\nineeufm
\scriptscriptfont\eufmfam=\seveneufm \def\frak{\fam\eufmfam\teneufm}%
\normalbaselineskip=16pt%
\setbox\strutbox=\hbox{\vrule height11.5pt depth4.5pt width0pt}%
\normalbaselines\rm}%
%
%
\def\Bbb#1{{\msb#1}}

%
\def\N{{\frak N}}
\def\re{\Bbb R}
%
\mathchardef\plussquare="0\hexa01
\mathchardef\nge="3\hexb0B
\mathchardef\maltesecross="0\hexa7A
\mathchardef\del="0\hexf01
%
%
%
%
\font\sc=cmcsc10
%
%
%
%
\def\sqr#1#2{{\vcenter{\vbox{\hrule  height.#2truept
	\hbox{\vrule width.#2truept height#1truept 
	\kern#1truept \vrule width.#2truept}
	\hrule height.#2truept}}}}
\def\sq{\sqr55}    
%
%
%
%
\newcount\sectionnumber            
\newcount\resultnumber             
\sectionnumber=0\resultnumber=1    
%
%
%
\def\section#1{\global\advance\sectionnumber by 1
\xdef\nextkey{\number\sectionnumber}
\vskip-\lastskip\penalty-800\vskip 20pt plus10pt minus5pt 
{\large\bf\number\sectionnumber\quad#1}         
\vskip 8pt plus4pt minus4pt
\nobreak\resultnumber=1}      
%
%
%
%
%
\def\sh#1{\vskip-\lastskip\ppar{\bf #1}\par\nobreak\medskip}         
%
%
%
%

%
\def\proc#1{\xdef\nextkey{\number\sectionnumber.\number\resultnumber}%
\vskip-\lastskip\ppar\bf%
\noindent#1\ \number\sectionnumber.\number\resultnumber
\stdspace\sl\global\advance\resultnumber by 1\ignorespaces}
\def\endproc{\rm\ppar} 
%
%
\def\prf{\vskip-\lastskip\ppar\noindent{\bf Proof}%
\stdspace\rm}                            
\def\qed{\hfill$\sq$\par\goodbreak\rm}   
\def\endprf{\unskip\stdspace\hbox{}
\hfill$\sq$\par\medskip}                 
%
%
%
%
%
%
%
%
\def\proclaim#1{\vskip-\lastskip\ppar\bf%
\noindent#1\stdspace\sl\ignorespaces} 

%
%
%
%
\def\rk#1{\vskip-\lastskip\ppar{\bf #1}\stdspace\ignorespaces}                

%
%
%
%
%
%
\def\label{\xdef\nextkey{\number\sectionnumber.\number\resultnumber}%
\number\sectionnumber.\number\resultnumber
\global\advance\resultnumber by 1}
%
%
%
%
%
%
%
%
%
%
%
%
%
%
%
%
\newcount\refnumber              
\refnumber=1                     
\long\def\reflist#1\endreflist{%
\long\def\thereflist{#1}{\def\refkey##1##2\par{\xdef##1{\number\refnumber}%
\global\advance\refnumber by 1}%
\def\key##1##2\par{\expandafter\xdef%
\csname##1\endcsname{\number\refnumber}%
\global\advance\refnumber by 1}#1\par}}
\long\def\references{%
\penalty-800\vskip-\lastskip\vskip 15pt plus10pt minus5pt 
{\large\bf References}\ppar 
{\leftskip=25pt\frenchspacing    
\small\parskip=3pt plus2pt       
\def\refkey##1##2\par{\noindent  
\llap{[##1]\stdspace}\ignorespaces##2\par}         
\def\key##1##2\par{\noindent  
\llap{[\ref{##1}]\stdspace}\ignorespaces##2\par}  
\def\,{\thinspace}\thereflist\par}}
%
%
%
\newcount\footnotenumber         
\footnotenumber=1                
\def\fnote#1{\xdef\nextkey{\number\footnotenumber}%
{\small\ifnum\footnotenumber>9\parindent=14pt%
\else\parindent=10pt\fi\footnote{$^{\number\footnotenumber}$}%
{\hglue-5pt#1}\global\advance\footnotenumber by 1}}
%
%
%
%
%
%
%
\newcount\figurenumber          
\figurenumber=1                 
\def\caption#1{\xdef\nextkey{\number\figurenumber}%
\cl{\small Figure \number\figurenumber: #1}%
\global\advance\figurenumber by 1}
\def\figurelabel{\xdef\nextkey{\number\figurenumber}%
\cl{\small Figure \number\figurenumber}%
\global\advance\figurenumber by 1}
\long\def\figure#1\endfigure{{\xdef\nextkey{\number\figurenumber}%
\let\captiontext\relax\def\caption##1{\xdef\captiontext{##1}}%
\midinsert\cl{\ignorespaces#1\unskip\unskip\unskip\unskip}\vglue6pt\cl{\small 
Figure \number\figurenumber\ifx\captiontext\relax\else: \captiontext
\fi}\endinsert\global\advance\figurenumber by 1}}
%
%
%
%
%
%
%
\def\nextkey{??}   
%
\def\key#1{\expandafter\xdef\csname #1\endcsname{\nextkey}}
\def\ref#1{\expandafter\ifx\csname #1\endcsname\relax
\immediate\write16{Reference {#1} undefined}??\else
\csname #1\endcsname\fi}
%
%
%
%
%
%
%
\newread\gtinfile
\newwrite\gtreffile
\def\useforwardrefs{
\openin\gtinfile\jobname.ref
\ifeof\gtinfile
\closein\gtinfile
\immediate\write16{No file \jobname.ref}
\else
\closein\gtinfile
\input \jobname.ref
\fi
\immediate\openout\gtreffile \jobname.ref
%
%
\def\key##1{{\def\\{\noexpand}%
\expandafter\xdef\csname ##1\endcsname{\nextkey}%
\immediate\write\gtreffile{\\\expandafter\\\def\\\csname ##1\\\endcsname%
{\nextkey}}}}
%
%
\long\def\reflist##1\endreflist{%
\long\def\thereflist{##1}{\def\refkey####1####2\par{\xdef####1{%
\number\refnumber}{\def\\{\noexpand}\immediate\write\gtreffile
{\\\def\\####1{\number\refnumber}}}\global\advance\refnumber by 1}%
\def\key####1####2\par{\expandafter\xdef%
\csname####1\endcsname{\number\refnumber}%
{\def\\{\noexpand}\immediate\write\gtreffile
{\\\expandafter\\\def\\\csname ####1\\\endcsname{\number\refnumber}}}
\global\advance\refnumber by 1}##1\par}}
\long\def\biblio##1\endbiblio{\reflist##1\endreflist\references}%
%
%
\def\numkey##1{{\def\\{\noexpand}%
\xdef##1{\number\sectionnumber.\number\resultnumber}
\immediate\write\gtreffile{\\\def\\##1%
{\number\sectionnumber.\number\resultnumber}}}}
\def\seckey##1{{\def\\{\noexpand}\xdef##1{\number\sectionnumber}
\immediate\write\gtreffile{\\\def\\##1{\number\sectionnumber}}}}
\def\figkey##1{\xdef##1{\number\figurenumber}%
{\def\\{\noexpand}\immediate\write\gtreffile%
{\\\def\\##1{\number\figurenumber}}}
\number\figurenumber\global\advance\figurenumber by 1}
}   
%
%
%
%
\def\figkey#1{\xdef#1{\number\figurenumber}%
\number\figurenumber\global\advance\figurenumber by 1}
\def\fig#1#2\endfig{%
\midinsert\cl{#2}\vglue6pt\cl{\small Figure #1}\endinsert}
\def\newfig{\number\figurenumber\global\advance\figurenumber by 1}
\def\numkey#1{\xdef#1{\number\sectionnumber.\number\resultnumber}}
\def\seckey#1{\xdef#1{\number\sectionnumber}}
%
%
%
%
%
%
%
%
%
\def\verb{\catcode`\"=\active}       
\def\brev{\catcode`\"=12}            
\brev                                
\verb                                
{\obeyspaces\gdef {\ }}              
{\catcode`\`=\active\gdef`{\relax\lq}}
\def"{%
\begingroup\baselineskip=12pt\def\par{\leavevmode\endgraf}%
\tt\obeylines\obeyspaces\parskip=0pt\parindent=0pt%
\catcode`\$=12\catcode`\&=12\catcode`\^=12\catcode`\#=12%
\catcode`\_=12\catcode`\~=12%
\catcode`\{=12\catcode`\}=12\catcode`\%=12\catcode`\\=12%
\catcode`\`=\active\let"\endgroup}
\brev      
%
%
%
%
%
%
\def\items{\par\leftskip = 25pt}           
\def\enditems{\par\leftskip = 0pt}         
\def\item#1{\par\leavevmode\llap{#1\stdspace}%
\ignorespaces}                             
%
%

%
%
\def\co{\colon\thinspace}    
\def\np{\vfil\eject}         
\def\nl{\hfil\break}         
\def\cl{\centerline}         
\def\gt{{\mathsurround=0pt\it $\cal G\mskip-2mu$eometry \&\ 
$\cal T\!\!$opology}}        
\def\gtm{{\mathsurround=0pt\it $\cal G\mskip-2mu$eometry \&\ 
$\cal T\!\!$opology $\cal M\mskip-1mu$onographs}}    
\def\agt{{\mathsurround=0pt\it$\cal A\mskip-.7mu$lgebraic \&\ 
$\cal G\mskip-2mu$eometric $\cal T\!\!$opology}}  
%
%
%

%
%
%
%
%
\def\title#1{\def\thetitle{#1}}

\def\author#1{\edef\previousauthors{\theauthors}
 \ifx\theauthors\relax\def\theauthors{#1}\else
 \def\theauthors{\previousauthors\par#1}\fi}

%
\def\address#1{\edef\previousaddresses{\theaddress}
 \ifx\theaddress\relax\def\theaddress{#1}\else
 \def\theaddress{\previousaddresses\par\vskip 2pt\par#1}\fi}
\def\secondaddress#1{\edef\previousaddresses{\theaddress}
 \ifx\theaddress\relax\def\theaddress{#1}\else
 \def\theaddress{\previousaddresses\par{\rm and}\par#1}\fi}   

\def\email#1{\edef\previousemails{\theemail}
 \ifx\theemail\relax\def\theemail{#1}\else
 \def\theemail{\previousemails\hskip 0.75em\relax#1}\fi}
\def\secondemail#1{\edef\previousemails{\theemail}
 \ifx\theemail\relax\def\theemail{#1}\else
 \def\theemail{\previousemails\hskip 0.75em{\rm and}\hskip 0.75em
 \relax#1}\fi}
\def\url#1{\edef\previousurls{\theurl}
 \ifx\theurl\relax\def\theurl{#1}\else
 \def\theurl{\previousurls\hskip 0.75em\relax#1}\fi}
\def\secondurl#1{\edef\previousurls{\theurl}
 \ifx\theurl\relax\def\theurl{#1}\else
 \def\theurl{\previousurls\hskip 0.75em{\rm and}\hskip 0.75em
 \relax#1}\fi}
\long\def\abstract#1\endabstract{\long\def\theabstract{#1}}
\def\primaryclass#1{\def\theprimaryclass{#1}}
\def\secondaryclass#1{\def\thesecondaryclass{#1}}
\def\keywords#1{\def\thekeywords{#1}}
%
%
\let\\\par\let\thetitle\relax\let\theshorttitle\relax
\let\theauthors\relax\let\theshortauthors\relax
\let\theaddress\relax\let\theshortaddress\relax
\let\theemail\relax\let\theurl\relax
\let\theabstract\relax\let\theprimaryclass\relax
\let\thesecondaryclass\relax\let\thekeywords\relax
%
%
%
%
\long\def\maketitlepage{    

\vglue 0.2truein   

%
{\parskip=0pt\leftskip 0pt plus 1fil\def\\{\par\smallskip}{\large
\bf\thetitle}\par\medskip}   

\vglue 0.15truein 

%
{\parskip=0pt\leftskip 0pt plus 1fil\def\\{\par}{\sc\theauthors}
\par\medskip}%
 
\vglue 0.1truein 

%
{\small\parskip=0pt
{\leftskip 0pt plus 1fil\def\\{\par}{\sl\theaddress}\par}
\ifx\theemail\relax\else  
\vglue 5pt \def\\{\stdspace{\rm and}\stdspace} 
\cl{Email:\stdspace\tt\theemail}\fi
\ifx\theurl\relax\else    
\vglue 5pt \def\\{\stdspace{\rm and}\stdspace} 
\cl{URL:\stdspace\tt\theurl}\fi\par}

\vglue 7pt 

{\bf Abstract}

\vglue 5pt

\theabstract

\vglue 7pt 

{\bf AMS Classification numbers}\quad Primary:\quad \theprimaryclass\par

Secondary:\quad \thesecondaryclass

\vglue 5pt 

{\bf Keywords:}\quad \thekeywords

\np  

}    
%
%
\long\def\makeshorttitle{    


%
{\parskip=0pt\leftskip 0pt plus 1fil\def\\{\par\smallskip}{\large
\bf\thetitle}\par\medskip}   

\vglue 0.05truein 

%
{\parskip=0pt\leftskip 0pt plus 1fil\def\\{\par}{\sc\theauthors}
\par\medskip}%
 
\vglue 0.03truein 

%
{\small\parskip=0pt
{\leftskip 0pt plus 1fil\def\\{\par}{\sl\ifx\theshortaddress\relax
\theaddress\else\theshortaddress\fi}\par}
\ifx\theemail\relax\else  
\vglue 5pt \def\\{\stdspace{\rm and}\stdspace} 
\cl{Email:\stdspace\tt\theemail}\fi
\ifx\theurl\relax\else    
\vglue 5pt \def\\{\stdspace{\rm and}\stdspace} 
\cl{URL:\stdspace\tt\theurl}\fi\par}

\vglue 10pt 


{\small\leftskip 25pt\rightskip 25pt{\bf Abstract}\stdspace\theabstract

{\bf AMS Classification}\stdspace\theprimaryclass
\ifx\thesecondaryclass\relax\else; \thesecondaryclass\fi\par
{\bf Keywords}\stdspace \thekeywords\par}
\vglue 7pt
}    
\let\maketitle\makeshorttitle        
%
%

\def\volumenumber#1{\def\thevolumenumber{#1}}
\def\volumename#1{\def\thevolumename{#1}}
\def\volumeyear#1{\def\thevolumeyear{#1}}
\def\pagenumbers#1#2{\def\startpage{#1}\def\finishpage{#2}}
\def\published#1{\def\publishdate{#1}}

\volumenumber{X}
\volumename{Volume name goes here}
\volumeyear{20XX}
\pagenumbers{1}{XXX}
\published{XX Xxxember 20XX}

\long\def\makegtmontitle{   

\count0=\startpage

\gtm\nl        
{\small Volume \thevolumenumber: \thevolumename\nl 
Pages \startpage--\finishpage\nl}

\vglue 0.1truein   

{\parskip=0pt\leftskip 0pt plus 1fil\def\\{\par\smallskip}{\large
\bf\thetitle}\par\medskip}   
\vglue 0.05truein 

%
{\parskip=0pt\leftskip 0pt plus 1fil\def\\{\par}{\sc\theauthors}
\par\medskip}%
 
\vglue 0.03truein 


{\small\leftskip 25pt\rightskip 25pt{\bf Abstract}\stdspace\theabstract

{\bf AMS Classification}\stdspace\theprimaryclass
\ifx\thesecondaryclass\relax\else; \thesecondaryclass\fi\par
{\bf Keywords}\stdspace \thekeywords\par}\vglue 7pt

}   

\long\def\makeagttitle{   
\agt\hfill      
\hbox to 60truept{\vbox to 0pt{\vglue -14truept{\bf [Logo here]}\vss}\hss}
\break
{\small Volume \thevolumenumber\ (\thevolumeyear)
\startpage--\finishpage\nl
Published: \publishdate}

\vglue .2truein

{\parskip=0pt\leftskip 0pt plus 1fil\def\\{\par\smallskip}{\large
\bf\thetitle}\par\medskip}   
\vglue 0.05truein 

%
{\parskip=0pt\leftskip 0pt plus 1fil\def\\{\par}{\sc\theauthors}
\par\medskip}%
 
\vglue 0.03truein 


{\small\leftskip 25truept\rightskip 25truept{\bf Abstract}\stdspace\theabstract

{\bf AMS Classification}\stdspace\theprimaryclass
\ifx\thesecondaryclass\relax\else; \thesecondaryclass\fi\par
{\bf Keywords}\stdspace \thekeywords\par}\vglue 7truept

}   


\def\Addresses{\bigskip
{\small \parskip 0pt \leftskip 0pt \rightskip 0pt plus 1fil \def\\{\par}
\sl\theaddress\par\medskip \rm Email:\stdspace\tt\theemail\par
\ifx\theurl\relax\else\smallskip \rm URL:\stdspace\tt\theurl\par\fi}}

\def\agtart{
\hoffset 14truemm
\voffset 31truemm
\font\phead=cmsl9 scaled 950
\font\pnum=cmbx10 scaled 913
\font\pfoot=cmsl9 scaled 950
\headline{\vbox to 0pt{\vskip -4.5mm\line{\small\phead\ifnum
\count0=\startpage ISSN numbers are printed here
\hfill {\pnum\folio}\else\ifodd\count0\def\\{ }%
\ifx\theshorttitle\relax\thetitle\else\theshorttitle\fi\hfill{\pnum\folio}
\else\def\\{ and }{\pnum\folio}\hfill\ifx\theshortauthors\relax\theauthors
\else\theshortauthors\fi\fi\fi}\vss}}
\footline{\vbox to 0pt{\vglue 0mm\line{\small\pfoot\ifnum\count0=\startpage
Copyright declaration is printed here\hfill\else
\agt, Volume \thevolumenumber\ (\thevolumeyear)\hfill\fi}\vss}}
\let\maketitle\makeagttitle\let\makeshorttitle\makeagttitle
\let\maketitlepage\makeagttitle}

\def\gtmonart{
\hoffset 14truemm
\voffset 31truemm
\font\phead=cmsl9 scaled 950
\font\pnum=cmbx10 scaled 913
\font\pfoot=cmsl9 scaled 950
\headline{\vbox to 0pt{\vskip -4.5mm\line{\small\phead\ifnum
\count0=\startpage ISSN numbers are printed here
\hfill {\pnum\folio}\else\ifodd\count0\def\\{ }%
\ifx\theshorttitle\relax\thetitle\else\theshorttitle\fi\hfill{\pnum\folio}
\else\def\\{ and }{\pnum\folio}\hfill\ifx\theshortauthors\relax\theauthors
\else\theshortauthors\fi\fi\fi}\vss}}
\footline{\vbox to 0pt{\vglue 0mm\line{\small\pfoot\ifnum\count0=\startpage
Copyright declaration is printed here\hfill\else
\gtm, Volume \thevolumenumber\ (\thevolumeyear)\hfill\fi}\vss}}
\let\maketitle\makegtmontitle\let\makeshorttitle\makegtmontitle
\let\maketitlepage\makegtmontitle}

\def\gtart{
\hoffset 14truemm
\voffset 31truemm
\font\phead=cmsl9 scaled 950
\font\pnum=cmbx10 scaled 913
\font\pfoot=cmsl9 scaled 950
\headline{\vbox to 0pt{\vskip -4.5mm\line{\small\phead\ifnum
\count0=\startpage ISSN numbers are printed here
\hfill {\pnum\folio}\else\ifodd\count0\def\\{ }%
\ifx\theshorttitle\relax\thetitle\else\theshorttitle\fi\hfill{\pnum\folio}
\else\def\\{ and }{\pnum\folio}\hfill\ifx\theshortauthors\relax\theauthors
\else\theshortauthors\fi\fi\fi}\vss}}
\footline{\vbox to 0pt{\vglue 0mm\line{\small\pfoot\ifnum\count0=\startpage
Copyright declaration is printed here\hfill\else
\gt, Volume \thevolumenumber\ (\thevolumeyear)\hfill\fi}\vss}}
\let\maketitle\maketitlepage\let\makeshorttitle\maketitlepage}

\input pictex
\chardef\newinsCatAt\the\catcode `\@
\catcode `\@=11
%
%
%
\newskip\insertskipamount\newskip\inserthardskipamount
\insertskipamount 12pt plus2pt  
\inserthardskipamount 4pt       
\def\insertskip{\vskip\insertskipamount}
%
%
\newskip\LastSkip
\def\SaveLastSkip{\LastSkip\lastskip}
\def\RestoreLastSkip{\nobreak\vskip-\LastSkip\vskip\LastSkip}
%
%
\newcount\SplitTest
\def\SetSplitTest{\SplitTest\insertpenalties
  \insert\topins{\floatingpenalty1}%
  \advance\SplitTest-\insertpenalties}
%
%
\def\midinsert{\par
 \SaveLastSkip\penalty-150\SetSplitTest\RestoreLastSkip
 \ifnum\SplitTest=-1
  \@midfalse\p@gefalse\else\@midtrue\fi\@ins}
\def\@ins{\par\begingroup\setbox\z@\vbox\bgroup%
  \vglue\inserthardskipamount}
\def\endinsert{\egroup 
  \if@mid \dimen@\ht\z@ \advance\dimen@\dp\z@
    \advance\dimen@\insertskipamount
    \advance\dimen@\pagetotal\advance\dimen@-\pageshrink
    \ifdim\dimen@>\pagegoal\@midfalse\p@gefalse\fi\fi
  \if@mid%
    \ifdim\lastskip<\insertskipamount\removelastskip\insertskip\fi
    \nointerlineskip\box\z@\penalty-200\insertskip
  \else%
    \SaveLastSkip
    \insert\topins{\penalty100 
    \splittopskip\z@skip
    \splitmaxdepth\maxdimen \floatingpenalty\z@
    \ifp@ge \dimen@\dp\z@
    \vbox to\vsize{\unvbox\z@\kern-\dimen@}
    \else \box\z@\nobreak\insertskip\fi}
    \RestoreLastSkip
   \fi\endgroup}
%
\catcode `\@=\newinsCatAt

\let\ninepoint\small
\hoffset 0.5truein     
\voffset 1truein       

\font\spec=cmtex10 scaled\magstephalf  
\font\braid=tiny


\def\bar{\overline}
\def\text{\hbox}

\def\row#1#2#3{(#1_#2,\ldots,#1_#3)}
\def\a{\longrightarrow}
\def\inv{^{-1}}
\def\co{\colon\thinspace}
\def\ep{\epsilon}

\def\d{\hbox{\spec \char'017\kern 0.05em}} 
\def\de{\delta}

\def\la{\lambda}
\let\\=\cr

\def\wreath{\hbox{\tenmsa\lower0.4ex\hbox{\char'161}\kern-0.13em\raise0.2ex\hbox{\char'170}}}

\def\N{{\Bbb N}}

\def\D{{\cal D}}

\def\E{{\Bbb E}}

\def\N{{\cal N}}

\def\Z{{\Bbb Z}}
\def\F{{\cal F}}
\def\re{{\Bbb R}}

\def\ha{{\scriptstyle{1\over2}}}

\def\cupprod{{\scriptscriptstyle\,\cup\,}}
\def\immersed{{\msb\char'156}}
\def\VL{{\cal VL}}
\def\SL{{\cal SL}}

\def\sqr#1#2{{\vcenter{\hrule width #1.#2pt height.#2pt
\hbox to #1.#2pt{\vrule width.#2pt height#1pt \hss
\vrule width.#2pt}
\hrule width #1.#2pt height.#2pt}}}

\def\sq{\sqr55}
\def\tsq{\sqr57}

\mathchardef\square="0\hexa03
\mathchardef\celt="0\hexa01
\mathchardef\twid="1218
\mathchardef\diam="027D       
\mathchardef\clock="0\hexa08
\mathchardef\aclock="0\hexa09
\mathchardef\rar="0225        
\mathchardef\lar="022D        
\mathchardef\ldar="022E       
\mathchardef\rdar="0226       

\def\mo{\hbox{\hskip 1pt\vrule height 8.5pt depth 2pt width 0.7pt\hskip 1pt}}

\def\b#1{\vbox{\braid\offinterlineskip #1}}
\def\r#1{\raise6pt\hbox{#1}}
\def\sone{\b{ii\nl oo\nl ah\nl fc\nl oo}}
\def\stwo{\b{ah\nl fc\nl oo\nl oo\nl ii}}
\def\sonei{\b{ii\nl oo\nl ed\nl bg\nl oo}}
\def\stwoi{\b{ed\nl bg\nl oo\nl oo\nl ii}}

\long\def\mbox#1{$$\vbox{#1}$$}
\newdimen\unitlength
\unitlength= 1.000cm
\def\makebox(0,0)[l]#1{\hbox{#1}}
\def\symbol#1{\char#1}
\font\tencirc=lcircle10
\def\SetFigFont#1#2#3#4{{\ninepoint$#4$}}
\font\thinlinefont=cmr7

%
%
\reflist

\refkey\BRS {\bf S Buoncristiano}, {\bf C Rourke}, {\bf B Sanderson},
 {\it A geometric approach to homology theory}, London
Math. Soc. Lecture Note Series, no.  18 C.U.P. (1976)

\refkey\Carteretal {\bf J\,S Carter}, {\bf D Jelsovsky}, {\bf S Kamada},
{\bf L Langford}, {\bf M Saito}, {\it Quandle cohomology and state-sum
invarants of knotted curves and surfaces}, {\tt arXiv: math.GT/9906115}

\refkey\CKS {\bf J\,S Carter}, {\bf S Kamada}, {\bf M Saito}, {\it Stable
equivalences of knots on surfaces and virtual knot cobordisms},  
{\tt arXiv:math.GT/0008118}

\refkey\Cerf {\bf J Cerf}, {\it Topologie de certains espaces de plongements},
Bull. Soc. Math. France, 98 (1961) 227--382

\refkey\Dou {\bf A Douady}, {\it Vari\'et\'es \`a bord anguleux et
voisinages tubulaires}, S\'eminaire Henri Cartan, no. 1 (1961--2)

\refkey\Fenn {\bf R Fenn}, {\it Techniques of Geometric Topology},  London
Math. Soc. Lecture Note Series, no.~57 C.U.P. (1983)

\refkey\FJK  
{\bf R Fenn}, {\bf M Jordan}, {\bf L\,H Kauffman}, {\it Biracks and
virtual links}, available from:
{\tt http://www.maths.sussex.ac.uk/Staff/RAF/Maths/loumerc.ps}

\refkey\FeRo {\bf R Fenn}, {\bf C Rourke}, {\it Racks and links in codimension 
two}, J. Knot Theory Ramifications, 1 (1992) 343--406

\refkey\Intro  {\bf R Fenn}, {\bf C Rourke}, {\bf B Sanderson},   {\it An 
Introduction to species and the rack space}, Topics in Knot Theory,
(M\,E Bozh\"uy\"uk, editor), Kluwer Academic Publishers
(1993) 33--55

\refkey\Trunks {\bf R Fenn}, {\bf C Rourke}, {\bf B Sanderson},  {\it 
Trunks and classifying spaces}, Applied Categorical Structures, 3
(1995) 321--356

\refkey\James {\bf R Fenn}, {\bf C Rourke}, {\bf B Sanderson}, 
{\it James bundles and applications}, preprint (1996) {\tt
http://www.maths.warwick.ac.uk/\char'176cpr/ftp/james.ps}

\refkey\JBun {\bf R Fenn}, {\bf C Rourke}, {\bf B Sanderson}, {\it James 
bundles}, submitted for publication, {\tt arXiv:math.AT/0301354}

\refkey\Class {\bf R Fenn}, {\bf C Rourke}, {\bf B Sanderson}, {\it A 
classification of classical links}, (to appear)

\refkey\Flower {\bf J Flower}, {\it Cyclic bordism and rack spaces},
PhD Thesis, University of Warwick (1995)

\refkey\Greene {\bf M Greene}, {\it PhD Thesis}, Warwick (1996), 
available from:\nl {\tt
http://www.maths.warwick.ac.uk/\char'176cpr/ftp/mtg.ps.gz}

\refkey\Gr
{\bf M Gromov}, {\it Partial differential relations}, Springer--Verlag (1986)

\refkey\Hilt {\bf P Hilton}, {\it On the homotopy groups of the union of
spheres}, J. London Math. Soc. 30 (1955) 154--172

\refkey\KamKam {\bf N Kamada}, {\bf S Kamada}, {\it Abstract link diagrams
and virtual knots}, to appear

\refkey\Kauf {\bf L\,H Kauffman},  {\it Virtual knot theory}, Euro. J.
Combin. 7 (1999) 663--690

\refkey\KS {\bf U Koschorke},  {\bf B Sanderson}, {\it Self-intersections and
higher Hopf invariants}, Topology, 17 (1978) 283--290

\refkey\Kup {\bf G Kuperberg}, {\it What is a virtual link?}, {\tt
arXiv:math.GT/0208039}

\refkey\LS {\bf R Lashof}, {\bf S Smale}, {\it Self-intersections of immersed
manifolds}, J. Math. Mech. 8 (1959) 143--157

\refkey\LN
{\bf R\,A Litherland}, {\bf S Nelson}, {\it The Betti numbers of
some finite racks}, J. Pure Appl. Algebra 178 (2003) 187--202,
{\tt arXiv:math.GT/0106165}

\refkey\Pastor {\bf J\,G Pastor}, {\it Bundle complexes and bordism of 
immersions}, PhD Thesis, University of Warwick (1982)

\refkey\PLTop {\bf C Rourke}, {\bf B Sanderson}, {\it 
Introduction to piecewise-linear topology}, Ergebnisse der Mathematik
und ihrer Grenzgebiete, Band 69, Springer--Verlag, New
York--Heidelberg (1972)

\refkey\Comp {\bf C Rourke}, {\bf B Sanderson}, {\it The compression 
theorem: I and II}, Geom. Topol, 5 (2001) 399--429, 431--440

\refkey\Tref {\bf C Rourke}, {\bf B Sanderson}, {\it There are two
2--twist spun trefoils}, {\tt arxiv:\break math.GT/0006062:v1}

\refkey\SMah {\bf B Sanderson}, {\it The geometry of Mahowald orientations},
Algebraic Topology, Aarhus, Springer Lecture Notes in Mathematics,
no. 763 (1978) 152--174

\refkey\SBord {\bf B Sanderson}, {\it Bordism of links in codimension two}, J.
London Math. Soc. 35 (1987) 367--376

\refkey\Whitehead {\bf J\,C\,H Whitehead}, {\it On adding relations to homotopy
groups}, Ann. of Math. 42 (1941) 409--428

\refkey\Wiest {\bf B Wiest}, {\it Rack spaces and loop spaces}, J. Knot 
Theory Ramifications 8 (1999) 99--114

\endreflist

\title{The rack space}
\author{Roger Fenn}
\address{Department of Mathematics, University of Sussex\\
Falmer, Brighton, BN1 9QH, UK\\\bigskip {\sc Colin Rourke\\Brian
Sanderson}\\\medskip Mathematics Institute, University of Warwick\\
Coventry, CV4 7AL, UK}
\email{R.A.Fenn@sussex.ac.uk, cpr {\rm and} bjs@maths.warwick.ac.uk}

\abstract

The main result of this paper is a new classification theorem for
links (smooth embeddings in codimension 2).  The classifying space is
the rack space (defined in [\Trunks]) and the classifying bundle is
the first James bundle (defined in [\JBun]).

We investigate the algebraic topology of this classifying space and
report on calculations given elsewhere.  Apart from defining many new
knot and link invariants (including generalised James--Hopf
invariants), the classification theorem has some unexpected
applications.  We give a combinatorial interpretation for $\pi_2$ of a
complex which can be used for calculations and some new
interpretations of the higher homotopy groups of the 3--sphere.  We
also give a cobordism classification of virtual links.

\endabstract

\primaryclass{55Q40, 57M25}
\secondaryclass{57Q45, 57R15, 57R20, 57R40}


\maketitle

\sectionnumber=-1
\section{Introduction}

The main result of this paper is a classification theorem
for links (smooth embeddings of codimension 2):

\proclaim{Classification Theorem}Let $X$ be a rack.  
Then the rack space $BX$ has the property that $\pi_n(BX)$ 
is in natural bijection with the set of cobordism classes of 
framed submanifolds $L$ of $\re^{n+1}$ of codimension 2
equipped with a homomorphism of the fundamental rack
$\Gamma(L\subset\re^{n+1})$ to $X$.

Moreover there is a smooth mock bundle $\zeta^1(BX)$ over $BX$
which plays the r\^ole of classifying bundle.\rm

\ppar
It is important to note that this is a totally new type of
classification theorem.  Classical cobordism techniques give a
bijection between cobordism classes of framed submanifolds of
$\re^{n+1}$ and $\pi_n(\Omega(S^2))$.  But these techniques are unable
to cope with the extra geometric information given by the homomorphism
of fundamental rack.  For readers unfamiliar with the power of the
rack concept (essentially a rack is a way of encapsulating the
fundamental group and peripheral group system in one simple piece of
algebra) here is a weaker result phrased purely in terms of the
fundamental group:

\proclaim{Corollary}Let $\pi$ be a group.  There is a classifying
space $BC(\pi)$ such that $\pi_n(BC(\pi))$
is in natural bijection with the set of cobordism classes of 
framed submanifolds $L$ of $\re^{n+1}$ of codimension 2
equipped with a homomorphism $\pi_1(\re^{n+1}-L)\to\pi$.

\prf Let $C(\pi)$ be the conjugacy rack of $\pi$ [\FeRo; example 1.3.1
page 349] then a homomorphism $\Gamma(L\subset Q\times\re)\to C(\pi)$
is equivalent to a homomorphism $\pi_1(Q\times\re-L)\to\pi$, see
[\FeRo; corollaries 2.2 and 3.3, pages 354 and 361]. \endprf

There are several ingredients of the proof of the classification theorem.
The classifying space (the rack space) $BX$ is defined in [\Trunks].
In addition we need the geometry of $\sq$--sets developed in 
[\JBun], and in particular the James bundles of a $\sq$--set.
The compression theorem [\Comp] is needed to reduce the codimension
2 problem to the codimension 1 problem of classifying diagrams up
to cobordism and finally we need to develop a theory of smooth
tranversality to, and smooth mock bundles over, a $\sq$--set.  
This is contained in the present paper.

Here is an outline of this paper:

In section 1 ``{\sl Basic definitions\/}'' we recall the definitions
of $\sq$--sets, $\sq$--maps, the rack space and the associated James
complexes of a $\sq$--set.  In section 2 ``{\sl Mock bundles and
transversality\/}'' we define the concept of a mock bundle over a
$\sq$--set (cf [\BRS]) and observe that the James complexes of a
$\sq$--set $C$ define mock bundles $\zeta^i(C)$ which embed as framed
mock bundles in $\mo C\mo\times\re$.  These are the James bundles of
$C$.  We define transversality for a map of a smooth manifold into a
$\sq$--set and prove that any map can be approximated by a transverse
map.  Mock bundles pull-back over transverse maps to yield mock
bundles whose total spaces are manifolds and in particular the first
James bundle pulls back to give a self-transverse immersion of
codimension 1 and the higher James bundles pull-back to give the
multiple point sets of this immersion.

The transversality theorem leads to our first classification theorem
in section~3 ``{\sl Links and diagrams\/}'' namely that a $\sq$--set
$C$ is the classifying space for cobordism classes of link diagrams
labelled by the cubes of $C$.  In the key example in which $C$ is the
rack space $BX$, there is a far simpler description and we deduce the
classification theorem stated above which interprets the homotopy
groups of $BX$ as bordism classes of links with representation of
fundamental rack in $X$.  There are similar interpretations for sets
of homotopy classes of maps of a smooth manifold in $C$ and $BX$ and
for the bordism groups of $C$ and $BX$.

In section 4 ``{\sl The classical case}'' we look in detail at the
lowest non-trivial dimension ($n=2$) where the cobordism classes can
be described as equivalence classes under simple moves.  This gives a
combinatorial description of $\pi_2(C)$ which can be used for
calculations.  To illustrate this we translate the Whitehead
conjecture [\Whitehead] into a conjecture about coloured link
diagrams.  We finish by classifying virtual links (Kauffman [\Kauf],
see also Kuperberg [\Kup]) up to cobordism, in terms of the
2--dimensional homology of the rack space.

The theory of James bundles gives invariants for knots and links for
the following reason.  If $\Gamma$ is the fundamental rack of a link
$L$ then any invariant of the rack space $B\Gamma$ is {\it a fortiori}
an invariant of $L$.  In particular any of the classical algebraic
toplogical invariants of rack spaces are link invariants.  Further
invariants are obtained by considering representations of the
fundamental rack in a small rack and pulling back invariants from the
rack space of this smaller rack.  In section 5 ``{\sl The algebraic topology
of rack spaces}'' we concentrate on calculating invariants of rack
spaces.  We describe all the homotopy groups of $BX$ where $X$ is the
fundamental rack of an irreducible (non-split) link in a 3--manifold
(this is a case in which the rack completely classifies the link
[\FeRo]). This description leads to the new geometric descriptions for
the higher homotopy groups of the 3--sphere mentioned earlier.  We
also calculate $\pi_2$ of $BX$ where $X$ is the fundamental rack of a
general link in $S^3$.  We show that $BX$ is always a simple space and
we compute the homotopy type of $BX$ in the cases when $X$ is a free
rack and when $X$ is a trivial rack with $n$ elements.  We also report
on further calculations given elsewhere [\Flower, \Greene, \Wiest].

It is important to note that there are invariants of the rack space
which are not homotopy invariants, but combinatorial ones.  These
include the James--Hopf invariants (defined by the James bundles) and
the characteristic classes and associated generalised cohomology
theories constructed in [\JBun].  Although not homotopy invariants,
these all yield invariants of knots and links.  Now the homotopy type
of the rack space does not contain enough information to reconstruct
the rack: there are examples where the fundamental rack is a
classifying invariant but the homotopy type of the rack space does not
classify.  See the remarks following theorem 5.4.  However the
combinatorics of the rack space contain all the information
needed to reconstruct the rack, so in principle combinatorial
invariants should give a complete set of invariants.

This program is explored further in [\Class] where we prove that in
the classical case of links in $S^3$ the rack together with the
canonical class in $\pi_2(B\Gamma)$ determined by any diagram is a
complete invariant for the link.  This leads to computable invariants
which can effectively distinguish different links.

This paper appeared in preliminary form as part of our 1996 preprint
[\James] and many of the results were announced with outline proofs in
1993 in [\Intro].  Since this early work of ours, other authors have
investigated rack and quandle cohomology, notably J Scott Carter et al
[\Carteretal].  Rack cohomology is the cohomology of the rack space and
quandle cohomology is a quotient, see Litherland and Nelson [\LN].

\section{Basic definitions}

We give here a minimal set of definitions for $\sq$--sets.
For more detail, other definitions and examples, see [\Trunks; sections
2 and 3] and [\JBun; section 1].

\sh{The category $\tsq$}

The {\sl $n$--cube} $I^n$ is the subset $[0,1]^n$ of $\re^n$.   

A {\sl $p$--face} of $I^n$ is a subset defined by choosing $n-p$
coordinates and setting some of these equal to 0 and the rest to 1.  
In particular there are $2n$
faces of dimension $n-1$ determined by setting $x_i=\epsilon$ where
$i\in\{1,2,\ldots,n\}$ and $\ep\in \{0,1\}$.

A 0--face is called a {\sl vertex} and corresponds to a point
of the form
$(\ep_1,\ep_2,\ldots,\ep_n)$ where $\ep_i=0$ or $1$ and
$i\in\{1,2,\ldots,n\}$.  The 1--faces are called {\sl edges} and the
2--faces are called {\sl squares}.

Let $p\leq n$ and let $J$ be a $p$--face of $I^n$.  Then there is
a canonical {\sl face map} $\lambda\co I^p\to I^n$, with $\lambda(I^p)=J$,
given in coordinate
form by preserving the order of the coordinates $(x_1,\ldots,x_p)$
and inserting $n-p$ constant coordinates which are either 0 or 1.  
If $\la$ inserts only $0$'s (resp.\ only $1$'s) we call it a {\sl front}
(resp.\
{\sl back}) face map.
Notice that any face map has a unique {\sl front--back} decomposition as $\la\mu$, say,
where $\la$ is a front face map and $\mu$ is a back face map. There is also a unique 
{\sl  back--front} decomposition.
There are $2n$ face maps defined by the $(n-1)$--faces which are denoted
$\delta_i^\ep\co I^{n-1}\to I^n$, and given by:
$$\delta_i^\ep(x_1,x_2,\ldots,x_{n-1})=
(x_1,\ldots,x_{i-1},\ep,x_{i},\ldots,x_{n-1}),\qquad\epsilon\in \{0,1\}.$$

The following relations hold:
$$
\de_i^\ep\de_{j-1}^\omega
=\de_j^\omega\de_i^\ep,\quad 1\leq i<j\leq n,\quad
\ep,\omega\in \{0,1\}.\leqno{\bf \label}$$

\rk{Definition} The category $\sq$ is the category whose objects are the $n$--cubes $I^n$ for $n=0,1,\ldots$ and whose morphisms are the face maps.
\ppar

\sh{$\tsq$--Sets and their Realisations}

A {\sl $\sq$--set} is a functor $C\co\sq^{op}\to Sets$ where
$\sq^{op}$ is the opposite category of $\sq$ and $Sets$ denotes the
category of sets.

A {\sl $\sq$--map} between $\sq$--sets is a natural transformation.

We write $C_n$ for $C(I^n)$, $\lambda^*$ for $C(\lambda)$ and we
write $\d_i^\ep$ for $C(\delta_i^\ep)=(\delta_i^\ep)^*$.

The {\sl realisation\/} $\mo C\mo$ of a $\sq$--set $C$ 
is given by making the identifications
$(\lambda^*x,t)\sim (x,\lambda t)$ in the disjoint union
$\coprod_{n\ge0} C_n\times I^n$. 

We shall call 0--cells (resp.\ 1--cells, 2--cells) of $\mo C\mo$
vertices (resp.\ edges, squares) and this is consistent with the
previous use for faces of $I^n$, since $I^n$ determines a $\sq$--set
with cells corresponding to faces, whose realisation can be
identified in a natural way with $I^n$.

Notice that $\mo C\mo$ is a CW complex with one $n$--cell for each
element of $C_n$ and that each $n$--cell has a canonical
characteristic map from the $n$--cube.  However, not every CW complex
with cubical characteristic maps comes from a $\sq$--set --- even if
the cells are glued by isometries of faces.  In $\mo C\mo$, where $C$
is a $\sq$--set, cells are glued by face {\it maps}, in other words by
canonical isometries of faces.

There is also a notion of a {\sl$\sq$--space\/} namely a functor $X\co
\sq^{op}\to Top$ (where $Top$ denotes the category of topological
spaces and continuous maps) and its realisation $\mo X\mo$ given by
the same formula as above.

\rk{Notation}We shall often omit the mod signs and use the notation
$C$ both for the $\sq$--set $C$ and its realisation $\mo C\mo$.  We
shall use the full notation whenever there is any possibility of
confusion.

\proclaim{Key Example}The rack space\rm

A {\sl rack} is a set $R$ with a binary operation written $a^b$ such
that $a\mapsto a^b$ is a bijection for all $b\in R$ and such that the
{\sl rack identity}
$$a^{bc}=a^{cb^c}$$ holds for all $a,b,c\in R$.  (Here we use the
conventions for order of operations derived from exponentiation
in arithmetic.  Thus $a^{bc}$ means $(a^b)^c$ and $a^{cb^c}$ means
$a^{c(b^c)}$.)

For examples of racks see [\FeRo].

If $R$ is a rack, the {\sl rack space} is the $\sq$--set denoted $BR$ and  
defined by:

$BR_n = R^n$\quad (the $n$--fold cartesian product of $R$ with itself).
$$\d_i^0(x_1,\ldots ,x_n)=(x_1,\ldots
,x_{i-1},x_{i+1},\ldots ,x_n),$$
$$\d_i^1(x_1,\ldots ,x_n)=
((x_1)^{x_i},\ldots ,(x_{i-1})^{x_i},x_{i+1},\cdots
,x_n)
\hbox {\quad for\quad}1\leq i\leq n.$$

More geometrically, we can think of $BR$ as the $\sq$--set with one
vertex, with (oriented) edges labelled by rack elements and with
squares which can be pictured as part of a
link diagram with arcs labelled by $a$, $b$ and $a^b$ (figure 
\figkey\RSpace).

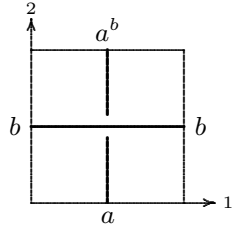
\fig{\RSpace: Diagram of a typical 2--cell of the rack space}
\beginpicture\ninepoint
\setcoordinatesystem units <.4truein,.4truein> point at 0 0
\put {$b$} [r] <-3pt,0pt> at 0 1
\put {$a$} [t] <0pt,-3pt> at 1 0
\put {$a^b$} [b] <0pt,3pt> at 1 2
\put {$b$} [l] <3pt,0pt> at 2 1
\setlinear
\plot 0 0 2 0 2 2 0 2 0 0 /
\arrow <4pt> [.2, .7] from 0 2 to 0 2.4
\arrow <4pt> [.2, .7] from 2 0 to 2.4 0
\setplotsymbol ({\ninerm .})
\plot 1 0 1 0.85 /
\plot 1 1.15 1 2 /
\plot 0 1 2 1 /
\put {$\scriptstyle1$} [l] <2pt,0pt> at 2.4 0
\put {$\scriptstyle2$} [b] <0pt,2pt> at 0 2.4
\endpicture
\endfig

The higher dimensional cubes are determined by the squares: roughly
speaking a cube is determined by its 2--skeleton, for more detail
see [\JBun; example 1.4.3].

Notice that the rack space of the rack with one element has
precisely one cube in each dimension.  This is a
description of the {\sl trivial $\sq$--set}.

\sh{Associated James complexes of a $\tsq$--set}

\rk{Projections}
An {\sl$(n+k,k)$--projection} is a function $\lambda\co I^{n+k}\to I^k$
of the form $$\lambda \co (x_1,x_2,\ldots,x_{n+k})\mapsto (x_{i_1},x_{i_2},
\ldots,x_{i_k}),\leqno{\bf \label}$$ where $1\le i_1<i_2<\ldots<i_k\le n+k$.

Let $P_k^{n+k}$ denote the set of $(n+k,k)$--projections.  Note that
$P_k^{n+k}$ is a set of size $n+k\choose k$.

Let $\lambda\in P_k^{n+k}$ and let $\mu\co I^l\to I^k,\quad l\le k$ be a face
map.  The projection $\mu^\sharp(\lambda)\in P_l^{n+l}$ and the face map
$\mu_\lambda\co I^{n+l}\to I^{n+k}$ are defined uniquely by the 
following pull-back diagram:
$$\matrix{
I^{n+l}&\buildrel\mu_\lambda\over\a&I^{n+k}\cr
\cr
{\scriptstyle\mu^\sharp(\lambda)}\downarrow&&{\scriptstyle\lambda}\downarrow\cr
\cr
I^l&\buildrel\mu\over\a&I^k\cr}
\leqno{\bf \label}$$

\rk{Definition}
Let $C$ be a $\sq$--set.  The {\sl$n^{\rm th}$ associated James
complex} of $C$, denoted $J^n(C)$ is defined as follows.  The
$k$--cells are given by   
$$J^n(C)_k=C_{n+k}\times P_k^{n+k}$$
and face maps by
$$\mu^*(x,\lambda)=(\mu_\lambda^*(x),\mu^\sharp(\lambda))$$
where $\mu\co I^l\to I^k,\quad l\le k$ is a face map. 

\rk{Notation}Let $\lambda\in P_k^{n+k}$ and $c\in C_{n+k}$ then
we shall use the notation $c_\lambda$ for the $k$--cube $(c,\lambda)\in
J^n(C)$.  When necessary, we shall use the full notation $(\la_1,\ldots,\la_n)$
for the projection $\lambda$ (given by formula 1.2) where
$\lambda _1<\lambda _2<\ldots<\lambda _n$
 and $\{\lambda _1,\ldots,\lambda _n\}=\{1,\ldots,k+n\}-
\{i_1,\ldots,i_k\}$.  In other words we index cubes of $J^n$ by the
$n$ directions (in order) which are collapsed by the defining
projection.

\sh{Picture for James complexes}

We think of $J^n(C)$ as comprising all the codimension $n$ central
subcubes of cubes of $C$.   For example a 3--cube $c$ of $C$ gives rise
to the three 2--cubes of $J^1(C)$ which are illustrated in
figure \figkey\JaComp.

\fig{\JaComp}
\beginpicture
\setcoordinatesystem units <.8mm,.8mm> point at 0 0
\ninepoint
\put {$\scriptstyle1$} [t] at 7 -0.5
\put {$\scriptstyle2$} [br] at 5 2.5
\put {$\scriptstyle3$} [r] at 0 7
\put {$c_{(2)}$} [tl] at 40 5
\put {$c_{(3)}$} [l] at 50 25
\put {$c_{(1)}$} [b] at 35 41
\putrule from 0 0 to 30 0
\putrule from 0 30 to 30 30
\putrule from 0 0 to 0 30
\putrule from 30 0 to 30 30
\putrule from 20 40 to 50 40
\putrule from 50 10 to 50 40
\plot 0 30  20 40 /
\plot 30 30  50 40 /
\plot 30 0  50 10 /
\setdashes
\putrule from 20 10  to 20 40
\putrule from 20 10  to 50 10
\plot 0 0 20 10 /
\setsolid
\linethickness=1.2pt
\putrule from 0 15 to 30 15
\putrule from 15 0 to 15 30
\putrule from 10 35 to 40 35
\putrule from 40 5 to 40 35
\putrule from 10 5 to 10 35
\putrule from 10 5 to 40 5
\putrule from 35 10 to 35 40 
\putrule from 20 25 to 50 25
\putrule from 25 5 to 25 35
\putrule from 10 20 to 40 20
\setlinear\setplotsymbol ({\tenrm .})
\plot 15 30  35 40 /
\plot 30 15  50 25 /
\plot 15 0  35 10 /
\plot 15 15  35 25 /
\plot 0 15  20 25 /
\endpicture
\endfig

In figure \JaComp\ we have used full notation for projections. 
Thus for example,
$c_{(2)}$ corresponds to the projection $(x_1,x_2,x_3)\mapsto(x_1,x_3)$,
($x_2$ being collasped).

This picture can be made more precise by considering the section
$s_\lambda\co I^k\to I^{n+k}$ of $\lambda$ given by
$$
\def\ha{{\scriptstyle{1\over2}}}
s_\lambda(x_1,x_2,\ldots,x_k)=(\ha,\ldots,\ha,x_1,\ha,\ldots,\ha,x_2,\ha,
\ldots,\ha,x_k,\ha,\ldots)
$$
where the non-constant coordinates are in places $i_1,i_2,\ldots,i_k$
and $\lambda$ is given by 1.2.

For example in the picture the image of $s_\lambda$ where 
$\lambda(x_1,x_2,x_3)\mapsto(x_1,x_3)$ is the 2--cube labelled $c_{(2)}$.

Now the commuting diagram (1.3) which defines the face maps implies
that the $s_\lambda$'s are compatible with faces and hence they fit 
together to define a map 
$$p_n\co\mo J^n(C)\mo\to\mo C\mo\quad\hbox{given by}\quad 
p_n[c_\lambda,t]=[c,s_\lambda(t)].$$  
In the next section we will see that $p_n$ is a mock bundle projection.

\section{Mock bundles and transversality}

In this section we define mock bundles over smooth CW complexes, which
include $\sq$--sets, and prove a transversality theorem for
$\sq$--sets.  This material is similar to material in [\BRS; Chapters
2 and 7].  However [\BRS] is set entirely in the PL category and deals
only with transversality with respect to a transverse CW complex.
Here we shall need to extend the work to the smooth category and prove
transversality with respect to a $\sq$--set (which is not quite a
transverse CW complex).  However many of the proofs are similar to
proofs in [\BRS] and therefore we omit details when appropriate.
Similar material can also be found in [\Fenn] and [\Pastor].

The main technicalities concern manifolds with corners, which is where
we start.  We shall use {\it smooth} to mean $C^\infty$.  

\rk{Definition}{\sl Manifold with corners}

For background material on smooth manifolds with corners, see Cerf or
Douady [\Cerf, \Dou].  In particular these references contain an
appropriate version of the tubular neighbourhood theorem.  There is a
uniqueness theorem for these tubular neighbourhoods.  The proof can be
obtained by adapting the usual uniqueness proof.

A {\sl smooth $n$--manifold with corners} $M$ is a space modelled on
$\E^n$ where $\E=[0,\infty)$.  In other words $M$ is equipped with a
maximal atlas of charts from open subsets of $\E^n$ such that overlap
maps are smooth.  There is a natural stratification for such a
manifold.  Define the {\sl index} of a point $p$ to be the minimum $q$
such that a neighbourhood of $p$ in $M$ is diffeomorphic to an open
subset of $\E^q\times\re^{n-q}$.  The {\sl stratum of index $q$}
denoted $M_{(q)}$ comprises all points of index $q$.  Note that the
dimension of $M_{(q)}$ is $n-q$ and that $M_{(0)}$ is the interior of
$M$, $M_{(1)}$ is an open codimension 1 subset of $\d M$ and in
general $M_{(i)}$ is an $(n-i)$--manifold lying in the closure of
${M_{(i-1)}}$.

\rk{Examples}$I^n$ and $\E^n$ are manifolds with corners.  
If $M$ and $Q$ are manifolds with corners then so is $M\times Q$.

\rk{Definition}{\sl Maps of manifolds with corners}

Let $M$ and $Q$ be manifolds with corners then a {\sl stratified map}
is a smooth map $f\co M\to Q$ such that $f(M_{(q)})\subset Q_{(q)}$
for each $q$.  This is the analogue for manifolds with corners of
a proper map for manifolds with boundary.

An {\sl embedding of manifolds with corners} is a smooth embedding
$i\co M\subset Q$ such that the pair is locally like the inclusion
of $\E^p\times \E^t$ in $\E^p\times \re^t\times \re^s$ where 
$p+t={\rm dim}(M)$ and $p+t+s={\rm dim}(Q)$.  Thus an embedding
of manifolds with corners which is proper (ie takes boundary to
boundary and interior to interior) must be
a stratified embedding.  However in general an embedding of
manifolds with corners allows corners on
$M$ not at corners of $Q$, see the examples below.

A {\sl face map of manifolds with corners} 
$\lambda\co M\to Q$ is a smooth embedding which is a proper
map of topological spaces (preimage of compact is compact) such that
$\lambda(M_{(q)})\subset Q_{(t+q)}$ 
for each $q$ where $t={\rm dim}(Q)
-{\rm dim}(M)$.  Thus a face map is a diffeomorphism of $M$ onto a union
of components of strata of $Q$.

\rk{Examples}The map $s_\lambda:I^k\to I^{n+k}$ of $\lambda$ given by
$$
s_\lambda(x_1,x_2,\ldots,x_k)=(\ha,\ldots,\ha,x_1,\ha,\ldots,\ha,x_2,\ha,
\ldots,\ha,x_k,\ha,\ldots)
$$
(see the end of section 1) is a stratified embedding.  

If $M$ is a manifold with
corners and $Q$ is a smooth manifold without boundary then the projection
$M\times Q\to M$ is a stratified map.

The inclusions $I^n\subset\E^n\subset\re^n$ are embeddings of
manifolds with corners.

A face map (in the sense of section 1) $\lambda\co I^k\to I^n$ is a
face map of manifolds with corners.

\rk{Definition}{\sl Smooth CW complex}

We refer to [\PLTop; pages 13--14] for the
definition and basic properties of convex linear cells in $\re^n$
(which we shall abbreviate to convex cells).  Convex cells
have well-defined faces and if $e'$ is a face of $e$ we write
$e'<e$.  Convex cells are smooth manifolds with corners.
A smooth face map between convex cells (in the sense of
manifolds with corners) is the same as a diffeomorphism 
onto a face. 

A {\sl smooth CW complex} is a collection of convex cells glued
by smooth face maps.  More precisely it comprises a CW complex
$W$ and for each cell $c\in W$ a preferred characteristic map
$\chi_c\co e_c\to W$ where $e_c$ is a convex cell such that
for each face $e'<e_c$ there is a cell $d\in W$ with preferred
characteristic map $\chi_d\co e_d\to W$ and a diffeomorphism $\mu\co e_d
\to e'$ such that the diagram commutes:
$$
\matrix{
e_d&\buildrel\chi_d \over \a&W\cr
&&\cr
\downarrow\!{\scriptstyle\mu}&&\uparrow\!{\scriptstyle\chi_c}\cr
&&\cr
e'&\buildrel\rm inc.\over \a&e_c\cr}
$$
\rk{Examples}A (realised) $\sq$--set or $\Delta$--set, a
simplicial complex or a convex linear cell complex are
all examples of smooth CW complexes.

We say that a smooth CW complex $W$ gives rise to a {\sl smooth
decomposition} of a smooth manifold $M$ (possibly with corners) if
there is a homeomorphism $h\co W\to M$ such that $h\circ\chi_c\co
e_c\to M$ is an embedding of manifolds with corners for each cell
$c\in W$.  Usually we identify $W$ and $M$ via $h$ in this situation
and say that $W$ is a smooth decomposition of $M$.

\proc{Definition}\rm {\sl Smooth mock bundle}

Let $W$ be a smooth CW complex.   A {\sl mock bundle $\xi$ over 
$W$} of codimension $q$ (denoted $\xi^q/W$) comprises
a {\sl total space} $E_\xi$ and a {\sl projection} 
$p_\xi\co E_\xi\to W$ with the following property.%
\fnote{Note that the notation used here for dimension of
a mock bundle, namely that $q$ is {\it codimension}, is the negative of that used
in [\BRS] where $\xi^q/W$ meant a mock bundle of {\it fibre dimension}
$q$ ie codimension $-q$.  The notation used here is consistent with the
usual convention for cohomology.}

Let $c$ be an $n$--cell of $W$ with characteristic map 
$\chi_c\co e_c\to  W$, then there is a smooth manifold with 
corners $B_c$ of dimension $n-q$ called the {\sl block} over $c$
and a stratified map $p_c\co B_c\to e_c$ and a map $b_c\co B_c\to E_\xi$
such that the following diagram is a pull-back:
$$
\matrix{
B_c&\buildrel b_c\over \a&E_\xi\cr
&&\cr
\downarrow\!{\scriptstyle p_c}&&\downarrow\!{\scriptstyle p_\xi}\cr
&&\cr
e_c&\buildrel \chi_c\over \a&W\cr}
$$

\proc{Amalgamation lemma}Suppose that $\xi^q/W$ is a smooth
mock bundle and that $W$ is a smooth decomposition of a smooth manifold
with corners $M^m$.  Then $E_\xi$ can be given the structure of a
smooth manifold with corners of dimension $m-q$ such that $p_\xi\co
E_\xi\to M$ is $\varepsilon$--homotopic through mock bundle projections to a
stratified map.

\prf By standard embedding theorems we may assume that each block 
$B_c$ of $\xi$ is a stratified submanifold of $e_c\times\re^N$ for
some $N$ and that these embeddings fit together to give an embedding
of $E_\xi$ in $M\times\re^N$ such that $p_\xi$ is the restriction of
projection on the first coordinate.  We shall isotope this embedding
so that $E_\xi$ becomes a smooth submanifold of $M\times\re^N$.  We
work inductively over the skeleta of $W$.  Assume inductively that
this isotopy has already been carried out over a neighbourhood of the
$(i-1)$--skeleton of $W$.

Now fix attention on the interior of a particular block $B^\circ_c$
which is the subset of $E_\xi$ lying over the interior $c^\circ$ of an
$i$--cell of $W$.  The standard tubular neighbourhood theorem applied
to the submanifold $c^\circ$ of $\mo W\mo=M$ yields a tubular
neighbourhood $\la$ of $c^\circ$ in $\mo W\mo=M$ formed by tubular
neighbourhoods in each of the incident cells.  Using the tubular
neighbourhood theorem for manifolds with corners we can construct a
(non-smooth) tubular neighbourhood $\mu$ of $B^\circ_c$ in $E_\xi$
formed by (smooth) tubular neighbourhoods in each of the incident
blocks, which extends to a tubular neighbourhood $\mu^+$ on
$c^\circ\times\re^N$ in $M\times \re^N$.  By inductively applying
uniqueness we can deform $\mu^+$ near $B^\circ_c$ to $\la\times\re^N$
by an $\varepsilon$--isotopy.  This carries $E_\xi$ to a smooth
submanifold of $M\times\re^N$ so that projection on $M$ is a
stratified map and moreover the isotopy determines a homotopy through
mock bundle projections of the projection on $\mo W\mo$.  By induction
$E_\xi$ is already smooth near $\d B$ and we can keep a neighbourhood
of $\d B$ fixed through the isotopy.  Do this for each $i$--cell of $W$
to complete the induction step.
\qed

\rk{Definition}{\sl Maps of smooth CW complexes}

A {\sl linear projection} of convex cells is a surjective map $f\co
e_1\to e_2$, where $e_1, e_2$ are convex cells, which is the
restriction of an affine map, and such that $f(e')<e_2$ for each face
$e'<e_1$.  Examples include simplicial maps of one simplex onto
another, projections $I^{n+k}\to I^n$ and projections of the form
$d\times e\to e$ where $d,e$ are convex cells.

A {\sl smooth projection} of convex cells is a linear projection
composed with a diffeomorphism of $e_1$.

A {\sl smooth map $f\co W\to Z$ of smooth CW complexes} is
a map such that for each cell $c\in W$ there is a cell
$d\in Z$ and a smooth projection $\phi\co e_c\to e_d$
such that the diagram commutes:
$$
\matrix{
e_c&\buildrel\phi\over \a&e_d\cr
&&\cr
\downarrow\!{\scriptstyle\chi_c}&&\downarrow\!{\scriptstyle\chi_d}\cr
&&\cr
W&\buildrel f\over \a&Z\cr}
$$

Examples include $\sq$--maps, $\Delta$--maps, simplicial maps
and projections $W\times Z\to W$, where $W,Z$ are any two smooth
CW complexes.

The following lemma follows from definitions:

\proc{Pull-back lemma}Let $\xi^q/Z$ be a smooth mock bundle
and $f\co W\to Z$ a smooth map of smooth CW complexes, then
the pull-back $f^*(E_\xi)\to W$ is a mock bundle of the same
codimension (denoted $f^*(\xi^q)/W$). \qed

\sh{Properties of mock bundles}

We shall summarise properties of mock bundles.  The details of all
the results stated here are analogous to the similar results for
PL mock bundles over cell complexes proved in [\BRS].  

The set of cobordism classes of mock bundles with
base a smooth CW complex forms an abelian group (under disjoint
union of total spaces) and there is a relative group 
(the total space is empty over
the subcomplex).  This all fits together with the pull-back
construction to define a cohomology 
theory which can be identified
with smooth cobordism (classified by the Thom spectrum $MO$).

If the base is a manifold then the amalgamation 
lemma defines a map from $q$--cobordism to $(n-q)$--bordism.
This map is the Poincar\'e duality isomorphism.

Given two mock bundles $\xi^q/W$ and $\eta^r/W$ we can define the mock bundle
$(\xi\cupprod\eta)^{q+r}/W$ in various equivalent ways analogous to the 
Whitney sum of bundles.  We can pull one bundle back over the total space
of the other and then compose.  This is equivalent to making the 
projection of the first bundle transverse to the projection of
the second and then pulling back.  We can take the
external product $\xi\times\eta/W\times W$ and restrict to the diagonal.
These equivalent constructions define the cup product in cobordism.
If the base is a manifold then pull-back (or transversality)
defines the cap product which coincides under Poincar\'e duality
with the cup product.

Mock bundles can be generalised and extended in a number of ways.  The
simplest is to use orientation.  If each block (and cell) is oriented
in a compatible way (cf [\BRS; page 82]) then the resulting theory of oriented
mock bundles defines oriented cobordism (classified by $MSO$). More
generally we can consider restrictions on the stable normal bundle
of blocks and this yields the corresponding cobordism theory.
A particular example of relevance here is the case when blocks are
stably framed manifolds; in this case the resulting theory is stable 
cobordism
classified by the sphere spectrum $\Bbb S$.  By considering manifolds
with singularities, the resulting theory can be further generalised
and such a mock bundle theory can represent the cohomology theory corresponding
to an arbitrary spectrum [\BRS; chapter 7].  The corresponding homology
theory is represented by the bordism theory given by using manifolds
with the same allowed singularities.  Coefficients and sheaves of
coefficients can also be defined geometrically (see [\BRS; chapters
3 and 6]).

\rk{Key example}{\sl James bundles}

At the end of the last section we defined a projection
$$p_n:\mo J^n(C)\mo\to\mo C\mo$$  
Where $J^n(C)$ is the $n$-th associated James complex of the $\sq$--set $C$.

Now  if we choose
a particular $(n+k)$--cell $\sigma$ of $C$ then the pull back of $p_n$
over $I^{n+k}$ (by the characteristic map for $\sigma$) is
a $k$--manifold (in fact it is the $n+k\choose k$ copies of 
$I^k$ corresponding to the elements of $P^{n+k}_k$).  Therefore
$p_n$ is the projection of a mock bundle of codimension
$n$, which we shall call the {\sl$n$--th James bundle} of $C$
denoted $\zeta^n(C)$.

\proc{Definition}Embedded mock bundle\rm

Let $W$ be a smooth CW complex and $Q$ a smooth manifold without boundary.
An {\sl embedded mock bundle} in $W\times Q$ is defined to a
mock bundle $\xi/W$ with an embedding $E_\xi\subset W\times Q$
such that $p_\xi$ is the restriction of the projection $W\times Q
\to W$ and such that for each cell $c\in W$ the induced embedding 
$B_c\subset e_c\times Q$ is a stratified embedding. The proof of the
amalgamation lemma (with $\re^N$ replaced by $Q$) implies that if $W$
is a smooth decomposition of a manifold then
$E_\xi$ can be $\varepsilon$--isotoped to a smooth submanifold of $W\times Q$
so that projection on $W$ is still a mock bundle projection.

We shall be particularly concerned with the case when $Q$ is $\re^t$
for some $t$ and each block is framed in $e_c\times\re^t$.  The theory
defined by mock bundles of this type is unstable cohomotopy.  See in
particular [\JBun; 3.6 and 5.1].

\rk{Key example}{\sl Embedding the James bundles in $\mo C\mo\times\re$}

Let $C$ be a $\sq$--set.  The James bundles can be embedded
in $\mo C\mo\times\re$.  This is done by ordering
the cubes of $J^n(C)$ over a particular cube of $C$ and 
lifting in that order.  Recall that the $k$--cubes of $J^n(C)$
lying over a $(k+n)$--cube are indexed by projections 
$\lambda=(i_1,\ldots,i_n)\in P^{n+k}_n$.  These may
be ordered lexicographically. The lexicographic order is compatible 
with face maps and can be used to define the required embedding
by induction on dimension of cells of $C$ as follows.  

Suppose inductively that the embedding has been defined 
over cells of $C$ of dimension $\le k+n-1$.  

Consider a $(k+n)$--cube $c\in C$ with characteristic map
$\chi_c\co I^{k+n}\to \mo C\mo$.  Pulling the embedding
back (where it is already defined) over $\chi_c$ gives an
embedding of $\zeta^n(\d I^{k+n})$ in $\d I^{k+n}\times\re$.
Now embed the centres of the
$k$--cubes of $J^n(I^{k+n})$ at $(\ha,\ldots,\ha)\times r_\lambda$
where $r_\lambda$ are real numbers for $\lambda\in P^{n+k}_n$ chosen
to increase strictly corresponding to the lexicographic order
on $P^{n+k}_n$.  

Now embed each $k$--cube of $J^n(I^{k+n})$ as the cone
on its (already embedded) boundary.  Finally smooth the resulting
embedding and push it forwards to
$\mo C\mo\times\re$ using $\chi_c\times {\rm id}$.

Precise smooth formul\ae\ for this embedding can be found in 
[\JBun; section 3] using a bump function.

The embedding is in fact framed.  This can be seen as follows.
Each $k$--cube $p_n(c,\la)$, where $\lambda=(i_1,\ldots,i_n)$, of
$J^n(I^{k+n})$ is framed in $I^{k+n}$ by the $n$ vectors parallel
to directions $i_1,\ldots,i_n$.  These lift to parallel vectors
in $I^{n+k}\times\re$ and the framing is completed by the vector
parallel to the positive $\re$ direction (vertically up).  This
framing is compatible with faces and defines a framing of $\zeta^n(C)$
in $\mo C\mo\times\re$.  The formul\ae\ in [\JBun; section 3] 
also give formul\ae\ for the framing.

For the special case $n=1$ the map of $\zeta^n(C)$ to $\re$ can be
simply described: the centre of $c_{(k)}$ is mapped to $k$. This
determines a map, linear on simplexes of $Sd_\Delta \zeta^1(C)$, to
$\re$. It follows that the centre of $\d_i^\epsilon c_{(k)}$ is mapped
to $k$ if $i\ge k$ and to $k-1$ if $i<k$.

In figure \JaComp\ we illustrated $J^1(C)$  for a 3--cube
$c\in C$.  The embeddings in $\mo C\mo\times\re$ (before smoothing) above
each of the three 2--cubes are illustrated in figure \figkey\JaEmb.

\fig{\JaEmb}
\beginpicture
\setcoordinatesystem units < 0.7cm, 0.7cm>
\unitlength= 0.7cm
\linethickness=1pt
\setplotsymbol ({\makebox(0,0)[l]{\tencirc\symbol{'160}}})
\setshadesymbol ({\thinlinefont .})
\setlinear
%
%
\linethickness= 0.500pt
\setplotsymbol ({\thinlinefont .})
\plot  0.921 20.034  2.826 20.415 /
\plot  2.826 20.415  4.731 20.034 /
\plot  4.731 20.034  4.921 21.654 /
\plot  4.921 21.654  5.683 22.892 /
\plot  5.683 22.892  3.778 23.273 /
\plot  3.778 23.273  1.873 22.892 /
\plot  1.873 22.892  1.111 21.654 /
\plot  1.111 21.654  0.921 20.034 /
\putrule from  0.921 20.034 to  0.921 20.034
%
%
\linethickness= 0.500pt
\setplotsymbol ({\thinlinefont .})
\plot  3.112 22.415  1.111 21.654 /
%
%
\linethickness= 0.500pt
\setplotsymbol ({\thinlinefont .})
\plot  3.112 22.415  1.873 22.892 /
%
%
\linethickness= 0.500pt
\setplotsymbol ({\thinlinefont .})
\plot  3.112 22.415  3.778 23.273 /
%
%
\linethickness= 0.500pt
\setplotsymbol ({\thinlinefont .})
\plot  3.112 22.415  5.683 22.892 /
%
%
\linethickness= 0.500pt
\setplotsymbol ({\thinlinefont .})
\plot  3.112 22.415  4.921 21.654 /
%
%
\linethickness= 0.500pt
\setplotsymbol ({\thinlinefont .})
\plot  3.112 22.415  2.826 20.415 /
\plot  3.112 22.415  4.731 20.034 /   
\plot  3.112 22.415  0.921 20.034 /   
%
%
\linethickness= 0.500pt
\setplotsymbol ({\thinlinefont .})
\putrule from  0.921 20.034 to  4.731 20.034
\plot  4.731 20.034  5.683 22.892 /
\putrule from  5.683 22.892 to  3.683 22.892
%
%
\linethickness= 0.500pt
\setplotsymbol ({\thinlinefont .})
\plot  0.921 20.034  1.492 21.654 /
%
%
\linethickness= 0.500pt
\setplotsymbol ({\thinlinefont .})
\plot  1.873 22.892  1.587 21.939 /
%
%
\linethickness= 0.500pt
\setplotsymbol ({\thinlinefont .})
\putrule from  1.873 22.892 to  3.302 22.892
%
%
\linethickness= 0.500pt
\setplotsymbol ({\thinlinefont .})
\putrule from 12.351 20.034 to 16.161 20.034
\plot 16.161 20.034 17.113 22.892 /
\putrule from 17.113 22.892 to 13.208 22.892
\plot 13.208 22.892 12.351 20.034 /
%
%
\linethickness= 0.500pt
\setplotsymbol ({\thinlinefont .})
\plot  6.636 20.034  8.636 20.415 /
\plot  8.636 20.415 10.446 20.034 /
%
%
\linethickness= 0.500pt
\setplotsymbol ({\thinlinefont .})
\plot  7.588 22.892  9.493 23.273 /
%
%
\linethickness= 0.500pt
\setplotsymbol ({\thinlinefont .})
\plot  9.493 23.273 11.398 22.892 /
%
%
\linethickness= 0.500pt
\setplotsymbol ({\thinlinefont .})
\plot  9.493 23.273  8.636 20.415 /
%
%
\linethickness= 0.500pt
\setplotsymbol ({\thinlinefont .})
\putrule from  9.684 22.892 to 11.398 22.892
\plot 11.398 22.892 10.446 20.034 /
\putrule from 10.446 20.034 to  6.636 20.034
\plot  6.636 20.034  7.588 22.892 /
\putrule from  7.588 22.892 to  9.112 22.892
%
%
\def\SetFigFont{\ninepoint}%
\put{\SetFigFont $c(2)$ } [lB] <0pt, 10pt> at  8.446 18.605
%
%
\put{\SetFigFont $c(3)$} [lB] <0pt, 10pt> at  2.445 18.605
%
%
\put{\SetFigFont $c(1)$} [lB] <0pt, 10pt> at 14.161 18.701
%
%
\put{\SetFigFont$\scriptstyle 1$} [lB] <0pt,2pt> at  1.302 19.558
%
%
\put{\SetFigFont$\scriptstyle  2$} [lB] <0pt, 5pt> at  1.302 20.701
%
%
\put{\SetFigFont$\scriptstyle  2$} [lB] <0pt, 0pt> at 12.732 19.653
%
%
\put{\SetFigFont$\scriptstyle  3$} [lB] <-4pt, 0pt> at 12.827 20.701
%
%
\put{\SetFigFont$\scriptstyle  1$} [lB] at  7.017 19.653
%
%
\put{\SetFigFont$\scriptstyle  3$} [lB] <-2pt, 0pt> at  7.017 20.606
\linethickness=0pt
\putrectangle corners at  0.921 23.273 and 17.113 18.605
\endpicture
\endfig

We finish this section with a discussion
of transversality with respect to a $\sq$--set, which will be important
for the main classification theorems of section 4.  A similar
treatment can be given for any smooth CW complex, but we shall
not need this in this paper.
\rk{Transversality}

\rk{Definition}{\sl Transverse map to a $\sq$--set}

Let $M$ be a smooth manifold (possibly with boundary) of dimension $m$
and $C$ a $\sq$--set.  Let $c$ be an
$n$--cell of $C$ with characteristic map $\chi_c\co I^n\to \mo C\mo$
denote $\chi_c(I^n-\d I^n)$ by $c^\circ$ (the interior of $c$) and
$\chi_c(\ha,\ldots,\ha)$ by $\hat c$ (the centre of $c$).

Let $f\co M\to\mo C\mo$ be a map.  Let $M_c$ denote the closure
of $f\inv(c^\circ)$ and $N_c$ the closure of ${f\inv(\hat c)}$.

We say $f$ is {\sl transverse} to $c$ if  $M_c$ 
is a smooth $m$--manifold with corners embedded
(as a manifold with corners) in $M$ and equipped with a
diffeomorphism 
$\iota_c\co N_c\times I^n\a M_c$
such that the diagram commutes:
$$
\matrix{
N_c\times I^n&\buildrel \iota_c\over \a&M_c\cr
&&\cr
\downarrow\!{\scriptstyle p_2}&&\downarrow\!{\scriptstyle f}\cr
&&\cr
I^n&\buildrel \chi_c\over \a&\mo C\mo\cr}
$$
where $p_2$ denotes projection on the second coordinate.   Thus $N_c$ 
is a framed submanifold (with corners) of codimension $n$
framed by copies of $I^n$ on each of which $f$ is the characteristic
map for $c$.

A map $f\co M\to\mo C\mo$ is {\sl transverse} if it is transverse to
each cell $c\in C$ and the framings are compatible with face maps in
the following sense.  Let $\lambda\co I^q\to I^n$ be a face map and 
$d=\lambda^*(c)$.  Then there is a face map (of manifolds with corners)
$\la^*_c\co N_c\to N_d$ and 
the following diagram commutes:
$$
\matrix{
N_c\times I^q&\buildrel {\rm id}\times\lambda\over \a&N_c\times I^n\cr
&&\cr
\downarrow\!{\scriptstyle \la^*_c\times{\rm id}}&&\downarrow\!{\scriptstyle i_c}\cr
&&\cr
N_d\times I^q&\buildrel i_d\over \a&M\cr}
$$
where $i_c={\rm inc.}\circ\iota_c$, using the notation established above. 
The last condition can be summarised by saying that $M$ is the realisation
of a $\sq$--space and $f$ is the realisation of a $\sq$--map.

To see this, define a $\sq$--space by $X_n=\coprod_{c\in C_n}
N_c$ and $\lambda^*=\coprod\la_c^*$ 
then the diffeomorphisms $\iota_c$ define a homeomorphism $\iota\co
\mo X\mo\to M$
and if we identify $M$ with $\mo X\mo$ via $\iota$ then the
commuting diagrams above imply that
$f$ is the realisation of a $\sq$--map $\mo X\mo\to\mo C\mo$.

\rk{Remark}The framing
compatibility condition in the definition of a transverse 
map is unnecessary.  If
$f\co M\to\mo C\mo$ is transverse to each cell of $C$ then the
framings can in fact be changed to become compatible (without
altering $f$).  However, the
full definition is the one that we shall need in practice.

\rk{Ellucidation}To help the reader understand the (somewhat complicated)
concept of a transverse map to a $\sq$--set we shall describe 
transversality for maps of closed surfaces and 3--manifolds (possibly
with boundary) into $C$.

A transverse map of a closed surface $\Sigma$ into a $\sq$--set C
meets only the 2--skeleton of $C$.
The pull-backs of the squares of $C$ are a number of disjoint little
squares in $\Sigma$ (each of which can be identified with the standard
square $I^2$ and maps by a characteristic map to a 2--cell of $C$).  
The pull-backs of the 1--cells are a number of bicollared
1--manifolds which are either bicollared closed curves or are
attached to edges of the little squares (and each edge of each
square is used in this way).  Each bicollar line can identified with
$I^1$ and is mapped by $f$ to a 1--cell of $C$ (by the characteristic map
for that cell).  Finally each component of $\Sigma -\{$little squares
and collared 1--manifolds$\}$ is mapped to a
vertex of $C$.

Thus we can think of the transverse map as defining a 
thickened diagram of self-transverse
curves (with the squares at the double points).  This is illustrated
in figure \figkey\Rings.  We shall explore the connection between
transverse maps and diagrams in the next section.

\fig{\Rings}
\beginpicture
\setcoordinatesystem units < 0.6cm, 0.6cm>
\unitlength= 1.000cm
\linethickness=1pt
\setplotsymbol ({\makebox(0,0)[l]{\tencirc\symbol{'160}}})
\setshadesymbol ({\thinlinefont .})
\setlinear
%
%
\linethickness= 0.500pt
\setplotsymbol ({\thinlinefont .})
\circulararc 233.924 degrees from  4.921 25.337 center at  3.936 23.400
%
%
\linethickness= 0.500pt
\setplotsymbol ({\thinlinefont .})
\circulararc 234.353 degrees from  4.921 21.463 center at  5.916 23.400
%
%
\linethickness=1pt
\setplotsymbol ({\makebox(0,0)[l]{\tencirc\symbol{'160}}})
\circulararc 221.429 degrees from  5.239 21.780 center at  5.839 23.368
%
%
\linethickness=1pt
\setplotsymbol ({\makebox(0,0)[l]{\tencirc\symbol{'160}}})
\circulararc 220.110 degrees from  4.604 24.987 center at  4.008 23.406
%
%
\linethickness= 0.500pt
\setplotsymbol ({\thinlinefont .})
\circulararc 203.206 degrees from  4.318 24.701 center at  4.025 23.438
%
%
\linethickness= 0.500pt
\setplotsymbol ({\thinlinefont .})
\circulararc 87.154 degrees from  4.318 24.701 center at  5.621 23.398
%
%
\linethickness=1pt
\setplotsymbol ({\makebox(0,0)[l]{\tencirc\symbol{'160}}})
\circulararc 99.776 degrees from  4.635 24.352 center at  5.390 23.360
%
%
\linethickness= 0.500pt
\setplotsymbol ({\thinlinefont .})
\circulararc 101.111 degrees from  4.985 24.003 center at  5.431 23.332
%
%
\linethickness= 0.500pt
\setplotsymbol ({\thinlinefont .})
\circulararc 103.130 degrees from  4.826 22.765 center at  4.414 23.447
%
%
\linethickness=1pt
\setplotsymbol ({\makebox(0,0)[l]{\tencirc\symbol{'160}}})
\circulararc 111.860 degrees from  5.143 22.447 center at  4.595 23.395
%
%
\linethickness= 0.500pt
\setplotsymbol ({\thinlinefont .})
\circulararc 199.624 degrees from  5.588 22.066 center at  5.824 23.334
%
%
\linethickness= 0.500pt
\setplotsymbol ({\thinlinefont .})
\circulararc 81.270 degrees from  5.556 22.035 center at  4.090 23.357
%
%
\linethickness= 0.500pt
\setplotsymbol ({\thinlinefont .})
\putrule from  1.841 23.368 to  2.635 23.368
%
%
\linethickness= 0.500pt
\setplotsymbol ({\thinlinefont .})
\putrule from  3.842 23.368 to  4.540 23.368
%
%
\linethickness= 0.500pt
\setplotsymbol ({\thinlinefont .})
\putrule from  5.334 23.400 to  5.969 23.400
%
%
\linethickness= 0.500pt
\setplotsymbol ({\thinlinefont .})
\putrule from  7.239 23.368 to  8.033 23.368
%
%
\linethickness= 0.500pt
\setplotsymbol ({\thinlinefont .})
\plot  2.477 24.955  2.985 24.384 /
%
%
\linethickness= 0.500pt
\setplotsymbol ({\thinlinefont .})
\plot  3.207 22.257  2.794 21.622 /
%
%
\linethickness= 0.500pt
\setplotsymbol ({\thinlinefont .})
\plot  6.477 22.162  6.858 21.495 /
%
%
\linethickness= 0.500pt
\setplotsymbol ({\thinlinefont .})
\plot  6.826 25.305  6.413 24.543 /
%
%
\linethickness= 0.500pt
\setplotsymbol ({\thinlinefont .})
\setdots < 0.0953cm>
\plot  4.572 21.812  5.175 22.447 /
%
%
\linethickness= 0.500pt
\setplotsymbol ({\thinlinefont .})
\plot  4.604 24.384  5.271 24.955 /
%
%
\linethickness= 0.500pt
\setplotsymbol ({\thinlinefont .})
\plot  5.271 24.257  4.572 24.987 /
%
%
\linethickness= 0.500pt
\setplotsymbol ({\thinlinefont .})
\plot  4.540 22.447  5.239 21.749 /
%
%
\linethickness= 0.500pt
\setplotsymbol ({\thinlinefont .})
\setsolid
\plot  4.985 24.003  5.620 24.606 /
\plot  5.620 24.606  4.921 25.305 /
\plot  4.921 25.305  4.318 24.701 /
\plot  4.318 24.701  4.985 24.003 /
\putrule from  4.985 24.003 to  4.985 24.003
%
%
\linethickness= 0.500pt
\setplotsymbol ({\thinlinefont .})
\plot  4.921 21.463  5.556 22.066 /
\plot  5.556 22.066  4.858 22.765 /
\plot  4.858 22.765  4.255 22.162 /
\plot  4.255 22.162  4.921 21.463 /
\putrule from  4.921 21.463 to  4.921 21.463
\linethickness=0pt
\putrectangle corners at  1.715 25.622 and  8.096 21.209
\endpicture
\endfig

A transverse map of a closed 3--manifold
into a $\sq$--set $C$ meets only the 3--skeleton of $C$.
The pullback of the 3--cubes are a number of little cubes each 
of which can be identified with $I^3$ and which map onto 3--cubes of 
$C$ by characteristic maps.  The pullback of squares
are framed 1--manifolds  (framed by a copy of the standard
square) and such that each such square is mapped onto a square of $C$
by a characteristic map.  These framed 1--manifolds are attached to
the square faces of the little cubes at their boundaries. 
The pullback of edges are framed sheets (framed by copies of $I^1$
and mapped by characteristic maps to edges of $C$).  The edges of
the sheets are attached to edges of the framed 1--manifold (ie
along 1--manifold $\times$ edge of square) and the two framings
are required to be the same here (this is the framing compatibility
condition in this case). The remainder of $M$ is then
a 3--manifold (with corners) each component of which is mapped to
a vertex of $C$.  A general view near one of the little cubes is
illustrated in figure \figkey\Cube.

\fig{\Cube}
\beginpicture
\setcoordinatesystem units < 0.7cm, 0.7cm>
\unitlength= 1.000cm
\linethickness=1pt
\setplotsymbol ({\makebox(0,0)[l]{\tencirc\symbol{'160}}})
\setshadesymbol ({\thinlinefont .})
\setlinear
%
%
\linethickness=1pt
\setplotsymbol ({\makebox(0,0)[l]{\tencirc\symbol{'160}}})
\putrule from  2.826 23.368 to  2.826 22.701
%
%
\linethickness=1pt
\setplotsymbol ({\makebox(0,0)[l]{\tencirc\symbol{'160}}})
\plot  2.826 23.368  3.493 23.749 /
%
%
\linethickness=1pt
\setplotsymbol ({\makebox(0,0)[l]{\tencirc\symbol{'160}}})
\plot  2.826 23.368  2.350 23.654 /
%
%
\linethickness=1pt
\setplotsymbol ({\makebox(0,0)[l]{\tencirc\symbol{'160}}})
\putrule from  2.826 23.368 to  2.826 26.003
%
%
\linethickness=1pt
\setplotsymbol ({\makebox(0,0)[l]{\tencirc\symbol{'160}}})
\plot  2.350 23.654  0.127 22.447 /
%
%
\linethickness=1pt
\setplotsymbol ({\makebox(0,0)[l]{\tencirc\symbol{'160}}})
\plot  2.826 23.368  0.159 21.939 /
%
%
\linethickness=1pt
\setplotsymbol ({\makebox(0,0)[l]{\tencirc\symbol{'160}}})
\plot  2.826 22.701  0.191 21.273 /
%
%
\linethickness=1pt
\setplotsymbol ({\makebox(0,0)[l]{\tencirc\symbol{'160}}})
\plot  2.826 22.701  5.366 21.018 /
%
%
\linethickness=1pt
\setplotsymbol ({\makebox(0,0)[l]{\tencirc\symbol{'160}}})
\plot  2.826 23.368  5.493 21.749 /
%
%
\linethickness=1pt
\setplotsymbol ({\makebox(0,0)[l]{\tencirc\symbol{'160}}})
\putrule from  2.350 23.654 to  2.350 26.035
%
%
\linethickness=1pt
\setplotsymbol ({\makebox(0,0)[l]{\tencirc\symbol{'160}}})
\putrule from  3.493 23.654 to  3.493 25.971
%
%
\linethickness= 0.500pt
\setplotsymbol ({\thinlinefont .})
\setdashes < 0.1270cm>
\plot  2.350 23.654  1.873 23.940 /
%
%
\linethickness= 0.500pt
\setplotsymbol ({\thinlinefont .})
\plot  1.873 23.940  1.873 23.178 /
%
%
\linethickness= 0.500pt
\setplotsymbol ({\thinlinefont .})
\plot  1.873 23.178  2.826 22.701 /
%
%
\linethickness= 0.500pt
\setplotsymbol ({\thinlinefont .})
\plot  2.350 23.654  2.350 22.987 /
%
%
\linethickness= 0.500pt
\setplotsymbol ({\thinlinefont .})
\plot  2.826 22.701  3.493 23.082 /
%
%
\linethickness= 0.500pt
\setplotsymbol ({\thinlinefont .})
\plot  3.493 23.685  3.493 23.082 /
%
%
\linethickness= 0.500pt
\setplotsymbol ({\thinlinefont .})
\plot  2.826 22.670  2.826 22.098 /
%
%
\linethickness= 0.500pt
\setplotsymbol ({\thinlinefont .})
\plot  2.826 22.098  1.873 22.543 /
%
%
\linethickness= 0.500pt
\setplotsymbol ({\thinlinefont .})
\plot  1.873 23.178  1.873 22.543 /
%
%
\linethickness= 0.500pt
\setplotsymbol ({\thinlinefont .})
\plot  2.350 22.892  2.350 22.352 /
%
%
\linethickness= 0.500pt
\setplotsymbol ({\thinlinefont .})
\plot  2.826 22.098  3.493 22.574 /
\plot  3.493 22.574  3.493 22.892 /
\plot  3.493 22.892  3.493 22.860 /
%
%
\linethickness= 0.500pt
\setplotsymbol ({\thinlinefont .})
\plot  3.493 23.749  4.064 24.066 /
\plot  4.064 24.066  4.064 22.955 /
\plot  4.064 22.955  3.493 22.574 /
%
%
\linethickness= 0.500pt
\setplotsymbol ({\thinlinefont .})
\plot  3.524 23.082  4.064 23.368 /
%
%
\linethickness= 0.500pt
\setplotsymbol ({\thinlinefont .})
\plot  1.873 23.940  3.048 24.543 /
\plot  3.048 24.543  4.032 24.066 /
%
%
\linethickness= 0.500pt
\setplotsymbol ({\thinlinefont .})
\plot  2.381 23.654  3.429 24.320 /
%
%
\linethickness= 0.500pt
\setplotsymbol ({\thinlinefont .})
\plot  2.603 24.289  3.493 23.717 /
%
%
\linethickness=1pt
\setplotsymbol ({\makebox(0,0)[l]{\tencirc\symbol{'160}}})
\setsolid
\plot  3.493 23.717  5.493 22.543 /
%
%
\linethickness= 0.500pt
\setplotsymbol ({\thinlinefont .})
\setdashes < 0.1270cm>
\plot  2.985 23.336  2.350 22.955 /
%
%
\linethickness= 0.500pt
\setplotsymbol ({\thinlinefont .})
\setsolid
\plot  2.350 22.955  0.159 21.749 /
%
%
\linethickness= 0.500pt
\setplotsymbol ({\thinlinefont .})
\plot  3.524 23.082  5.461 21.907 /
%
%
\linethickness= 0.500pt
\setplotsymbol ({\thinlinefont .})
\putrule from  2.985 24.035 to  2.985 26.003
%
%
\linethickness= 0.500pt
\setplotsymbol ({\thinlinefont .})
\setdashes < 0.1270cm>
\plot  2.985 24.003  2.985 23.336 /
\plot  2.985 23.336  3.493 23.050 /
%
%
\linethickness= 0.500pt
\setplotsymbol ({\thinlinefont .})
\setdots < 0.0953cm>
\plot  2.826 22.098  0.159 20.669 /
%
%
\linethickness= 0.500pt
\setplotsymbol ({\thinlinefont .})
\plot  2.826 22.098  5.397 20.574 /
%
%
\linethickness= 0.500pt
\setplotsymbol ({\thinlinefont .})
\plot  4.064 22.955  5.461 22.162 /
%
%
\linethickness= 0.500pt
\setplotsymbol ({\thinlinefont .})
\plot  1.873 22.543  0.127 23.400 /
%
%
\linethickness= 0.500pt
\setplotsymbol ({\thinlinefont .})
\plot  1.873 22.511  0.191 21.590 /
%
%
\linethickness= 0.500pt
\setplotsymbol ({\thinlinefont .})
\plot  4.064 24.066  5.461 24.924 /
%
%
\linethickness= 0.500pt
\setplotsymbol ({\thinlinefont .})
\plot  3.048 24.543  5.207 25.940 /
%
%
\linethickness= 0.500pt
\setplotsymbol ({\thinlinefont .})
\plot  1.873 23.940  0.095 25.051 /
%
%
\linethickness= 0.500pt
\setplotsymbol ({\thinlinefont .})
\plot  3.048 24.543  0.572 25.971 /
%
%
\linethickness= 0.500pt
\setplotsymbol ({\thinlinefont .})
\plot  4.064 24.035  5.493 23.241 /
%
%
\linethickness= 0.500pt
\setplotsymbol ({\thinlinefont .})
\plot  1.873 23.940  0.127 23.019 /
%
%
\linethickness= 0.500pt
\setplotsymbol ({\thinlinefont .})
\setdashes < 0.1905cm>
\plot  0.159 22.606  2.350 24.797 /
%
%
\linethickness= 0.500pt
\setplotsymbol ({\thinlinefont .})
\plot  2.350 25.273  0.064 22.987 /
%
%
\linethickness= 0.500pt
\setplotsymbol ({\thinlinefont .})
\plot  2.350 24.289  0.826 22.796 /
%
%
\linethickness= 0.500pt
\setplotsymbol ({\thinlinefont .})
\plot  3.461 25.622  5.397 23.622 /
%
%
\linethickness= 0.500pt
\setplotsymbol ({\thinlinefont .})
\plot  3.493 25.114  5.397 23.178 /
%
%
\linethickness= 0.500pt
\setplotsymbol ({\thinlinefont .})
\plot  3.493 24.606  5.461 22.638 /
%
%
\linethickness= 0.500pt
\setplotsymbol ({\thinlinefont .})
\plot  3.493 24.130  4.445 23.178 /
%
%
\linethickness= 0.500pt
\setplotsymbol ({\thinlinefont .})
\plot  2.350 22.447  3.270 22.415 /
%
%
\linethickness= 0.500pt
\setplotsymbol ({\thinlinefont .})
\plot  1.397 21.939  4.000 21.939 /
%
%
\linethickness= 0.500pt
\setplotsymbol ({\thinlinefont .})
\plot  0.572 21.495  4.731 21.463 /
%
%
\linethickness= 0.500pt
\setplotsymbol ({\thinlinefont .})
\plot  2.350 25.749  0.064 23.463 /
\linethickness=0pt
\putrectangle corners at  0.064 26.035 and  5.493 20.574
\endpicture
\endfig

For a 3--manifold $M$ with boundary, transversality has a similar 
description.  In this case the little cubes are all in
the interior of $M$, $f|\d M$ is a transverse map of a
surface in $C$ and the framed 1--manifolds can terminate at little squares
in $\d M$ as well as faces of little cubes.

\proc{Theorem}{\rm Transversality for $\sq$--sets}

Let $M$ be a smooth manifold (possibly with boundary)
and $C$ a $\sq$--set.  Let $f\co M\to \mo C\mo$ be a map.
Then $f$ is homotopic to a transverse map.

If $f|\d M$ is already transverse, then the homotopy can be assumed
to keep $f|\d M$ fixed.

\prf We shall first prove the theorem in the case that
$M$ is a closed surface $\Sigma$, as this case contains all
the ideas for the general case.

We start by using standard cellular approximation techniques to homotope
$f$ to meet only the 2--skeleton.  Next we make $f$ transverse to the
centres of the squares of $C$.  The result is that $f$ maps a number of
small squares in $\Sigma$ diffeomorphically onto neighbourhoods of
centres of squares in $C$.  By radial homotopies, we can assume that 
each of these small squares in fact maps onto the whole of a square in
$C$ and that the rest of $\Sigma$ now maps to the 1--skeleton.  It is clear
how to identify each little square with $I^2$ so that $f$ maps each
by a characteristic map.  Now let $\Sigma'$ denote the 
closure of $\Sigma-\{$small squares$\}$.  Then
$\Sigma'$ is a surface with boundary (with corners) and $f|\d\Sigma'$
is transverse to the centres of the edges of $C$.  By relative 
transversality and further radial homotopies, we can homotope
$f$ rel $\d\Sigma'$ so that the preimages of the centres of 
the edges of $C$ are framed 1--manifolds in $\Sigma'$ with framing
lines mapped onto the relevant edges of $C$.  Moreover we can
identify each framing line with $I^1$ so that it is mapped by
a characteristic map.  If $\Sigma_0$ now denotes the closure
of $\Sigma-\{$small squares and framed 1--manifolds$\}$
then $\Sigma_0$ is mapped to the 0--skeleton, ie each component
of $\Sigma_0$ is mapped to a vertex of $C$.
The map $f$ is now transverse.

For the general case of an $n$--manifold (perhaps with boundary), 
we can assume inductively that $f|\d M$ is already transverse and
use cellular approximation rel $\d M$ to ensure that $f$ meets 
only the $n$--skeleton.  We next make  $f$ transverse to the
centres of $n$--cells.  By radial homotopies as in the 2--dimensional
case we can assume that the closure of the preimage of the interiors
of the $n$--cells denoted $M_n$ is a collection of disjoint little $n$--cubes in
$M-\d M$ each of which can be identified with $I^n$ and is mapped 
by a characteristic map, ie $f$ is now transverse to the $n$--cells
of $C$.  Now let $M'$ be the closure of $M-M_n$
then $f|\d M'$ is transverse to the $(n-1)$--cells of $C$.  By relative 
transversality and radial homotopies we can homotope $f$ rel $\d M'$ to
be transverse to the $(n-1)$--cells.  We then proceed by downward induction
to complete the construction of a transverse map homotopic to
$f$ rel $\d M$.  Notice that the process produces compatible framings
automatically. \qed

\sh{Pulling back mock bundles by transverse maps}

Now suppose that $f\co M^m\to \mo C\mo$ is a transverse map and
$\xi^q/C$ is a mock bundle.  Then by choosing smooth CW decompositions
of each manifold $N_c$ so that the face maps $\la_c^*$ are inclusions
of subcomplexes (notation from the definition of a transverse map,
above), then the product structure on $M_c$ for each $c\in C$ defines
a smooth decomposition of $M$ so that $f$ is a smooth map (of smooth
CW complexes).  It follows from 2.2 and 2.3 that $f^*(\xi)$ is a mock
bundle of codimension $q$ over $M$ and that $E(f^*(\xi))$ is a smooth
manifold of dimension $m-q$ representing the Poincar\'e dual to
$f^*(\xi)$.  Moreover if $\xi$ is embedded in $\mo C\mo\times Q$ then
$E(f^*(\xi))$ is $\varepsilon$--isotopic to a smooth submanifold of
$M\times Q$.

The case when $\xi$ is a James bundle $\zeta^i$ of $C$ will be
particularly important for the rest of the paper.  Let $P_i$ denote
$E(f^*(\zeta^i))$.  Then $P_i$ is a smooth manifold of dimension $m-i$
which can be assumed to be embedded smoothly in $M\times\re$.
Moreover we can see from construction that the image of $P_1$ in $M$
is a framed immersed self-transverse submanifold $V$ of $M$ of
codimension 1. (This is illustrated for the case $m=2$ in figure
\Rings\ above.)  Moreover the images of $P_i$ are the $i$--tuple
points of $V$ and this is illustrated in figure \Cube.  In this figure
the image of $P_2$ is the immersed 1--manifold defined by the double
lines and $P_3$ is the 0--manifold of triple points.

The choice of terminology is explained in [\JBun; section 3] where
James bundles are related to classical James--Hopf invariants.  There
is then a connection with the results of [\KS] which also relate
generalised James--Hopf invariants to multiple points of immersions
see [\JBun; remark 3.7].  In the next section we shall establish the
connection with links and diagrams suggested by figures \Rings\ and
\Cube.

\section{Links and diagrams}

In this section we use the transversality theorem proved above to
deduce the main classification theorem stated at the start of the
paper, together with several related classification results.

We start by defining diagrams in arbitary dimensions. First we need
the concept of a self transverse immersion.  Let $\celt^p$ be the
$p$--cube $I^p$ together with the $p$ hyperplanes $x_i=\ha,
i=1,\ldots,p$ and let $T^p$ denote the union of the $p$
hyperplanes. Then $\celt^p\times I^{n-p}$ consists of an $n$--cube
with $p$ central hyperplanes meeting in an $n-p$ dimensional subspace
called the {\sl core}. An immersed smooth manifold $M$ of dimension
$n-1$ in a manifold $Q$ of dimension $n$ is called {\sl self
transverse} if each point $x$ of $M$ has a neighbourhood in $Q$ like
$\celt^p\times I^{n-p}$ in which the point $x$ corresponds to an
interior point of the core and in which $T^p\times I^{n-p}$
corresponds to the image of $M$.  The integer $p$ is called the {\sl
index} of $x$ and the set of points of index $p$ is called the {\sl
stratum of index $p$}.  The stratum of index one is locally embedded
and the closures of components of index one points will be called the
{\sl sheets} of the immersed manifold $M$.  Sheets can be locally
continued through points of higher index.  At a point of index $p$
there are locally $p$ such extended sheets which meet in a manifold of
codimension $p$.

\rk{Notation}We write $M\immersed Q$ if $M$ is a self-transverse
immersed submanifold of codimension 1 of $Q$.

\rk{Remark}Any immersed manifold of codimension 1 can be regularly 
homotoped so that it is self transverse (see Lashof and Smale [\LS]).

\rk{Definition}{\sl (Framed) diagram}

A {\sl diagram} $D$ in an $n$--manifold $Q$ is a framed immersed self
transverse submanifold $M$ such that the sheets are locally totally
ordered and this ordering is preserved in a neighbourhood of points of
higher index.  The ordering is thought of as ``vertical'' and we speak
of a sheet being ``above'' another if it follows in the order.

We call the components of $Q-M$ the {\sl regions} of the diagram and
also refer to these as the stratum of index 0 the diagram.  In general
the stratum of index $p$ of $D$ is the stratum of index $p$ of $M$ as
described above.

\rk{Example 1}{\sl Diagram on a surface}

A diagram $D$ on a closed surface $\Sigma$ is a familiar concept.  It
comprises a collection of framed immersed circles in general position
in $\Sigma$ such that at each crossing one of components is locally
regarded as the overcrossing curve and the other as the undercrossing
curve.
$$
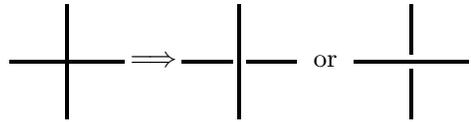
\beginpicture
\setcoordinatesystem units <.15truein,.15truein> point at 6 0
\linethickness=1pt
\putrule from -2 0 to 2 0
\putrule from 0 -2 to 0 2
\put{$\Longrightarrow$} at 3 0
\setcoordinatesystem units <.15truein,.15truein> point at 0 0
\linethickness=1pt
\putrule from -2 0 to -.25 0
\putrule from .25 0 to 2 0
\putrule from 0 -2 to 0 2
\put{\ninerm or} at 3 0
\setcoordinatesystem units <.15truein,.15truein> point at -6 0
\linethickness=1pt
\putrule from -2 0 to 2 0
\putrule from  0 .25 to 0 2
\putrule from 0 -2 to 0 -.25
\endpicture
$$
The framing can be pictured as a transverse arrow for each arc
of the diagram the direction of which is preserved through crossings.
$$
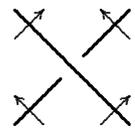
\beginpicture\ninepoint
\setcoordinatesystem units <.3truein,.3truein> point at 0 0
\setplotarea x from -1 to 1, y from -1 to 1
\arrow <3pt> [0.3, 1] from 0.5 -1 to 1 -0.5
\arrow <3pt> [0.3, 1] from -1 0.5 to -0.5 1
\arrow <3pt> [0.3, 1] from -0.5 -1 to -1 -0.5
\arrow <3pt> [0.3, 1] from 1 0.5 to 0.5 1
\setplotsymbol ({\tenrm .})
\setlinear
\plot -1 -1 -.2 -.2 /
\plot 1 1 .2 .2 /
\plot -1 1 1 -1 /
\endpicture
$$
The component of $\Sigma- D$ are the {\sl regions} of the diagram and
we call the components of $D-\{\hbox{double points}\}$ the {\sl arcs}
of $D$.  These are what we called sheets above.

If the surface is oriented then the framing determines an orientation
on the arcs of the diagram (and conversely) by the left-hand rule
illustrated:
$$
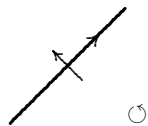
\beginpicture\ninepoint
\setcoordinatesystem units <.3truein,.3truein> point at 0 0
\setplotarea x from -1 to 1, y from -1 to 1
\arrow <4pt> [0.3, 1] from 0.5 0.5 to 0.55 0.55
\arrow <3pt> [0.3, 1] from 0.25 -0.25 to -0.25 0.25
\setplotsymbol ({\tenrm .})
\setlinear
\plot -1 -1 1 1 /
\put {$\aclock$} [lb] at 1 -1
\endpicture
$$
This framing is usually called the ``blackboard framing''.

\rk{Example 2}{\sl Diagram in a 3--manifold}

A diagram in a closed 3--manifold is a self-transverse immersed surface (ie 
transverse double curves and triple points) equipped with a
compatible ordering of sheets at double curves and triple points.
Compatible means that the ordering of sheets at the triple points
restricts to give the ordering at adjacent double arcs.  It follows
that the ordering of sheets at double arcs is preserved in a neighbourhood of 
a triple point.   So, for example, if in the following diagram 
sheet 2 is above sheet 1 then sheet $2'$ is above sheet $1'$.
%
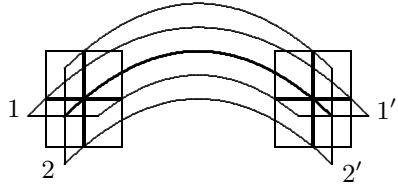
$$
\beginpicture
\setcoordinatesystem units <.2truein,.25truein> point at 0 0
\setplotarea x from -4 to 4, y from -1.5 to 2
\small
\put {$1$} [r] at -4.6 -.2
\put {$1'$} [l] at 4.6 -.2
\put {$2$} [rt] at -3.7 -1.3
\put {$2'$} [lt] at 3.7 -1.3
\putrule from -4 -1 to -2 -1
\putrule from -2 -1 to -2 1
\putrule from -2 1 to -4 1
\putrule from -4 1 to -4 -1
\putrule from 4 -1 to 2 -1
\putrule from 2 -1 to 2 1
\putrule from 2 1 to 4 1
\putrule from 4 1 to 4 -1
\putrule from -3.5 -1.3611 to -3.5 .6389
\putrule from 3.5 -1.3611 to 3.5 .6389
\putrule from -4.455528 -.3611 to -2.624652 -.3611
\putrule from 4.455528 -.3611 to 2.624652 -.3611
\linethickness=1truept
\putrule from -3 -1 to -3 1
\putrule from 3 -1 to 3 1
\putrule from -4 0 to -2 0
\putrule from 4 0 to 2 0
\setquadratic
\plot -3.5 .6389  0 2  3.5 .6389 /
\plot -3.5 -1.3611  0 0  3.5 -1.3611 /
\plot -2.624652 -.3611  0 .5  2.624652 -.3611 /
\plot -4.455528 -.3611  0 1.5  4.455528 -.3611 /
\setplotsymbol ({.})
\plot -3.5 -.3611  0 1  3.5 -.3611 /
\endpicture
$$

\rk{Remarks}(1)\qua
If $M\immersed Q$ and both $M$ and $Q$ are oriented then the
orientations determine a canonical framing of $M$ in $Q$ analogous to
the blackboard framing for a diagram on a surface.

(2)\qua In the case that $M$ or $Q$ have boundary, we assume that
diagrams are {\sl proper} ie that $M$ meets $\d Q$ in its boundary
(which thus defines a diagram in $\d Q$).

\sh{Framed embeddings and diagrams}

A diagram $D: M \immersed Q$ determines a framed embedding of $M$ in
$Q\times\re$ (ie a {\sl link\/}) by lifting the sheets at multiple
points in the $\re$ direction in the order given.  In other words a
sheet ``above'' another sheet is lifted higher.  The framing is given
by using the given framing of the diagram and taking vertically
upwards as the last framing vector.  It is clear that the resulting
link, called the {\sl lift of $D$} is well-defined up to an isotopy
which moves points vertically.

There is a converse to this process.  A framed link determines a
diagram.  This is a consequence of the compression theorem.  Let $M$
be a framed embedding in $Q\times\re$.  We call $M$ {\sl horizontal}
if the last framing vector is always vertically up.  Note that a
horizontal embedding covers an immersion in $Q$.

\proclaim{Compression Theorem}
Let $M$ be a framed embedding in $Q\times\re$, then $M$ is isotopic
(by a small isotopy) to a horizontal embedding; moreover if $M$
already horizontal in the neighbourhood of some compact set, then the
isotopy can be assumed fixed on that compact set.\endproc

The theorem follows from a deep result of Gromov [\Gr].  Direct proofs
are given in [\Comp].

\proc{Corollary}Any framed link of codimension 2 is isotopic to the
lift of some diagram.\endproc

\prf  This follows at once from the Compression Theorem and Lashof--Smale
[\LS]. \endprf

There are relative versions of both results used in the corollary so
that, for example, if $M$ and $Q$ have boundary and $M$ is embedded
properly in $Q\times\re$ such that the embedding of the boundary is
the lift of a diagram, then the diagram determined by $M$ can be taken
to extend the given diagram of the boundary.

\sh{The fundamental rack of a diagram and a link}

For details here see [\FeRo; pages 369--375].  
A diagram in $\re^2$ determines a {\sl fundamental rack} by labelling
the arcs by generators and reading a relation at each crossing:
$$
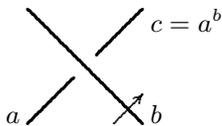
\beginpicture\ninepoint
\setcoordinatesystem units <.3truein,.3truein> point at 0 0
\put { } [rB] at -1.1 1.2
\put {$a$} [rB] at -1.1 -1
\put {$b$} [lB] at 1.1 -1
\put {$c=a^b$} [lt] at 1.1 1
\arrow <3pt> [0.3, 1] from 0.5 -1 to 1 -0.5
\setplotsymbol ({\tenrm .})
\setlinear
\plot -1 -1 -.2 -.2 /
\plot 1 1 .2 .2 /
\plot -1 1 1 -1 /
\endpicture\leqno{\bf \label}\key{cross-rel}
$$
More generally any diagram determines a rack in a similar way.
Components of index 1 of $D$ are labelled by generators and a relation
is read from each component of index 2 by the same rule as in the
2--dimensional case (think of a perpendicular slice); points of higher
index are not used (cf [\FeRo; remark (2), page 375]).

(Note that there are really two versions of diagram \ref{cross-rel}
because the framing of the understring is immaterial.)

\rk{Notation}If $D$ is a diagram, then we denote the fundamental rack
of $D$ defined as above by $\Gamma(D)$.\medskip

A framed codimension 2 embedding $L$ also determines a {\sl
fundamental rack} denoted $\Gamma(L)$, see the next definition for
details; moreover if $D$ is a diagram in $\re^n$ and $L$ is a lift in
$\re^{n+1}$, then $\Gamma(D)$ can be naturally identified with
$\Gamma(L)$ [\FeRo; theorem 4.7 and remark (2), page 375].

We need to interpret the rack $\Gamma(D)$ in the case that the diagram
is not in $\re^n$.  This case was not covered in [\FeRo].

\rk{Definition}{\sl The reduced fundamental rack}

Let $L:M\subset Q\times\re$ be a framed codimension two embedding (ie
a link) and choose a basepoint $*\in Q\times\re-M$.  Consider paths
$\alpha$ from the parallel (framing) manifold of $M$ to $*$ in
$Q\times\re-M$.  Recall that the fundamental rack $\Gamma(L)$
comprises the set of homotopy classes of these paths with the rack
operation $a^b$, where $a=[\alpha], b=[\beta]$, given by the class of
the composition of $\alpha$ with the ``frying pan'' loop determined by
$\beta$, namely $\bar{\beta}\circ \mu\circ \beta$ where $\mu$ is the
meridian at the start of $\beta$ [\FeRo; page 359].  To define the
reduced fundamental rack $\bar{\Gamma}(L)$ we kill the action of
$\pi_1(Q)$.  More precisely, two paths $\alpha,\beta$ starting from
the same point are equivalent if $\bar{\beta}\circ\alpha$ is a product
of conjugates of elements of $\pi_1(Q)$ in $\pi_1(Q\times\re-M)$,
where $Q$ is identified with $Q\times {1}$ and we assume that the link
lies below level $1$.  This is extended to homotopy classes in the
obvious way.

It can be checked that this is an equivalence relation and that the
rack operation is defined on equivalence classes.  The resulting rack is
the {\sl reduced fundamental rack\/} $\bar{\Gamma}(L)$.

There is a simple interpretation of the reduced rack.  Replace $\re$
by $[-1,1]$ (assume that the link lies between levels $-1$ and 1).
Now define $\bar Q$ to be $Q\times[-1,1]/Q\times 1$ (ie a copy of the
cone on $Q$) based at the image of $Q\times 1$.  Note that the
definition of the fundamental rack does not use that the bigger
manifold is actually a manifold.  Thus we can define the fundamental
rack $\Gamma(M\subset \bar Q)$.  It is not hard to check that
$\bar{\Gamma}(L)= \Gamma(M\subset \bar Q)$.

Note that if $Q$ is simply connected then the reduced fundamental rack
coincides with the usual fundamental rack.

\proc{Lemma}\key{reduced-ident}Let $D: M\immersed Q$ be a diagram 
and let $L$ be the link given by lifting $D$ to $Q\times\re$.  Then
$\Gamma(D)$ can be naturally identified with
$\bar{\Gamma}(L)$.\endproc

\prf  Using the interpretation of $\bar{\Gamma}(L)$ as 
$\Gamma(M\subset \bar Q)$, the proof given in [\FeRo; pages 372--375]
adapts with obvious changes. \endprf

\sh{Labelling diagrams by racks}

Let $D$ be a diagram.  We say that $D$ is {\sl labelled} by a rack $X$
if each component of the stratum of index $1$ of $D$ is labelled by an
element of $X$ with compatibility at strata of index 2 given by 
rule \ref{cross-rel} on a perpendicular slice, where $a,b,c$ now
denote elements of $X$ (rather than generators of $\Gamma(D)$) and
$c=a^b$ in $X$.

There is a natural labelling of any diagram by its fundamental
rack, which we call the {\sl identity labelling}.  More generally
we have the following result which follows immediately from
definitions:

\proc{Lemma}A labelling of a diagram $D$ by a rack $X$ is
equivalent to a rack homomorphism $\Gamma(D)\to X$.\qed

The lemma implies that labelling is functorial in the sense that
a rack homomorphism $X\to Y$ induces a labelling by $Y$, it also
implies that labelling is really a property of the lift, not the
diagram:

\proc{Corollary}A labelling of a diagram by a rack $X$ is equivalent
to a rack homomorphism to $X$ of the reduced fundamental rack of 
the lifted framed embedding. \qed\endproc

We often speak of a framed link having a {\sl representation in $X$}
to mean that the the fundamental rack has a homomorphism to $X$.

\sh{Labelling diagrams by $\tsq$--sets}

A diagram is {\sl labelled by a $\sq$--set $C$} if, for each $p$, each
component of the stratum of index $p$ is labelled by a $p$--cube of
$C$ with compatibility conditions.  We shall explain these in detail
for diagrams in surfaces and 3--manifolds.  The general case is a
straightforward extension.

A diagram $D$ on a surface is {\sl labelled} by a $\sq$--set $C$ if:

(1)\qua The regions are labelled by vertices of $C$.

(2)\qua The arcs are labelled by edges of $C$ compatibly with
adjacent regions.  This means that the edge labelling an arc
$\alpha$ is attached (in $C$) to the two vertices labelling the
adjacent regions.  To decide which vertex labels which side, 
identify the edge with the transverse framing arrow:
$$
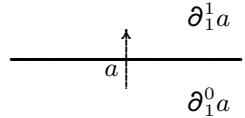
\beginpicture\ninepoint
\setcoordinatesystem units <.3truein,.3truein> point at 0 0
\setplotarea x from -2 to 2, y from -0.5 to 0.5
\arrow <3pt> [0.3, 1] from 0 -0.5 to 0 0.5
\setplotsymbol ({\tenrm .})
\setlinear
\plot -2 0 2 0 /
\put {$a$} [rt] at -0.1 -0.1
\put {$\d^1_1a$} [bl] at 1 0.5
\put {$\d^0_1a$} [tl] at 1 -0.5
\endpicture
$$
(3)\qua The double points are labelled by squares of $C$
compatibly with adjacent arcs.  This means that the four
adjacent arcs are labelled by the four 1--faces of the square.
The rule for determining which face labels which arc is this:
position the a copy of the standard square (ie $I^2$)
at the double point with faces oriented
correctly by the framing arrows and axis 1 parallel to the
overcrossing.  Now the four faces intersect the appropriate
adjacent arcs, see the diagram below, where we have
drawn both possible orientations for the square:
$$
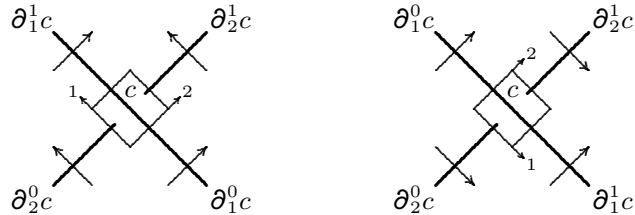
\beginpicture\ninepoint
\setcoordinatesystem units <.4truein,.4truein> point at 0 0
\setplotarea x from -1 to 1, y from -1 to 1
\plot -0.5 0 0 0.5 0.5 0 0 -0.5 -0.5 0 /
\arrow <2pt> [0.3, 1] from 0 0.5 to 0.15 0.65
\arrow <2pt> [0.3, 1] from 0 -0.5 to 0.15 -0.65
\put {$\scriptstyle 1$} [lt] at 0.15 -0.65
\put {$\scriptstyle 2$} [lb] at 0.15 0.65
\put {$c$} [b] <0pt, 4pt> at 0 0
\put {$\d^1_2c$} [bl] at 1 1
\put {$\d^0_2c$} [tr] at -1 -1
\put {$\d^1_1c$} [tl] at 1 -1
\put {$\d^0_1c$} [br] at -1 1
\arrow <3pt> [0.3, 1] from 0.5 -1 to 1 -0.5
\arrow <3pt> [0.3, 1] from -1 0.5 to -0.5 1
\arrow <3pt> [0.3, 1] from -1 -0.5 to -0.5 -1
\arrow <3pt> [0.3, 1] from 0.5 1 to 1 0.5
\setplotsymbol ({\tenrm .})
\setlinear
\plot -1 -1 -.2 -.2 /
\plot 1 1 .2 .2 /
\plot -1 1 1 -1 /
\setcoordinatesystem units <.4truein,.4truein> point at 5 0
\setplotarea x from -1 to 1, y from -1 to 1
\setplotsymbol ({\fiverm .})
\plot -0.5 0 0 0.5 0.5 0 0 -0.5 -0.5 0 /
\arrow <2pt> [0.3, 1] from 0.5 0 to 0.65 0.15 
\arrow <2pt> [0.3, 1] from -0.5 0 to -0.65 0.15
\put {$\scriptstyle 1$} [rb] at  -0.65 0.15
\put {$\scriptstyle 2$} [lb] at  0.65 0.15
\put {$c$} [b] <0pt, 4pt> at 0 0
\put {$\d^1_2c$} [bl] at 1 1
\put {$\d^0_2c$} [tr] at -1 -1
\put {$\d^0_1c$} [tl] at 1 -1
\put {$\d^1_1c$} [br] at -1 1
\arrow <3pt> [0.3, 1] from 0.5 -1 to 1 -0.5
\arrow <3pt> [0.3, 1] from -1 0.5 to -0.5 1
\arrow <3pt> [0.3, 1] from -0.5 -1 to -1 -0.5 
\arrow <3pt> [0.3, 1] from 1 0.5 to 0.5 1 
\setplotsymbol ({\tenrm .})
\setlinear
\plot -1 -1 -.2 -.2 /
\plot 1 1 .2 .2 /
\plot -1 1 1 -1 /
\endpicture
$$
Let $D$ be a diagram in a 3--manifold $M$.  Call the components of
the double curves minus triple points {\sl double arcs}, the components
of the surface minus the double curves {\sl sheets} and the components
of $M$ minus the surface {\sl regions}.  

A labelling of $D$ by a $\sq$--set $C$ is a labelling of regions
(resp.\ sheets, double arcs, triple points)
by vertices (resp.\ edges, squares, 3--cubes) of $C$ subject to
compatibility conditions at sheets, double curves and triple points.

The compatibility conditions at sheets and double curves are the same as
for a 2--dimensional diagram (imagine working in a transverse slice)
whilst the compatibility condition at a triple point 
can be described as follows.  Say the \sl positive \rm side of a sheet
coincides with the head of its framing arrow. Suppose a triple point
is labelled by $x\in C_3$.  Suppose a nearby double curve is labelled by
$y\in C_2$ and the missing sheet is number $i\in\{1,2,3\}$. Then
$\d_i^\epsilon x=y$ where $\epsilon=1$ if the double curve is on the
positive side  of the missing sheet and $\epsilon=0$ otherwise. 
(Notice that a small copy of the standard cube $I^3$ can be placed at 
a triple point by orienting axes according to framing of sheets
and ordering axes so that axis $i$ is perpendicular to 
the $i^{th}$ sheet.  If this is done then 
the small cube meets an adjacent double arc in the face given by
the above labelling rules.)

These compatibility conditions extend to a general diagram in the
obvious way.  If a component of the stratum of index $p$ is labelled
by $c\in C_p$ then the neighbouring components of the index $(p-1)$
stratum are labelled by $\d_i^\ep(c)$ with the rule for determing $i$
and $\ep$ being analogous to the 3--dimensional case.  

\rk{Remarks}

(1)\qua If $C=BX$, where $X$ is a rack, then labelling in $C$ is
precisely the same as labelling in $X$.  This is the same as labelling
index 1 components by elements of $X$ with the usual compatibility
requirement at index 2 points (diagram \ref{cross-rel}).  Points of
higher index play no part because a 3--cube of $C=BX$ is determined by
its faces (see the discussion in section 1, the key example).  Also
there is no need to label index 1 components explicitly since a square
in $BX$ is determined by its edges.  The compatibility condition
ensures that the required square exists.

(2)\qua Every diagram has the {\sl trivial} labelling, namely
labelling by $T$, the trivial $\sq$--set (with one cell of each 
dimension).  This can also be regarded as labelling by the
trivial rack (with one element).

(3)\qua Labelling is functorial:  Given a diagram labelled in $C$
and a $\sq$--map $f\co C\to D$ then $f$ transforms the labels to
a labelling in $D$.

(4)\qua There is also the concept of labelling by a {\sl trunk} $T$ in
other words labelling by the nerve $\N T$ (see [\Trunks]).  In this
case regions are labelled by vertices and index 1 components by edges
(between the vertices labelling adjacent regions) such that at index
two components the adjacent index 1 components form a preferred
square.  As for racks, points of higher index play no part in the
labelling.

(5)\qua The case of labelling by the {\sl action rack space\/} $B_YX$
[\Trunks; example 3.1.2] is worth describing in detail.  This is the
same as labelling the diagram by the rack $X$ togther with a {\sl
regional labelling} by $Y$.  Here $Y$ is a set on which $X$ acts (see
[\Trunks; above 1.4]).  In other words regions are labelled in $Y$
compatibly with the labelling on sheets: if a region is labelled $a\in
Y$ and a sheet labelled $b\in X$ is crossed (in the framing direction)
then the region on the other side is labelled $a^b$ (the action of $b$
on $a$).  An important special case is when $X=Y$ and the action is
the rack operation and all labels lie in $X$.  This type of labelling
was used in [\Tref], with X the three colour rack, to distingush left
and right trefoils.  The space $B_XX$ is called the {\sl extended rack
space\/}.

\sh{Labelling and transversality}

In section 2 we defined a transverse map of a manifold $M$ in a
$\sq$--set $C$ and at the end of the section we observed that such a
map defines a framed self-transverse immersed submanifold $V$ of $M$
(see figures \Rings, \Cube), which can be seen as the image of the
pull-back of the first James bundle $\zeta^1(C)$.  We also observed
that $\zeta^1(C)$ embeds in $\mo C\mo\times\re$ and hence this
immersed submanifold is covered by an embedding in $M\times\re$ in
other words it is a {\sl diagram}.  Moreover this diagram is labelled
by $C$ in a natural way.  Recall that each component of index $p$ is
surrounded by a $p$--cube bundle, the fibres of which are mapped to a
$p$--cube of $C$.  Label the component by this cube.  It can readily
be checked that this labelling is compatible, indeed the definition of
compatibility (above) was set up precisely in this way.

Thus a transverse map {\sl to} $C$ determines a diagram labelled {\sl
by} $C$.  We now describe the converse process which is given by the
construction of a {\sl neighbourhood system} for the diagram.  To help
understand the somewhat complicated construction, we shall deal with
the cases $n=2,3$ in detail first.

\rk{2--dimensional case}

Suppose that we are given a diagram in a surface $\Sigma$ labelled by
$C$.  We construct a neighbourhood system for the diagram by drawing
little squares around the double points and continue to construct
bi-collars around the arcs.
$$
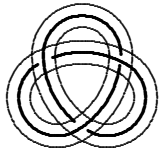
\beginpicture
\setcoordinatesystem units <.07mm,.07mm> point at 0 0
\put { } at 0 150
\put { } at 0 -120
\circulararc 180 degrees from 87.0 50.0 center at 0 50.0
\circulararc 180 degrees from 47.0 50.0 center at 0 50.0
\circulararc 180 degrees from -87.0 50.0 center at -43.5 -25.0
\circulararc 180 degrees from -67 15.2 center at -43.5 -25.0
\circulararc -180 degrees from 87.0 50.0 center at 43.5 -25.0
\circulararc -180 degrees from 67  15.2 center at 43.5 -25.0
\circulararc 60 degrees from 87.0 50.0 center at 0 -100.0
\circulararc 60 degrees from 67  15.2 center at 0 -100.0
\circulararc -60 degrees from 87.0 50.0 center at -87.0 50.0
\circulararc -60 degrees from 47.0 50.0 center at -87.0 50.0
\circulararc 60 degrees from -87.0 50.0 center at 87.0 50.0
\circulararc 60 degrees from -47.0 50.0 center at 87.0 50.0
\setplotsymbol ({$\tenrm .$})
\circulararc -55 degrees from 65.0 30.0 center at -87.0 50.0
\circulararc 180 degrees from -77 32.6 center at -43.5 -25.0
\circulararc 50 degrees from 77 32.6 center at 0 -100.0
\circulararc -180 degrees from 77 32.6 center at 43.5 -25.0
\circulararc 50 degrees from -67.0 50.0 center at 87.0 50.0
\circulararc 180 degrees from 67.0 50.0 center at 0 50.0
\endpicture
$$
This determines a transverse map into $C$ by mapping the regions
outside the squares and bi-collars to the labelling vertex, collapsing
the bi-collars onto fibres and mapping to the labelling edge and
finally mapping the little squares to the labelling squares for the
double points, using the orientations for edges and squares described
in the definition of labelling by a $\sq$--set (above).

It is clear that this construction is unique up to minor choices which
only affect the neighbourhood system up to an ambient isotopy fixing
the diagram setwise.  To be precise, define two transverse maps of
$\Sigma$ in $C$ to be {\sl diagram isotopic} if they differ by an
ambient isotopy fixing the pull-back diagram setwise then we have a
well-defined process for turning a labelled diagram into a diagram
isotopy class of transverse maps.

\rk{3--dimensional case}

Given a diagram in a 3--manifold $M$ we construct a neighbourhood
system as follows.  We choose little 3--cubes around the triple points
which meet neighbouring sheets in the three central 2--cubes.  Each of
these can be identified with $I^3$ in a canonical way using the
ordering of sheets as described above.  We call the portions of the
double arcs outside these cubes, {\sl reduced double arcs}.  Next
construct a trivial bundle with fibre a square around each reduced
double arc such that each square meets neighbouring sheets in the
central cross and which fits onto a face of the relevant 3--cube at
boundary points (precisely how to construct these trivial bundles will
be explained in lemma 3.6 below).  Finally construct bi-collars
around the reduced sheets (outside the square bundles) which fit onto
edges of the squares.  (See figure \Cube\ for a view of part of this
construction.)

Now map the 3--cubes to the labelling 3--cube of $C$,
collapse the square bundles onto a single square and map to
the labelling square, and likewise collapse the bi-collars and map
to labelling edges and finally map the reduced regions to labelling 
vertices.  The result is the required transverse map to $C$.

The only element of choice in the contruction was the choice of the
neighbourhood system.  In the next lemma we show how to choose this
system using collaring arguments.  Using uniqueness of collars we can
then see that the system is unique up to ambient isotopy fixing the
diagram setwise.

\proc{Lemma}{\rm Constructing neighbourhood systems by collars}

Let $D$ be a diagram in a 3--manifold $M$.  Then $D$ has a neighbourhood
system and any two such systems differ by an ambient isotopy of $M$
fixing $D$ setwise.

\prf
Start by choosing a collar for each triple point in each double arc
(ie an interval).  Now concentrate on a particular triple point $p$ and
label the extended sheets near $p$ 1,2,3 (in the order given by the diagram)
and use the label 12 for example for the extended double arc $1\cap 2$ etc.
The collar on $p$ in 12 is a double interval $J$ say.  Choose a bi-collar on
$J$ in 1 extending the chosen collar in 13.  This defines a square $S$ say
in $1$.  Do the same for 12 in 2 and 13 in 3.  We now have three mutually
perpendicular squares.  Complete the little cube at $p$ by choosing a bi-collar
on $S$ in $M$ extending the chosen (partially defined) collars on 12 in 2
and 13 in 3.

To define the square bundle on a reduced double arc $\alpha$ say, choose
bi-collars in both intersecting sheets (extending collars given by the
little cubes on $\d J$) and then extend one of these to a bi-collar on
the total space of the other collar (extending the collars given by the
squares over $\d J$).   Finally construct the bi-collars over the
reduced sheets extending the already constructed collars (given by the
square bundles) over the boundary.

It is clear that a neighbourhood system defines all the above collars
and the uniqueness part of the lemma now follows from uniqueness of
collars.\qed

\rk{The general case}

The extension of the case $n=3$ to the general case is
straightforward: a neighbourhood system for a diagram $D$ labelled in
$C$ is constructed by choosing little $n$--cubes around each point of
index $n$ meeting nearby sheets in central $(n-1)$--cubes.  Then the
faces are extended to trivial $(n-1)$--cubes bundles around the
reduced index $n-1$ strata and the process is completed by downward
induction on index.  Uniqueness is proved in the same way as the case
$n=3$.  The neighbourhood system and the labelling in $C$ determines
a transverse map to $C$ unique up to {\sl diagram isotopy} that is up to an
isotopy fixing the diagram setwise.

It is clear that the two processes: transverse map to labelled diagram
and labelled diagram to transverse map (via neighbourhood system) are
inverse and we can summarise this in the following result.

\proc{Proposition}There is a bijection between 
labelled diagrams in $M$ labelled by a $\sq$--set $C$ and 
diagram isotopy classes of transverse maps of $M$ in $C$.

The bijection is given by pulling back the first James bundle
$\zeta^1(C)$.\qed

\rk{Remarks}(1)\qua Note that the lift of the diagram in $M\times\re$
is also obtained by pulling back the embedded first James bundle
$\zeta^1(C)\subset C\times\re$.

(2)\qua There is a relative version of the proposition for the case
$M$ has boundary: the restriction of the bijection to the boundary
coincides with the bijection for the boundary. \medskip

In order to interpret $\pi_n$ of a cubical set, we need a
based version of the proposition.  Choose basepoint $*\in S^n$ and
base vertex $*\in C_0$ and identify $S^n-\{*\}$ with $\re^n$.  

\proc{Proposition}
There is a bijection between labelled diagrams in $\re^n$ labelled by
a $\sq$--set $C$ such that the non-compact region is labelled by the
vertex $*\in C_0$ and diagram isotopy classes of based transverse maps of
$S^n$ in $C$.

The bijection is given by pulling back the first James bundle
$\zeta^1(C)$.\qed

\sh{The classification theorems}

We now interpret homotopy classes of maps into $C$ and in particular
$\pi_n(C)$.  To do this we shall need the following definition.

\rk{Definition}{\sl Cobordism of diagrams}

Diagrams $D_0$, $D_1$ in $M$ are {\sl cobordant} if there is
a diagram $D$ in $M\times I$ which meets $M\times \{0,1\}$ in
$D_0$, $D_1$ respectively.  We call $D$ the {\sl cobordism} between
$D_0$ and $D_1$.  It can readily be checked that diagram cobordism
is an equivalence relation.

There is a similar notion of cobordism for labelled diagrams and we
denote the set of cobordism classes of diagrams in $M$ labelled in a
$\sq$--set $C$ by $\D(M,C)$.  

Now let $C$ be a $\sq$--set with a basepoint $*\in C_0$.  Define the
set $\D(n,C)$ to be the set of labelled cobordism classes of diagrams
in $\re^n$ (labelled by $C$) such that the non-compact region is
labelled by the vertex $*$.

The set of diagrams in $\re^n$ has an addition given by
disjoint union.  To be precise, given diagrams $D_1,D_2$ choose 
copies in disjoint half spaces, then define $ D_1+D_2$ to be
$D_1\amalg D_2$.   This addition is well-defined up to cobordism,
is compatible with cobordism and makes $\D(n,C)$ into an
abelian group.

In the case when $C=BX$, a rack space, we abbreviate $\D(M,BX)$ and
$\D(n,BX)$ to $\D(M,X)$ and $\D(n,X)$ respectively.

\proc{Theorem}{\rm Classification of labelled diagrams}\key{diag-class}

Let $C$ be a $\sq$--set.  There is a natural bijection between the set
of homotopy classes of maps $[M,|C|]$ and $\D(M,C)$.  If $C$ has
basepoint $*\in C_0$.  There is a natural isomorphism between
$\pi_n(C)$ and $\D(n,C)$.

The isomorphisms can be described as given by pulling back the first
James bundle $\xi^1(C)$ which thus plays the r\^ole of classifying
bundle for labelled diagrams.
\endproc

\prf By proposition 3.7 a labelled diagram in $M$ determines a 
(transverse) map of $M$ in $C$.  Similarly by a cobordism determines a
homotopy.  Thus we have a function $\Phi\co\D(M,C)\to [M,|C|]$.  By
transversality (theorem 2.5) $\Phi$ is a bijection.  For the second
bijection we use 3.8 instead of 3.7.  The rest of the theorem follows
from definitions.
\endprf

We now specialise to the case when $C$ is the rack space $BX$ and
labelling in $C$ is the same as labelling in $X$ and is a property of
the link which covers the diagram.  Note that $BX$ is based at the
unique 0--cell.

We need to define cobordism of links:

\rk{Definition}{\sl Cobordism of links}

We say that framed links $L_0, L_1$ in $W$ are {\sl cobordant} if
there is a framed link $L$ properly embedded in $W\times I$ which
meets $W\times \{0,1\}$ in $L_0$, $L_1$ respectively.

Let $X$ be a rack, then there is an analogous notion of cobordism of
links with representation in $X$, namely a cobordism with a
representation in $X$ (ie a homomorphism of the fundamental
rack in $X$) extending the given representations on the ends.

It can readily be checked that cobordism is an equivalence relation
and we denote the set of cobordism classes of framed links in $W$
with representation in $X$ by $\F(W,X)$.

In the case that $W=M\times\re$ then we can consider representations
using {\sl reduced} fundamental rack as defined earlier; we use the
notation $\bar{\F}(M\times\re,X)$ for the set of framed links in
$M\times\re$ with homomorphism of reduced fundamental rack in $X$ up
to cobordism also with homomorphism of reduced fundamental rack in
$X$.

In the case that $W=\re^{n+1}$ we abbreviate the notation
$\F(\re^{n+1},X)$ to $\F(n+1,X)$.

There is an addition on $\F(n+1,X)$ given by disjoint union.  To be
precise, given links $L_1,L_2$ choose copies in disjoint half spaces,
then define $ L_1+L_2$ to be $L_1\amalg L_2$.  We need to explain how
to represent $L_1\amalg L_2$ in $X$.  The simplest way to do this is
to use diagrams.  If $L_i$ is given by the diagram $D_i$ then
$L_1\amalg L_2$ is given by $D_1\amalg D_2$ and representations
correspond to labellings of $D_1,D_2$ which define a labelling of
$D_1\amalg D_2$ in the obvious way.  We can also define the required
representation using the fact that the fundamental rack of $L_1\amalg
L_2$ is the free product of the fundamental racks of $L_1,L_2$ (see
[\FeRo; page 357]).  The homomorphisms of the two factors determine a
homomorphism on the entire rack.

The addition  on $\F(n+1,X)$ is well-defined up to cobordism,
is compatible with cobordism and makes $\F(n+1,X)$ into an
abelian group.

\proc{Proposition}\key{diag-link}There is a bijection 
$\D(M,X)\to\bar{\F}(M\times\re,X)$ and an isomorphism
$\D(n,X)\to\F(n+1,X)$ both induced by lifting diagrams.\endproc

\prf This follows from definitions and corollaries 3.1 and 3.5. \endprf

Combining the last two results, we deduce our main classification
theorem:

\proc{Theorem}{\rm Classification of links}\key{link-class}

Let $X$ be a rack. There is a natural bijection between $[M,\mo
BX\mo]$ and $\bar{\F}(M\times\re,X)$ and there is a natural isomorphism
between $\pi_n(BX)$ and $\F(n+1,X)$.  The embedded first James bundle
plays the r\^ole of classifying bundle in both cases.
\endproc

\sh{Double cobordism}

Diagrams are {\sl doubly cobordant\/} if they are cobordant by a
simultaneous cobordism of diagram and the containing manifold.  More
precisely, suppose that $D_i: V_i\immersed M_i$ is a diagram for
$i=1,2$.  Then $D_1$ is doubly cobordant to $D_2$ if there is a
diagram $U\immersed W$ with boundary, such that the boundary is
the disjoint union $D_1\amalg D_2$.

There is a similar notion of double cobordism of diagrams with
labelling in a $\sq$--set or a rack and there is an analogous notion
of double cobordism of links possibly with representation in a rack.
We shall consider the special case of {\sl product links\/} by which
we mean links in $L:V\subset M\times\re$ up to double cobordism of $V$
and $M$ with representation of reduced fundamental rack in a given
rack $X$.

These sets all form abelian groups under disjoint union.

The proofs of theorems \ref{link-class} and \ref{diag-class} extend to
prove the following theorem.

\proc{Theorem}{\rm Classification up to double cobordism}
\key{double-class}

The set of double cobordism classes of diagrams in $n$--manifolds
labelled in the $\sq$--set $C$ is in natural bijection with ${\frak
N}_n(C)$ (the unoriented bordism group).  If the containing manifolds
are oriented and we use oriented bordism, then the set is in bijection
with ${\Omega}_n(C)$ (the oriented bordism group).

The set of double cobordism classes of product links in $n$--manifolds
cross $\re$ with representation of reduced fundamental rack in a given
rack $X$ is in natural bijection with ${\frak N}_n(BX)$.  If all manifolds
are oriented then the set is in bijection
with ${\Omega}_n(BX)$.\endproc

\sh{Calculations}

Later in the paper (section 5) we shall report on calculations of
homotopy and homology of rack spaces.  Any such calculation gives an
immediate calculation of a link group using the appropriate
classification theorem above.  We shall not spell out all such
corollaries, but here are a couple of sample results:

\medskip
(1)\qua{\sl Let $K$ be the trefoil knot with any framing.  Any link in
$\re^3$ with representation in $\Gamma(K)$ is cobordant (with
representation) to the disjoint union of $n$ copies of $K$ with
identity representation.}

\medskip
(2)\qua{\sl The group of double cobordism classes of three coloured
product links in oriented 3--manifolds is isomorphic to
$\Z\oplus\Z_3$.  Hence there is a particular three coloured link which
is non-trivial under double cobordism, but for which the disjoint
union of three copies is trivial.}

\medskip
Here {\sl three colouring\/} means representation in the {\sl three
colour rack\/} $D_3=\{0,1,2\mid a^b=c$ iff $a,b,c$ are all the same
or all different$\}$.

The first result follows from theorems \ref{link-class} and 5.4.
Theorem 5.4 implies that $\pi_2$ of the rack space of the fundamental
rack of any irreducible link in $\re^3$ is $\Z$.  (Result (1) is
therefore true with the trefoil replaced by any irreducible link.)

The second result follows from theorem \ref{double-class}, the fact that
$\Omega_3(X)\cong H_3(X)$ for any space $X$, and the calculation of
$H_3(BD_3)$ as $\Z\oplus\Z_3$ given in [\Tref].

\section{The classical case}

We now turn to the lowest non-trivial dimension ($n=2$).  In this case
the cobordism classes can readily be described as equivalence classes
under simple moves and this gives a combinatorial description of
$\pi_2(C)$ which can be used for calculations.

\sh{Framed embeddings and diagrams}

Up to isotopy any framed embedding in $\re^3$ is the lift of a diagram
in $\re^2$.  This is seen by choosing any diagram to represent the
(unframed) link and then correcting the framing by introducing twists
(Reidemeister 1--move):
$$
\beginpicture
\setcoordinatesystem units <.12truein,.12truein> point at 3 0
\setplotarea x from -2 to 2, y from -2 to 2
\put { } at 0 0
\savelinesandcurves on "omega1.dat"
\writesavefile {unlabelled omega1 move}
\setquadratic
\setplotsymbol ({\tenrm .})
\plot 0 2  .0625 1.75  .125 1.25  .25 .5  .45 -.05  .75 -.375  1.25 -.5  
1.75 -.375  2 0  1.75 .375  1.25 .5  .75 .375  .6 .25 /
\plot .25 -.25  .125 -1.25  0 -2 /
\put {$\Leftrightarrow$} at 3.5 0
\setcoordinatesystem units <.12truein,.12truein> point at -2 0
\setplotarea x from -2 to 2, y from -2 to 2
\put {$\Leftrightarrow$} at 2 0
\linethickness=1pt
\putrule from 0 -2.1 to 0 2.1 
\setcoordinatesystem units <.12truein,.12truein> point at -6 0
\setplotarea x from -2 to 2, y from -2 to 2
\setquadratic
\setplotsymbol ({\tenrm .})
\plot 0 -2  .0625 -1.75  .125 -1.25  .25 -.5  .45 .05  .75 .375  1.25 .5  
1.75 .375  2 0  1.75 -.375  1.25 -.5  .75 -.375  .6 -.25 /
\plot .25 .25  .125 1.25  0 2 /
\put { } at 3 0
\endpicture\leqno{\bf R1}
$$
(This argument is the proof of the Compression Theorem in this easy
case.)  The resulting diagram is unique up to regular homotopy (or
equivalently Reidemeister 2 and 3 moves) together with the double
Riedemeister 1--move illustrated below. For a proof see [\FeRo; pages
369-370].  We shall refer to these moves as R2, R3 and R11 (pronounced
r-one-one) respectively:

{\bf R2}
\hglue 1in\lower10pt\hbox{{\hsize=68pt\b{ahed\nl fcbg\nl oooo}
\smash{\raise10pt\hbox{$\Leftrightarrow$}}\ \ \qquad \b{iiii\nl oooo\nl iiii}}}

{\bf R3}
\hglue 0.7in\lower15pt\hbox{\hsize=24pt\stwoi\sone\stwo\qquad 
\smash{\raise15pt\hbox{$\Leftrightarrow$}} \qquad 
\sone\stwo\sonei}
 $$
\beginpicture
\setcoordinatesystem units <.12truein,.08truein> point at 3 0
\setplotarea x from -2 to 2, y from -2 to 2
\linethickness=1pt
\putrule from 6.5 -2.1 to 6.5 6.1
\put { } at 16 0 
\put {$\Leftrightarrow$} at 4 2 
\put {$\Leftrightarrow$} at 9 2 
\setquadratic
\setplotsymbol ({\tenrm .})
\plot 0.1 2  .0625 1.75  .125 1.25  .25 .5  .45 -.05  .75 -.375  1.25 -.5  
1.75 -.375  2 0  1.75 .375  1.25 .5  .75 .375  .6 .25 /
\plot .25 -.25  .125 -1.25  0 -2 /
\setcoordinatesystem units <.12truein,.08truein> point at 3 -4
\setplotarea x from -2 to 2, y from -2 to 2
\setquadratic
\setplotsymbol ({\tenrm .})
\plot 0.1 -2  .0625 -1.75  .125 -1.25  .25 -.5  .45 .05  .75 .375  1.25 .5  
1.75 .375  2 0  1.75 -.375  1.25 -.5  .75 -.375  .6 -.25 /
\plot .25 .25  .125 1.25  0 2 /
\setcoordinatesystem units <.12truein,.08truein> point at -8.5 -4
\setquadratic
\setplotsymbol ({\tenrm .})
\plot 0.1 2  .0625 1.75  .125 1.25  .25 .5  .45 -.05  .75 -.375  1.25 -.5  
1.75 -.375  2 0  1.75 .375  1.25 .5  .75 .375  .6 .25 /
\plot .25 -.25  .025 -1.25  0 -2 /
\setcoordinatesystem units <.12truein,.08truein> point at -8.5 0
\setplotarea x from -2 to 2, y from -2 to 2
\setquadratic
\setplotsymbol ({\tenrm .})
\plot 0.1 -2  .0625 -1.75  .125 -1.25  .25 -.5  .45 .05  .75 .375  1.25 .5  
1.75 .375  2 0  1.75 -.375  1.25 -.5  .75 -.375  .6 -.25 /
\plot .25 .25  .025 1.25  0 2 /
\endpicture\leqno{\bf R11}
$$
There is a similar result with the same proof for any surface.

\sh{Classification by cobordism}

Recall that $\D(2,C)$ is the (group) of labelled cobordism classes of
diagrams in $\re^2$ labelled in the based $\sq$--set $C$ (such that
the infinite region is labelled by $*$).  We recall Theorem
\ref{diag-class} in this case:

\proc{Theorem}{\rm Classification of labelled diagrams}\key{2diag-class}

Let $C$ be a $\sq$--set with basepoint $*\in C_0$.  There is a natural
isomorphism between $\pi_2(C)$ and $\D(2,C)$. \endproc

The theorem cuts both ways.  It classifies diagrams up to cobordism
and also provides an interpretation of $\pi_2(C)$ which can be used
for calculations.  For this purpose we need to break a cobordism into
a sequence of combinatorial moves which we now describe.

\sh{Cobordism by moves}

Suppose that we are given a cobordism (a diagram $D$ in $\re^2\times I$)
between diagrams $D_0,D_1$.  Think of $\re^2\times I$ as a sequence
of copies of $\re^2$.  This breaks $D$ into a sequence
of slices.   By general position this is an isotopy apart from a finite
number of critical slices which are maxima, minima and saddles of the 
sheets, maxima and minima of the double curves
and triple points.   Corresponding to these are the diagram moves
listed below:

\rk{Diagram moves}

{\bf BD}\quad  Births and deaths of little circles: $D\Leftrightarrow D\cup O$.

{\bf Br}\quad  Bridge between arcs with compatible framing:
$$
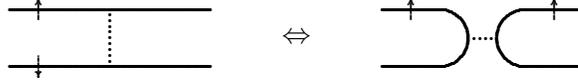
\beginpicture
\setcoordinatesystem units <.15truein,.15truein> point at 0 0
\setplotarea x from 0 to 20, y from -1 to 1
\arrow <2pt> [0.3, 1] from 14 0.65 to 14 1.35
\arrow <2pt> [0.3, 1] from 19 0.65 to 19 1.35
\arrow <2pt> [0.3, 1] from 1 0.65 to 1 1.35
\arrow <2pt> [0.3, 1] from 1 -0.65 to 1 -1.35
\put {$\Leftrightarrow$} at 10 0
\linethickness=0.8pt
\setplotsymbol({\tenrm .})
\plot 0 1 7 1 /
\plot 0 -1 7 -1 /
\plot 13 -1 15 -1 /
\plot 13 1 15 1 /
\plot 18 1 20 1 /
\plot 18 -1 20 -1 /
\circulararc 180 degrees from 15 -1 center at 15 0
\circulararc 180 degrees from 18 1 center at 18 0
\linethickness=0.4pt
\setdots <2pt>
\plot 3.5 -1 3.5 1 /
\plot 16 0 17 0 /
\endpicture
$$
{\bf R2} move

{\bf R3} move

If the cobordism is labelled then the labelling before and after
a move satisfies the following compatibility conditions:

{\bf BD}\quad  Births must correspond to an edge of $C$:
\font\specs=cmtex8  
\def\ds{\hbox{\specs \char'017\kern 0.05em}}
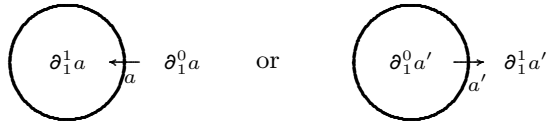
$$
\beginpicture
\setcoordinatesystem units <.3truein,.3truein> point at 0 0
\setplotarea x from -1 to 7, y from -1 to 1
\arrow <2pt> [0.3, 1] from 1.25 0 to 0.75 0
\arrow <2pt> [0.3, 1] from 6.75 0 to 7.25 0
\put {$\scriptstyle a$} [lt] <-1pt,0pt> at 1 -0.2
\put {$\scriptstyle a'$} [lt] <-1pt,0pt> at 7 -0.2
\put {$\scriptstyle \ds^0_1a'$} at 6 0
\put {$\scriptstyle \ds^1_1a'$} at 8 0
\put {$\scriptstyle \ds^1_1a$} at 0 0
\put {$\scriptstyle \ds^0_1a$} at 2 0
\put {\ninerm or} at 3.5 0
\setplotsymbol({\tenrm .})
\circulararc 360 degrees from 1 0 center at 0 0
\circulararc 360 degrees from 5 0 center at 6 0
\endpicture
$$
There is no condition for deaths.

{\bf Br}\quad  Bridges must be between arcs with the same
label.

An {\bf R2} move must involve two double points labelled by
the same square (with opposite orientations):
$$
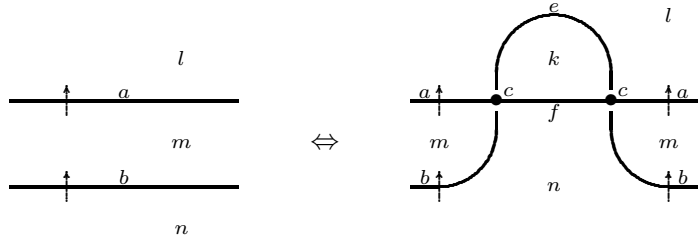
\beginpicture
\setcoordinatesystem units <.3truein,.3truein> point at 0 0
\setplotarea x from -2 to 2, y from -1.25 to 0.25
\arrow <2pt> [0.3, 1] from -1 -0.25 to -1 0.25
\arrow <2pt> [0.3, 1] from -1 -1.75 to -1 -1.25
\linethickness=0.8pt
\putrule from -2 0 to 2 0
\putrule from -2 -1.5 to 2 -1.5
\put {$\scriptstyle a$} [b] <0pt, 1pt> at 0 0
\put {$\scriptstyle b$} [b] <0pt, 1pt> at 0 -1.5
\put {$\scriptstyle l$} at 1 0.75
\put {$\scriptstyle m$} at 1 -0.75
\put {$\scriptstyle n$} at 1 -2.25
\put {$\Leftrightarrow$} at 3.5 -0.75
\setcoordinatesystem units <.3truein,.3truein> point at -7.5 0
\setplotarea x from -2 to 2, y from -1.25 to 1.5
\arrow <2pt> [0.3, 1] from -2 -0.25 to -2 0.25
\arrow <2pt> [0.3, 1] from -2 -1.75 to -2 -1.25
\arrow <2pt> [0.3, 1] from 2 -0.25 to 2 0.25
\arrow <2pt> [0.3, 1] from 2 -1.75 to 2 -1.25
\linethickness=0.8pt
\putrule from -2.5 0 to 2.5 0
\putrule from -2.5 -1.5 to -2 -1.5
\putrule from 2 -1.5 to 2.5 -1.5
\putrule from -1 -0.5 to -1 -0.2
\putrule from 1 -0.5 to 1 -0.2
\putrule from 1 0.2 to 1 0.5
\putrule from -1 0.2 to -1 0.5
\setplotsymbol ({\tenrm .})
\circulararc 90 degrees from -2 -1.5 center at -2 -0.5
\circulararc 90 degrees from 1 -0.5 center at 2 -0.5
\circulararc 180 degrees from 1 0.5 center at 0 0.5
\put {$\scriptstyle a$} [b] <0pt, 1pt> at -2.25 0
\put {$\scriptstyle a$} [b] <0pt, 1pt> at 2.25 0
\put {$\scriptstyle b$} [b] <0pt, 1pt> at 2.25 -1.5
\put {$\scriptstyle b$} [b] <0pt, 1pt> at -2.25 -1.5
\put {$\scriptstyle l$} at 2 1.5
\put {$\scriptstyle m$} at -2 -0.75
\put {$\scriptstyle m$} at 2 -0.75
\put {$\scriptstyle n$} at 0 -1.5
\put {$\bullet$} at -1 0
\put {$\bullet$} at 1 0
\put {$\scriptstyle f$} [t] <0pt, -1pt> at 0 0
\put {$\scriptstyle e$} [b] <0pt, 1pt> at 0 1.5
\put {$\scriptstyle k$} at 0 0.75
\put {$\scriptstyle c$} [lb] <1.5pt, 1.5pt> at -1 0
\put {$\scriptstyle c$} [lb] <1.5pt, 1.5pt> at 1 0
\endpicture
$$

{\leftskip=0.5truein\rightskip=0.5truein\ninepoint
In the figure the following face equalities hold: $n=\d^0_1b=\d^0_1f$,
$m=\d^1_1b=\d^0_1a$, $l=\d^1_1a$, $a=\d^1_1c$, $f=\d^0_1c$,
$e=\d^1_2c$, $b=\d^0_2c$, $k=\d^1_1f=\d^0_1e.$\par}

An {\bf R3} move must correspond to a 3--cube of $C$.

\proc{Proposition}Diagrams labelled in $C$ are cobordant iff
they differ by compatible moves BD, Br, R2, R3 (above).\qed

\rk{Note}If the labelling is by a rack then R2 and R3 moves are always
compatible (this is essentially what the rack laws are designed
for [\FeRo; section 4]), and births can have arbitrary labels.

\sh{Digression on $\pi_2$}

\rk{Remark}Theorem \ref{2diag-class} and the description of cobordisms 
in terms of moves allows us to make calculations of $\pi_2$. See for
example theorem 5.15 for a report of calculations made by this method.
Note at once that the writhe of a diagram is invariant under moves and
hence we always have a map to $\Z$.  The {\sl writhe} is defined (as
usual) to be the number of double points counted with sign --- a
right-hand crossing counting $+1$ and a left-hand one $-1$.

This map to $\Z$ is not always onto as the illustrative 
example below shows.  The writhe has the following
interpretation:  Let $t_C:C\to T$ be the constant map as usual.  Then
writhe is the same as $\pi_2(C)\buildrel (t_C)_* \over \a \pi_2(T)=\pi_2(
\Omega(S^2))=\pi_3(S^2)=\Z$.

\rk{Example}{\sl Calculation of $\pi_2$(torus)}\quad  

As an illustration of the use of diagrams to calculate $\pi_2$ of a
$\sq$--set we calculate $\pi_2$(2--torus).
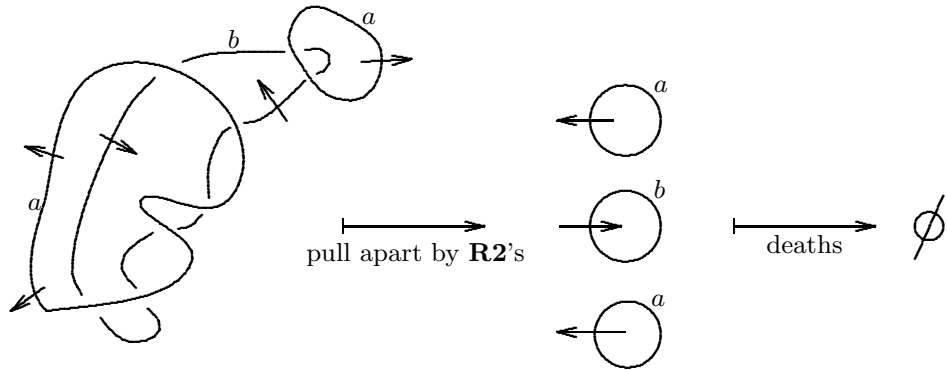
\mbox{\beginpicture
\setcoordinatesystem units < 0.900cm, 0.900cm>
\unitlength= 1.000cm
\linethickness=1pt
\setplotsymbol ({\makebox(0,0)[l]{\tencirc\symbol{'160}}})
\setshadesymbol ({\thinlinefont .})
\setlinear
%
%
\linethickness= 0.500pt
\setplotsymbol ({\thinlinefont .})
\ellipticalarc axes ratio  0.476:0.476  360 degrees 
	from  9.970 24.797 center at  9.493 24.797
%
%
\linethickness= 0.500pt
\setplotsymbol ({\thinlinefont .})
\ellipticalarc axes ratio  0.476:0.476  360 degrees 
	from  9.970 23.368 center at  9.493 23.368
%
%
\linethickness= 0.500pt
\setplotsymbol ({\thinlinefont .})
\ellipticalarc axes ratio  0.444:0.444  360 degrees 
	from  9.970 21.907 center at  9.525 21.907
%
%
\linethickness= 0.500pt
\setplotsymbol ({\thinlinefont .})
\ellipticalarc axes ratio  0.191:0.191  360 degrees 
	from 13.780 23.368 center at 13.589 23.368
%
%
\linethickness= 0.500pt
\setplotsymbol ({\thinlinefont .})
\plot  5.937 25.590  6.604 25.622 /
%
%
\plot  6.353 25.547  6.604 25.622  6.347 25.674 /
%
%
%
\linethickness= 0.500pt
\setplotsymbol ({\thinlinefont .})
\putrule from 10.954 23.368 to 12.859 23.368
%
%
\plot 12.605 23.304 12.859 23.368 12.605 23.431 /
%
%
%
\linethickness= 0.500pt
\setplotsymbol ({\thinlinefont .})
\putrule from  5.683 23.368 to  7.588 23.368
%
%
\plot  7.334 23.304  7.588 23.368  7.334 23.431 /
%
%
%
\linethickness= 0.500pt
\setplotsymbol ({\thinlinefont .})
\putrule from  5.683 23.463 to  5.683 23.273
%
%
\linethickness= 0.500pt
\setplotsymbol ({\thinlinefont .})
%
%
\plot  4.739 25.166  4.540 25.337  4.635 25.092 /
\plot  4.540 25.337  4.921 24.797 /
%
%
\linethickness= 0.500pt
\setplotsymbol ({\thinlinefont .})
%
%
\plot  2.635 24.416  2.889 24.352  2.695 24.528 /
\plot  2.889 24.352  2.413 24.606 /
%
%
\linethickness= 0.500pt
\setplotsymbol ({\thinlinefont .})
\plot  1.905 24.289  1.397 24.448 /
%
%
\plot  1.658 24.432  1.397 24.448  1.620 24.311 /
%
%
%
\linethickness= 0.500pt
\setplotsymbol ({\thinlinefont .})
\plot  1.651 22.543  1.206 22.225 /
%
%
\plot  1.376 22.424  1.206 22.225  1.450 22.321 /
%
%
%
\linethickness= 0.500pt
\setplotsymbol ({\thinlinefont .})
\putrule from  9.335 24.797 to  8.572 24.797
%
%
\plot  8.827 24.860  8.572 24.797  8.827 24.733 /
%
%
%
\linethickness= 0.500pt
\setplotsymbol ({\thinlinefont .})
%
%
\plot  9.176 23.304  9.430 23.368  9.176 23.431 /
\putrule from  9.430 23.368 to  8.604 23.368
%
%
\linethickness= 0.500pt
\setplotsymbol ({\thinlinefont .})
\putrule from  9.493 21.939 to  8.572 21.939
%
%
\plot  8.827 22.003  8.572 21.939  8.827 21.876 /
%
%
%
\linethickness= 0.500pt
\setplotsymbol ({\thinlinefont .})
\putrule from 10.954 23.463 to 10.954 23.273
%
%
\linethickness= 0.500pt
\setplotsymbol ({\thinlinefont .})
\plot 13.780 23.781 13.399 22.924 /
\linethickness= 0.500pt
\setplotsymbol ({\thinlinefont .})
%
%
\plot  1.683 23.812 	 1.700 23.900
	 1.717 23.998
	 1.734 24.106
	 1.753 24.220
	 1.775 24.340
	 1.800 24.464
	 1.829 24.589
	 1.863 24.714
	 1.904 24.838
	 1.953 24.958
	 2.009 25.072
	 2.074 25.179
	 2.149 25.276
	 2.236 25.363
	 2.334 25.436
	 2.445 25.495
	 2.559 25.537
	 2.676 25.565
	 2.795 25.581
	 2.914 25.584
	 3.034 25.575
	 3.152 25.555
	 3.269 25.526
	 3.384 25.487
	 3.496 25.439
	 3.604 25.384
	 3.707 25.321
	 3.805 25.251
	 3.897 25.176
	 3.982 25.096
	 4.059 25.012
	 4.128 24.924
	 4.169 24.860
	 4.208 24.792
	 4.242 24.720
	 4.272 24.645
	 4.297 24.567
	 4.317 24.487
	 4.331 24.405
	 4.339 24.323
	 4.340 24.241
	 4.335 24.159
	 4.321 24.079
	 4.300 24.000
	 4.271 23.924
	 4.232 23.851
	 4.185 23.782
	 4.128 23.717
	 4.062 23.667
	 3.989 23.634
	 3.909 23.618
	 3.823 23.615
	 3.734 23.623
	 3.643 23.639
	 3.551 23.661
	 3.460 23.687
	 3.372 23.714
	 3.287 23.739
	 3.209 23.760
	 3.137 23.774
	 3.021 23.774
	 2.953 23.717
	 2.958 23.633
	 3.032 23.546
	 3.152 23.454
	 3.222 23.406
	 3.296 23.357
	 3.369 23.307
	 3.440 23.256
	 3.563 23.151
	 3.641 23.040
	 3.651 22.924
	 3.609 22.815
	 3.542 22.720
	 3.455 22.637
	 3.349 22.566
	 3.228 22.505
	 3.163 22.478
	 3.095 22.453
	 3.025 22.431
	 2.952 22.410
	 2.878 22.391
	 2.803 22.374
	 2.727 22.359
	 2.650 22.345
	 2.573 22.332
	 2.496 22.320
	 2.419 22.309
	 2.344 22.300
	 2.270 22.291
	 2.197 22.283
	 2.126 22.275
	 2.057 22.268
	 1.992 22.261
	 1.929 22.254
	 1.813 22.240
	 1.715 22.225
	 1.683 22.225
	 1.610 22.311
	 1.554 22.403
	 1.514 22.502
	 1.487 22.605
	 1.473 22.712
	 1.470 22.821
	 1.475 22.932
	 1.489 23.043
	 1.508 23.153
	 1.532 23.262
	 1.559 23.367
	 1.587 23.469
	 1.615 23.565
	 1.642 23.655
	 1.665 23.738
	 1.683 23.812
	/
\linethickness= 0.500pt
\setplotsymbol ({\thinlinefont .})
%
%
\plot  3.143 25.368 	 3.069 25.296
	 3.005 25.233
	 2.906 25.131
	 2.838 25.052
	 2.794 24.987
	 2.737 24.888
	 2.677 24.775
	 2.613 24.651
	 2.581 24.585
	 2.548 24.517
	 2.516 24.446
	 2.483 24.374
	 2.450 24.300
	 2.418 24.225
	 2.386 24.148
	 2.354 24.070
	 2.323 23.992
	 2.293 23.913
	 2.264 23.833
	 2.236 23.753
	 2.209 23.674
	 2.184 23.594
	 2.160 23.515
	 2.137 23.437
	 2.117 23.359
	 2.098 23.282
	 2.081 23.207
	 2.067 23.133
	 2.054 23.061
	 2.044 22.991
	 2.037 22.923
	 2.033 22.857
	 2.032 22.733
	 2.062 22.624
	 2.101 22.547
	 2.159 22.447
	/
\linethickness= 0.500pt
\setplotsymbol ({\thinlinefont .})
%
%
\plot  2.413 22.130 	 2.475 22.054
	 2.530 21.990
	 2.627 21.896
	 2.712 21.839
	 2.794 21.812
	 2.907 21.807
	 3.030 21.827
	 3.138 21.880
	 3.207 21.971
	 3.195 22.096
	 3.140 22.168
	 3.048 22.257
	/
\linethickness= 0.500pt
\setplotsymbol ({\thinlinefont .})
%
%
\plot  2.889 22.543 	 2.793 22.639
	 2.732 22.719
	 2.702 22.790
	 2.699 22.860
	 2.730 22.965
	 2.803 23.062
	 2.927 23.162
	 3.011 23.215
	 3.112 23.273
	/
\linethickness= 0.500pt
\setplotsymbol ({\thinlinefont .})
%
%
\plot  3.461 23.336 	 3.565 23.350
	 3.620 23.368
	 3.682 23.418
	 3.778 23.527
	/
\linethickness= 0.500pt
\setplotsymbol ({\thinlinefont .})
%
%
\plot  3.873 23.749 	 3.868 23.838
	 3.864 23.915
	 3.861 24.038
	 3.864 24.127
	 3.873 24.194
	 3.902 24.298
	 3.948 24.423
	 4.005 24.545
	 4.064 24.638
	 4.159 24.701
	/
\linethickness= 0.500pt
\setplotsymbol ({\thinlinefont .})
%
%
\plot  4.350 24.765 	 4.474 24.778
	 4.540 24.797
	 4.651 24.861
	 4.714 24.907
	 4.784 24.964
	 4.864 25.034
	 4.954 25.118
	 5.057 25.218
	 5.175 25.337
	/
\linethickness= 0.500pt
\setplotsymbol ({\thinlinefont .})
%
%
\plot  3.524 25.590 	 3.642 25.630
	 3.730 25.658
	 3.842 25.686
	 3.914 25.696
	 3.996 25.704
	 4.093 25.710
	 4.206 25.715
	 4.270 25.716
	 4.340 25.718
	 4.415 25.718
	 4.496 25.719
	 4.583 25.719
	 4.678 25.719
	 4.780 25.718
	 4.889 25.718
	/
\linethickness= 0.500pt
\setplotsymbol ({\thinlinefont .})
%
%
\plot  4.953 25.908 	 4.968 25.817
	 5.012 25.720
	 5.066 25.629
	 5.112 25.559
	 5.176 25.459
	 5.262 25.337
	 5.360 25.223
	 5.461 25.146
	 5.572 25.101
	 5.638 25.083
	 5.707 25.071
	 5.778 25.066
	 5.846 25.070
	 5.969 25.114
	 6.092 25.237
	 6.137 25.315
	 6.172 25.401
	 6.195 25.490
	 6.207 25.579
	 6.205 25.667
	 6.191 25.749
	 6.159 25.826
	 6.109 25.897
	 6.044 25.962
	 5.971 26.022
	 5.893 26.075
	 5.816 26.121
	 5.744 26.161
	 5.683 26.194
	 5.586 26.246
	 5.463 26.301
	 5.398 26.322
	 5.332 26.335
	 5.207 26.321
	 5.105 26.251
	 5.026 26.145
	 4.975 26.024
	 4.953 25.908
	/
\linethickness= 0.500pt
\setplotsymbol ({\thinlinefont .})
%
%
\plot  5.175 25.749 	 5.299 25.763
	 5.388 25.758
	 5.493 25.686
	 5.495 25.565
	 5.429 25.463
	 5.334 25.432
	/
%
%
%
\put{\SetFigFont{12}{14.4}{rm}{\hbox{pull apart by {\bf R2}'s}}} [lB] <-7pt, 5pt> at  5.429 22.701
%
%
\put{\SetFigFont{12}{14.4}{rm}{\hbox{deaths}}} [lB] <0pt, 8pt> at 11.366 22.701
%
%
\put{\SetFigFont{12}{14.4}{rm}a} [lB] at  5.906 26.099
%
%
\put{\SetFigFont{12}{14.4}{rm}a} [lB] at  9.842 25.178
%
%
\put{\SetFigFont{12}{14.4}{rm}a} [lB] at  1.397 23.590
%
%
\put{\SetFigFont{12}{14.4}{rm}b} [lB] at  4.096 25.781
%
%
\put{\SetFigFont{12}{14.4}{rm}a} [lB] at  9.811 22.289
%
%
\put{\SetFigFont{12}{14.4}{rm}b} [lB] at  9.842 23.749
\linethickness=0pt
\putrectangle corners at  1.206 26.575 and 13.780 21.463
\endpicture}

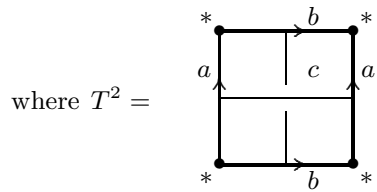
$$
\beginpicture\ninepoint
\setcoordinatesystem units <.35truein,.35truein> point at 0 0
\setplotarea x from 0 to 2, y from 0 to 2
\arrow <4pt> [0.3, 1] from 1 0 to 1.3 0
\arrow <4pt> [0.3, 1] from 1 2 to 1.3 2
\arrow <4pt> [0.3, 1] from 0 1 to 0 1.3
\arrow <4pt> [0.3, 1] from 2 1 to 2 1.3
\putrule from 0 1 to 2 1
\putrule from 1 0 to 1 0.8
\putrule from 1 1.2 to 1 2
\put {where $T^2=$} [r] at -1 1
\linethickness=0.8pt
\putrule from 0 0 to 2 0
\putrule from 0 0 to 0 2
\putrule from 0 2 to 2 2
\putrule from 2 0 to 2 2
\put {$a$} [r] <-2pt, 0pt> at 0 1.4
\put {$a$} [l] <2pt, 0pt> at 2 1.4
\put {$b$} [t] <0pt, -2pt> at 1.4 0
\put {$b$} [b] <0pt, 2pt> at 1.4 2
\put {$c$} at 1.4 1.4
\multiput {$\bullet$} at 0 0 2 0 0 2 2 2 /
\multiput {$*$} at -0.2 -0.2 2.2 -0.2 -0.2 2.2 2.2 2.2 /
\endpicture
$$

\rk{The Whitehead conjecture}

To illustrate the applicability of this
method in general, we now give a translation of the Whitehead conjecture
[\Whitehead] into a conjecture about coloured diagrams. 

Recall that the Whitehead conjecture states that if $K$ is a
subcomplex of the 2--dimensional complex $L$ then $\pi_2(L)=0$ implies
$\pi_2(K)=0$.  By [\JBun; proposition 1.2] it is sufficient to
establish the conjecture for 2--dimensional $\sq$--sets and subsets.

We consider plane diagrams up to moves Br, BD and R2. No R3's are
allowed because the dimension is at most 2.

A {\sl colouring} of a diagram is a colouring of arcs, regions and
crossings.  The outside (infinite) region is always coloured white say.  A
{\sl colour scheme} is a list of allowable colours (three lists one each for
regions, arcs and crossings) and rules about neighbouring colours.
The rules prescribe the two neighbouring colours for a given colour
on an edge and all the neighbouring colours for a given colour on a
crossing.

Moves are {\sl allowable} if Br and BD moves respect the colour scheme and
R2 moves involve two crossings with the same colour (but opposite
orientation).

A diagram is {\sl reducible} if it can be changed to the empty (all white)
diagram by allowable moves.

Using proposition 4.9, 
the Whitehead conjecture now has the following equivalent statement:

\proc{Conjecture}   Suppose that for a given colour scheme 
all diagrams are reducible, then any diagram can be reduced without
using any colours not already used in the diagram.\rm

\sh{Classification of links using moves}

Let $X$ be a rack; recall that $\F(3,X)$ is the group of cobordism
classes of links in $\re^3$ with representation in $X$; recall also
that $\F(3,X)\cong\D(2,X)\cong\pi_2(BX)$ as a special case of 
\ref{diag-class} and \ref{diag-link}.  It is possible to analyse 
cobordisms combinatorially in this  dimension   and this provides   an
alternative proof of the first of these two isomorphisms.

\proc{Theorem}\key{spec-case}{\rm (Special case of \ref{diag-link})}\qua
$\F(3,X)\cong\D(2,X)$

\prf Any diagram labelled in $X$ lifts to a framed link with representation
in $X$ and a cobordism of diagrams to a cobordism of links.  Thus we
have a homomorphism:
$$\Psi\co\D(2,X)\a\F(3,X)$$ $\Psi$ is surjective since any framed link
is the lift of a diagram (unique up to moves R11, R2, R3).  To see
that $\Psi$ is injective we have to show that labelled diagrams whose
lifts are cobordant are themselves cobordant (as diagrams) or
equivalently (using proposition 4.2) that they differ by moves BD, Br,
R2 and R3.  We shall show this by analysing the cobordism.  Since the
representation (or labelling) in $X$ plays no real r\^ole in the proof
it will henceforth be suppressed.

Now the cobordism is a framed 3--manifold in $\re^3\times I$ and
by slicing by parallel $\re^3$'s we obtain (using general position)
a sequence of framed links with
the following critical stages (where a slice contains a critical point
of the projection to $I$):  Bridge moves, births and
deaths of small circles.  Moreover by rotating small neighbourhoods
of these critical points (if necessary) we can assume that near the 
critical point nearby slices are lifts of diagrams;  ie that the R1
moves needed to correct the framing are away from the critical points.

Thus we can choose diagrams near these critical stages  whose lifts
are isotopic to the nearby slices and which differ by diagram moves BD or Br.
Now away from critical stages each link in the sequence corresponds
to a diagram unique up to moves R2, R3, R11 and combining the 
two sets of moves, we see that the cobordism corresponds to diagram
moves BD, Br, R2, R3, R11. But the following sequence of
pictures shows how to achieve an R11 as a combination of a
bridge move, an R2 and a death.  The rest of the theorem
is clear.\qed 
\mbox{\def\SetFigFont#1#2#3{\ninepoint\bf}
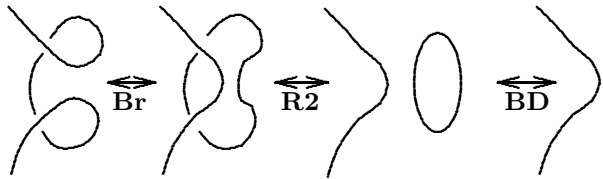
\beginpicture
\setcoordinatesystem units < 1.000cm, 1.000cm>
\unitlength= 1.000cm
\linethickness=1pt
\setplotsymbol ({\makebox(0,0)[l]{\tencirc\symbol{'160}}})
\setshadesymbol ({\thinlinefont .})
\setlinear
%
%
\linethickness= 0.500pt
\setplotsymbol ({\thinlinefont .})
\ellipticalarc axes ratio  0.286:0.603  360 degrees 
	from  5.937 25.273 center at  5.652 25.273
%
%
\linethickness= 0.500pt
\setplotsymbol ({\thinlinefont .})
%
%
\plot  1.905 25.337  1.651 25.273  1.905 25.209 /
\putrule from  1.651 25.273 to  2.223 25.273
%
%
\plot  1.969 25.209  2.223 25.273  1.969 25.337 /
%
%
%
\linethickness= 0.500pt
\setplotsymbol ({\thinlinefont .})
%
%
\plot  3.937 25.337  3.683 25.273  3.937 25.209 /
\putrule from  3.683 25.273 to  4.318 25.273
%
%
\plot  4.064 25.209  4.318 25.273  4.064 25.337 /
%
%
%
\linethickness= 0.500pt
\setplotsymbol ({\thinlinefont .})
%
%
\plot  6.636 25.337  6.382 25.273  6.636 25.209 /
\putrule from  6.382 25.273 to  7.080 25.273
%
%
\plot  6.826 25.209  7.080 25.273  6.826 25.337 /
\linethickness= 0.500pt
\setplotsymbol ({\thinlinefont .})
%
%
\plot  0.857 24.670 	 0.915 24.589
	 0.962 24.533
	 1.048 24.479
	 1.143 24.469
	 1.249 24.485
	 1.349 24.521
	 1.429 24.575
	 1.479 24.638
	 1.517 24.720
	 1.535 24.808
	 1.524 24.892
	 1.478 24.967
	 1.406 25.026
	 1.322 25.066
	 1.238 25.082
	 1.160 25.066
	 1.080 25.022
	 1.009 24.968
	 0.953 24.924
	 0.884 24.866
	 0.802 24.788
	 0.724 24.707
	 0.667 24.638
	 0.597 24.514
	 0.564 24.442
	 0.540 24.384
	 0.515 24.307
	 0.498 24.246
	 0.476 24.162
	/
\linethickness= 0.500pt
\setplotsymbol ({\thinlinefont .})
%
%
\plot  0.444 26.194 	 0.531 26.095
	 0.606 26.010
	 0.671 25.938
	 0.727 25.877
	 0.818 25.783
	 0.889 25.718
	 0.961 25.648
	 1.053 25.562
	 1.157 25.491
	 1.270 25.463
	 1.374 25.484
	 1.473 25.534
	 1.550 25.613
	 1.587 25.718
	 1.574 25.842
	 1.515 25.949
	 1.420 26.027
	 1.302 26.067
	 1.216 26.063
	 1.134 26.023
	 1.049 25.942
	 1.003 25.883
	 0.953 25.813
	/
\linethickness= 0.500pt
\setplotsymbol ({\thinlinefont .})
%
%
\plot  0.857 25.622 	 0.792 25.543
	 0.762 25.495
	 0.739 25.426
	 0.718 25.339
	 0.704 25.250
	 0.699 25.178
	 0.713 25.076
	 0.732 24.997
	 0.762 24.892
	/
\linethickness= 0.500pt
\setplotsymbol ({\thinlinefont .})
%
%
\plot  2.730 25.622 	 2.665 25.543
	 2.635 25.495
	 2.612 25.426
	 2.591 25.339
	 2.577 25.250
	 2.572 25.178
	 2.586 25.076
	 2.606 24.997
	 2.635 24.892
	/
\linethickness= 0.500pt
\setplotsymbol ({\thinlinefont .})
%
%
\plot  2.762 24.670 	 2.835 24.590
	 2.893 24.534
	 2.985 24.479
	 3.051 24.469
	 3.130 24.471
	 3.207 24.486
	 3.270 24.511
	 3.367 24.610
	 3.429 24.733
	 3.437 24.792
	 3.436 24.862
	 3.397 24.987
	 3.270 25.051
	 3.246 25.108
	 3.236 25.178
	 3.238 25.305
	 3.258 25.418
	 3.302 25.527
	 3.399 25.586
	 3.493 25.654
	 3.512 25.712
	 3.518 25.781
	 3.493 25.908
	 3.384 26.028
	 3.311 26.074
	 3.238 26.099
	 3.137 26.080
	 3.048 26.035
	 2.975 25.974
	 2.924 25.920
	 2.857 25.845
	/
\linethickness= 0.500pt
\setplotsymbol ({\thinlinefont .})
%
%
\plot  2.318 24.162 	 2.341 24.233
	 2.361 24.295
	 2.394 24.392
	 2.445 24.511
	 2.494 24.593
	 2.563 24.692
	 2.636 24.788
	 2.699 24.860
	 2.754 24.904
	 2.827 24.952
	 2.953 25.051
	 3.014 25.155
	 3.048 25.273
	 3.039 25.336
	 3.015 25.407
	 2.953 25.527
	 2.844 25.626
	 2.730 25.718
	 2.634 25.828
	 2.540 25.940
	 2.453 26.029
	 2.383 26.099
	 2.286 26.194
	/
\linethickness= 0.500pt
\setplotsymbol ({\thinlinefont .})
%
%
\plot  4.318 24.130 	 4.341 24.201
	 4.361 24.263
	 4.394 24.360
	 4.445 24.479
	 4.495 24.561
	 4.563 24.660
	 4.636 24.756
	 4.699 24.829
	 4.754 24.872
	 4.827 24.921
	 4.953 25.019
	 5.014 25.123
	 5.048 25.241
	 5.039 25.304
	 5.016 25.375
	 4.953 25.495
	 4.844 25.595
	 4.731 25.686
	 4.635 25.796
	 4.540 25.908
	 4.453 25.997
	 4.383 26.067
	 4.286 26.162
	/
\linethickness= 0.500pt
\setplotsymbol ({\thinlinefont .})
%
%
\plot  6.890 24.162 	 6.913 24.233
	 6.933 24.295
	 6.966 24.392
	 7.017 24.511
	 7.066 24.593
	 7.135 24.692
	 7.208 24.788
	 7.271 24.860
	 7.326 24.904
	 7.399 24.952
	 7.525 25.051
	 7.586 25.155
	 7.620 25.273
	 7.611 25.336
	 7.587 25.407
	 7.525 25.527
	 7.416 25.626
	 7.303 25.718
	 7.206 25.828
	 7.112 25.940
	 7.025 26.029
	 6.955 26.099
	 6.858 26.194
	/
%
%
\put{\SetFigFont{12}{14.4}{rm}Br} [lB] <-4pt, 4pt> at  1.841 24.797
%
%
\put{\SetFigFont{12}{14.4}{rm}R2} [lB] <0pt, 4pt> at  3.747 24.797
%
%
\put{\SetFigFont{12}{14.4}{rm}BD} [lB] <-4pt, 4pt> at  6.604 24.797
\linethickness=0pt
\putrectangle corners at  0.444 26.194 and  7.620 24.130
\endpicture}

\sh{Virtual links}

Virtual links have been introduced by Kauffman [\Kauf] and studied by
several authors including Carter--Saito--Kamada [\CKS], Fenn--Jordan--Kauffman 
[\FJK], Kamada--Kamada
[\KamKam] and Kuperberg [\Kup].  Here we shall show that the second
homology group of the rack space $BX$ classifies framed virtual links
with representation in a rack $X$, up to cobordism.

\rk{Oriented Gauss codes and oriented crossing graphs}

There are several equivalent definitions of a virtual link.  Kauffman
[\Kauf] introduced the subject and defined a virtual link as an
equivalence class of oriented Gauss codes up to Riedemeister moves
(R1, R2 and R3).  An oriented Gauss code is the same as a 4--valent
graph such that each vertex can be indentified with a standard
crossing in the plane up to rotation through $\pi$.  In other words,
we know which of the edges at the vertex are the over-crossing arcs
and which are the undercrossing arcs and we also have a cyclic
ordering of the arcs at the vertex --- over, under, over, under.  We
call such a 4--valent graph an {\sl oriented crossing graph}.  So a
virtual link is an equivalence class of oriented crossing graphs under
Riedemeister moves.

\rk{Virtual link diagrams}

Now an oriented crossing graph can be immersed in the plane with the
vertices forming crossings with the correct orientation.  Such an
immersion is well-defined up to changing the immersion on edges and
this leads to the more usual definition of virtual links in terms of
diagrams.  The crossings which come from the vertices of the graph are
the real crossings and the ones which come from crossings of the
immersed edges are the virtual crossings.  The result is a {\sl
virtual link diagram}.  The equivalence relation on virtual link
diagrams is generated by Reidemeister moves on real crossings together
with the ability to move an arc containing only virtual crossings to
any other position with the same endpoints.  This latter move can be
replaced by a set of {\sl extended Reidemeister moves\/} ---
Reidemeister moves R1, R2 and R3 for virtual crossings and one mixed
move, namely an R3 move with two virtual and one real crossing.  See
Kauffman [\Kauf; figures 2 and 3].

\rk{Framed virtual links}

We need to extend these definitions to {\sl framed virtual links}.  To
frame a virtual link, we orient the components and use the blackboard
framing convention.  Moreover in the above definitions we replace the
R1 move for real crossings with the double Reidemeister 1--move (the
R11 move).  Thus a framed virtual link is an equivalence class of
oriented crossing graphs, with edges oriented compatibly with
crossings, under R2, R3 and R11.  Equivalently it is an equivalence
class of oriented virtual link diagrams under moves R2, R3 and R11 for
real crossings plus the extended moves, {\it including R1 for virtual
crossings\/}. 

\rk{Stable equivalence}

We need to interpret (framed or unframed) virtual links as equivalence
classes of genuine (framed or unframed) links in oriented surfaces cross
$I$.  

\rk{Definitions}Let $\Sigma$ be an oriented surface and let $L$ be a link 
in $\Sigma\times I$.  Suppose that $D_1,D_2\subset \Sigma$ are discs
disjoint from the projection of $L$.  If we add an oriented handle to
$\Sigma$ with feet at $D_1,D_2$ to form $\Sigma'$ then we say that
$L\subset\Sigma'\times I$ is a {\sl stabilization\/} of
$L\subset\Sigma\times I$.

Links $L_1\subset\Sigma_1\times I$ and $L_2\subset\Sigma_2\times I$
are {\sl stably equivalent\/} if they differ by a sequence of
stabilizations and their inverses.\medskip  

The following result is ``well-known''.  It was suggested by Kauffman
in [\Kauf] and can be deduced from [\CKS; proposition 3.4].  However
the proof of this last result spreads over three papers [\CKS, \KamKam,
\Kauf] so it seems worthwhile to include here a short direct proof.

\proc{Theorem}\key{stab-equiv}Virtual links (respectively framed virtual 
links) are in bijective correspondence with stable equivalence classes
of links (respectively framed links) in oriented surfaces cross $I$.

\prf  Let $\VL$ denote the set of virtual links and let $\SL$ denote stable
equivalence classes of links in oriented closed surfaces cross $I$.

There is a function $\phi\co \SL\to \VL$ given as follows: Let
$L\subset\Sigma\times I$ be a link and project it to form a diagram in
$\Sigma$.  The diagram determines an oriented crossing graph and hence
a virtual link.  It is clear that the result is unaltered by
stabilization.  An isotopy can be replaced by Reidemeister moves
which correspond to Reidemeister moves on the oriented crossing graph.
Hence $\phi$ is well-defined.

There is another function $\psi\co \VL\to \SL$ given as follows.  Use
the definition of a virtual link as an equivalence class of oriented
crossing graphs.  Let $\Gamma$ be such a graph; consider an immersion
of $\Gamma$ in the plane with the vertices forming crossings with the
correct orientation.  Let $N$ be an induced regular neighbourhood.  It
is easy to see that $N$ depends only on the original graph.  Form an
oriented surface $\Sigma$ by capping the boundary circles of $N$ with
discs.  We obtain a link diagram in $\Sigma$ and hence a link in
$\Sigma\times I$.  The result is unchanged by an R1 or R3 move on the
original graph whilst an R2 may add (or delete) a handle disjoint from
the diagram.  Thus $\psi$ is also well-defined.

It is clear that $\phi\circ\psi$ is the identity on $\VL$.  To see
that $\psi\circ\phi$ is the identity on $\SL$ observe that the result
of applying $\psi\circ\phi$ to a link in $\Sigma\times I$ with diagram
$D$ is the same as surgering along the circles which form the boundary of
a regular neighbourhood of $D$.   The surgery is a stable equivalence.

For the framed case, the proof is exactly the same with R1 moves replaced
by R11 moves throughout.\endprf

\rk{Remark}Kuperberg [\Kup] has proved a far stronger result: a virtual
link has a unique {\sl irreducible\/} respresentation as a link in an oriented
surface cross $I$, where irreducible means that no destabilizations
are possible.

\rk{The fundamental rack of a virtual link}

There is a notion of a fundamental rack of a framed virtual link
obtained from any diagram by labelling arcs with generators (ignoring
all virtual crossings) and reading a relation at real crossings by the
usual rule (\ref{cross-rel}).  Since this is the same as reading the
fundamental rack of any corresponding diagram in an oriented surface,
it follows from lemma \ref{reduced-ident} that this rack coincides
with the reduced fundamental rack of any corresponding link in an
oriented surface cross $I$.

Thus we have a good notion of a representation in a rack $X$ for a
framed virtual link, namely a homomorphism of fundamental rack in $X$.

\rk{Cobordism of virtual links}

Virtual links are {\sl cobordant} if they differ by the allowed moves
plus the two cobordism moves BD (birth--death) and Br (bridge)
introduced earlier.  Note that R11 moves can be obtained from R2, Br
and BD moves and do not need to be included here.  Combining the
proofs of \ref{spec-case} and \ref{stab-equiv} we see that this
corresponds to double cobordism of any corresponding link in surface
cross $I$.  (Notice that cobordism of surfaces is generated by
stabilization.)

There are similar definitions of cobordism of framed links and of
links with representation in a rack $X$.

Cobordism classes form an abelian group under disjoint union and we
denote the group of cobordism classes of framed virtual links with
representation in a rack $X$ by $\VL(X)$.

\proc{Theorem}{\rm (Classification of virtual links)}

There is a natural isomorphism between $\VL(X)$ and the second
homology group of the rack space $H_2(BX)$.

\prf Theorem \ref{double-class} shows that there is a natural bijection
between the set of double cobordism classes of links in oriented
surfaces cross $I$ with representation in $X$ and the oriented bordism
group $\Omega_2(BX)$.  But $H_2(BX)=\Omega_2(BX)$ and the result
follows.
\qed

\sh{Calculations}

In section 5 we report on calculations of $H_2$ of rack spaces, in
particular Greene's results in theorem 5.16.  Using the last result
any of these calculations implies results about virtual links.

Below are some samples.  By {\sl $n$--colouring} we mean a
representation in the dihedral rack $D_n := \{0,1,\ldots,n-1\}$ with
$i^j=2j-i \ \hbox{mod}\ n$ for all $i,j$. The {\sl writhe} of a framed
virtual link is the number of crossings counted algebraically (the
framing gives signs to crossings in the usual way).  Writhe is a
cobordism invariant and can be interpreted as the element of
$H_2(T)\cong\Z$ determined by the link.  Here $T$ is the trivial
$\sq$--set and is the rack space of the trivial rack (ie no
labelling).

\medskip
(1)\qua {\sl Any odd coloured or uncoloured virtual link is
cobordant to any other with the same writhe.}

\medskip
(2)\qua {\sl $\VL(D_{2n})$ for $n>0$ is at least $\Z^4$, thus there
are at least three infinite families of even coloured cobordism
classes with any given writhe.  Furthermore $\VL(D_{4})$ has
2--torsion as well.}

\section{The algebraic topology of rack spaces}

In this section we turn to invariants of rack spaces. These are
important because any invariant of rack spaces  automatically becomes
a knot or link invariant by calculating it for the rack space of the
fundamental rack of the knot or link. Moreover, further invariants
can be derived by considering representations of the fundamental rack
in other racks, for example, finite racks. 

All the invariants considered in this section are homotopy type
invariants; however the homotopy type of the rack space is not a
complete link invariant for links in $S^3$ even when the rack itself
is (see the remarks on theorem 5.4 below).  In this context it is
worth reiterating that the combinatorial structure of the rack space
is equivalent to the rack itself and therefore it is valuable to
construct combinatorial invariants of a $\sq$--set which are not
homotopy type invariants. The main new invariants introduced in
[\JBun] (the James--Hopf invariants), and also the associated
generalised cohomology theories, are combinatorial invariants of this
type.

Here we start by identifying the fundamental group of rack spaces and
proving that they are simple. We then turn to calculations of
homotopy groups. We calculate all the homotopy groups of the rack
space of an irreducible link in a 3--manifold and the second homotopy
group of the rack space of any link in $S^3$. We determine the
homotopy type of the rack space of an irreducible link in an
irreducible 3--manifold with infinite fundamental group. We also
determine the homotopy type of the rack space of a free rack and of
the trivial rack with $n$ elements. (Note that Wiest
[\Wiest] has also determined the homotopy type of the rack  space of an
irreducible link in a general 3--manifold.)

We conclude with some results on homology groups and a review of
results of Flower [\Flower] and Greene [\Greene].

\sh{Fundamental group}

The fundamental groupoid of a $\sq$-set is discussed in [\Trunks].  We
repeat the computation of the fundamental group of the rack space.

Recall from [\FeRo] that the {\sl associated group} $As(X)$ to a rack
$X$ is the group generated by the elements of $X$ subject to the 
relations $a^b=b\inv ab$.

\proc{Proposition}The fundamental group $\pi_1(BX)$ of the rack space
of a rack $X$ is isomorphic to the associated group $As(X)$ of $X$.

\prf Recall that the rack space $BX$ has a single vertex and edges in
bijection with the elements of $X$ which therefore generate
$\pi_1(BX)$. Moreover the relations given by the squares of $BX$ are
$a^b=b\inv ab$ for $a,b\in X$. The result now follows from the
definition of $As(X)$. \qed 

\sh{Simplicity of the rack space}

We next prove that $BX$ is a simple space for any rack $X$.

\proc{Proposition}Let $X$ be any rack then the action of the
fundamental group  $\pi_1(BX)$ on $\pi_n(BX)$ is trivial. 

\prf Let $\beta\in\pi_n(BX)$, and let $\alpha\in\pi_1(BX)$ be a
generator corresponding to $x\in X$ as above. We may represent
$\beta$ by a labelled diagram $D$ in $\re^n$. Then $\alpha\cdot\beta$
is represented by the diagram which comprises a  framed standard
sphere in $\re^n$ labelled by $x$ and containing $D$ in its interior.
But the following diagram cobordism shows that the two diagrams are
equivalent. Pull the sphere under $D$ to one side (without changing
any labels on $D$) and then eliminate it.\qed

\rk{Notation}Let $X$ be a rack. Recall from section 3 that 
$\D(n,X)$, $\F(n+1,X)$ denote the group of cobordism classes of diagrams
in $\re^n$ labelled by $X$ and the group of framed cobordism classes
of framed codimension 2 links in $\re^{n+1}$  with representation in
$X$, respectively. 

We shall use the notation $[D,\lambda]$ for an element of $\D(n,X)$,
where $D$ denotes a diagram and $\lambda$ a labelling and we shall
use the notation $[L,F,\rho]$ for an element of $\F(n+1,X)$ where $L$
is a codimension 2 link  (an $(n-1)$--manifold) in $\re^{n+1}$, $F$
is a framing of $L$ and $\rho$ a representation in $X$, ie a
homomorphism of the fundamental rack of $(L,F)$ to $X$. 

If $x\in X$ is an element of $X$, then we denote by $\lambda^x$ the
labelling obtained from $\lambda$ by operating on all labels by $x$.
That this is also a labelling follows from the rack law. Similarly we
denote by $\rho^x$ the representation obtained by composing $\rho$
with the automorphism $a\mapsto a^x$ of $X$.

\proc{Corollary}Let $X$ be a rack. Let $[D,\lambda]\in\D(n,X)$ and
let $[L,F,\rho]\in\F(n+1,X)$ and $x\in X$. Then
$[D,\lambda]=[D,\lambda^x]$, and $[L,F,\rho]=[L,F,\rho^x]$. \prf To
see that $[D,\lambda]=[D,\lambda^x]$ simply pull the sphere  over
instead of under in the proof of 5.2. The other result now follows
from the isomorphism of the two groups (\ref{diag-link}).\qed

\rk{Remark}Notice that in the case that $X$ is the fundamental rack
of a link $L$ in $S^n$ then the action of $x$ is given by $a^x=a\cdot\d
x$ where $\d$ is the augmentation to $\pi_1(L):= \pi_1(S^3-L)$, 
hence the corollary also implies invariance of $\D(n,X)$ and
$\F(n+1,X)$ under the action of $\pi_1(L)$ in this case.

\sh{Homotopy groups of the rack space of an irreducible link}

Let $L$ be a framed link in a 3--manifold $M^3$. We say that $L$ is 
{\sl irreducible} if each embedded 2--sphere in $M^3-L$ bounds a
3--ball (ie $M-L$ is an irreducible 3--manifold). We say that $L$ is
{\sl trivial} if $M=S^3$ and $L$ is equivalent to the unknot  in
$S^3$ with zero framing. 

Let $\Lambda$ be a framed submanifold of $\re^{n+1}$ with framing
$F\co \Lambda\times D^2\to\re^{n+1}$. Write $N(\Lambda)={\rm im}(F)$
for the tubular neighbourhood of $\Lambda$ and
$\Lambda^c=\bar{\re^{n+1}-N(\Lambda)}$ for the closure of the
complement. Also write $\Lambda^+=F(\Lambda\times\{(1,0)\})$ for the
{\sl parallel manifold} to $\Lambda$.

\proc{Theorem}Let $M^3$ be a $3$--manifold and  let $L$ be a framed
non-trivial irreducible link in $M^3$.  Let $\Gamma(L)$ denote the
fundamental rack of $L$. Then for $n>1$, 
$$ \pi_n(B\Gamma(L))\cong\pi_{n+1}(M^3).$$ 

\prf Recall from theorem \ref{link-class} that  $\pi_n(B\Gamma(L))\cong
\F(n+1,\Gamma(L))$. 

Define the homomorphism $$\phi\co \pi_{n+1}(M^3) \to
\F(n+1,\Gamma(L))$$ as follows. Let $\alpha\in\pi_{n+1}(M^3)$ be
represented by a map  $f\co {\Bbb R}^{n+1}\to M^3$ which is constant
outside a compact set. Homotope $f$ to be transverse to $L$. Then
$\Lambda=f^{-1}L$ is a framed submanifold with framing $F$ such that
$f$ is compatible  with the framings. Let
$\rho\co\Gamma(\Lambda)\to\Gamma(L)$ be the induced homomorphism of
racks. We define $\phi(\alpha)=[\Lambda,F,\rho]$. If $f$ and $g$ are
both transverse representatives of $\alpha$ then a homotopy between
$f$ and $g$ can also be made transverse producing a bordism and we
see that $\phi$ is well defined.

\bf$\phi$ is surjective\rm 

For each component $\Lambda^+_i$ of $\Lambda^+$ let $L^+_i$ be the
corresponding component of $L$ which is in the image of
$\Lambda^+_i$. Notice that the $L^+_i$'s may not be distinct. For
each $i$ choose an embedded path $p_i$ from $\Lambda^+_i$ to the base
point representing an element of $\Gamma(\Lambda)$ so that the chosen
paths only meet at the base point. Similarly choose a path $q_i$ for
$L^+_i$, where $\rho[p_i]=[q_i]$. Corresponding to our choices we
have longitudinal and meridinal subgroups $\Lambda_i(l),\Lambda_i(m)$
of $\pi_1( \Lambda^c)$ and subgroups $L_i(l),L_i(m)$ of $\pi_1(L^c)$.
Now $As(\Gamma(\Lambda))\cong\pi_1(\Lambda^c)$ and 
$As(\Gamma(L))\cong\pi_2(M,L^c)$.  Then
$\rho\co\Gamma(\Lambda)\to\Gamma(L)$ induces a homomorphism 
$\rho_1\co\pi_1(\Lambda^c)\to  \pi_1(L^c)$ by composing with the
boundary map in the homotopy exact sequence  of the pair $(M,L^c)$.
See [\FeRo; page 360].\fnote{\ninepoint It has been pointed out by
Wiest that the proof of the result used here [\FeRo; proposition
3.2] contains a misleading statement.  To be precise the
statement made at the top of page 361 in [\FeRo] is open to 
misinterpretation.  The paths can also be
varied by an isotopy which moves one little disc around another --
essentially a pure braid automorphism. This can be realised by two
interchanges and is implicit in the subsequent text.} The
homomorphism has restrictions  $\rho_i(l)\co \Lambda_i(l)\to L_i(l)$
and  $\rho_i(m)\co \Lambda_i(m)\to L_i(m)$. 

We claim that the hypotheses on $L$ imply that the group $L_i(l)$  is
infinite cyclic. To see this, suppose not. Then some multiple of the
longitude is null homotopic in $L^c$. By the loop theorem there is a
closed essential curve in the neighbourhood of the longitude $L_i^+$
which bounds a disc in $L^c$. But this neighbourhood is an annulus
and the only possibility is that $L_i^+$ itself bounds a disc. By
irreducibility, we then see that $L=L_i$ and is trivial,
contradicting our hypotheses.

The space $L^+_i$ together with its embedded `tail' $q_i$ is an
Eilenberg-MacLane space so there is a unique map, up to homotopy,
$f_0\co \Lambda^+_i\cup p_i(I)\to L^+_i\cup q_i(I)$ inducing the
homomorphism $\rho_i(l)$. Further we can assume, for later
convenience, $f_0\Lambda^+_i\subset L^+_i$ and $f_0 p_i=q_i$. Since
$\rho_1$ is induced by sending meridians to meridians there is a
unique extension $f_1\co N(\Lambda)\to N(L)$ which preserves framing.
We can finally extend $f$ over $\Lambda^c$ since $L^c$  is also an
Eilenberg-MacLane space. 

\bf $\phi$ is injective \rm 

First we observe that the map $f$ constructed above is unique up to
based homotopy. To see this suppose alternative choices $p'_i$ and
$q'_i$ are made in place of $p_i$ and $q_i$ respectively and assume
end points agree. Consider $f'_0\co \Lambda^+_i\cup p'_i(I)\to
L^+_i\cup q'_i(L)$. We can assume the set of $p_i$'s together with
the $p'_i$'s do not meet on interiors. Now suppose $f_0\circ
(\bar{p_i}\cdot r\cdot p_i) \simeq \bar{q_i}\cdot l^{n_i}\cdot q_i$
where $r$ is some loop in $\Lambda^+_i$ and $l$ is `once round'
$L^+_i$.
Notice $p'_i \simeq p_i\cdot(\bar{p_i}\cdot p'_i)$. Then from the
fact that a homomorphism of racks preserves the action of the
associated groups we have, after an easy calculation,
$\rho_1[\bar{p'}_i\cdot r\cdot p'_i]= 
[\bar{q'}_i\cdot l^{n_i}\cdot q'_i]$. It follows that $f_0$ and $f'_0$
can be taken to agree on $\Lambda^+_i$ and $f_1$ and $f'_1$ agree on
$N(\Lambda)$.

Now in completing the constructing of $f$ consider the construction
over the $1$-skeleton. If we have $fp_i=q_i$ then essentially the
same argument shows we can assume $fp'_i=q'_i$. We are now ready to
prove $\phi$ is injective. Suppose $(W,F,\rho)$ is a bordism
between $(W_0,F_0,\rho_0)$ and $(W_1,F_1,\rho_1)$ and the latter
are associated with maps $g_0$ and $g_1$ respectively. Now repeat
the proof that $\phi$ is onto but this time using $S^{n+1}\times I$
and $W$ and $\rho$ in place of $S^{n+1}$ and $\Lambda$ and $\rho$.
The resulting bordism then must give $g_0$ and $g_1$ on the two ends
by the above observations.\qed

\sh{Remarks on and consequences of theorem 5.4}

It is worth pointing out that irreducibility of the {\sl link}
$L$ in $M$ does not imply irreducibility of the {\sl3--manifold}
$M$. Indeed any connected 3--manifold contains an irreducible
link (see [\FeRo; page 380 remark (2)]). The higher homotopy
groups of a general (non-irreducible) 3--manifold can be very
complicated. Each separating 2--sphere potentially determines
a copy of $\pi_*(S^2)$ which is then subject to action by $\pi_1$
of the manifold.

The theorem therefore implies that the higher homotopy
groups of rack spaces of irreducible links can in general be
complicated. In the case that the 3--manifold {\sl is} irreducible
its higher homotopy groups are either all zero (the case when the
fundamental group of $M$ is infinite) or coincide with the homotopy groups
of the 3--sphere. Turning first to infinite fundamental group
case, we deduce:

\proc{Corollary}Let $L$ be a framed non-trivial irreducible link in 
an irreducible 3--manifold $M$ with infinite fundamental group
and let $\Gamma(L)$ denote the fundamental rack of $L$.
Then $B\Gamma(L)$ is a $K(\pi,1)$ where $\pi$ is the kernel
of $\pi_1(M-L)\to \pi_1(M)$.

\prf By the theorem and the remarks above, $B\Gamma(L)$ is a $K(\pi,1)$
where $\pi$ is the associated group of $\Gamma(L)$ by proposition
5.1. But since $\pi_2(M)=0$, [\FeRo; proposition 3.2] implies that
the associated group is the kernel of $\pi_1(M-L)\to \pi_1(M)$. \qed

In the case when $M$ is irreducible and has finite fundamental
group (eg when $M=S^3$), then (as remarked earlier) the theorem
implies that the higher homotopy groups of $B\Gamma(L)$ coincide
with the higher homotopy groups of $S^3$ (with an index shift of 1).
Thus the theorem gives a plentiful supply of new geometric
interpretations for these groups. We now describe one such
interpretation which does not need the concept of fundamental
rack for its statement.

We need to define the writhe of a framed link in a higher
dimensional sphere. Let $M^{n-1}$ be a framed submanifold
of $S^{n+1}$. Let $M^+$ denote the parallel manifold to $M$
as usual and let $\phi\co M^c\to S^1$ be the map defined by any
Seifert bounding manifold for the codimension 2 submanifold $M$. Then
the composition $M\to M^+{\buildrel\phi|M^+\over\longrightarrow} S^1$
defines an element of $H^1(M)$ called the {\sl writhe}.
There is a similar notion for the writhe of a cobordism and a
cobordism which preserves writhe.

\proc{Corollary}Fix an integer $w>0$ then $\pi_{n+1}(S^3)$ is 
isomorphic to the set of equivalence classes of framed $(n-1)$--manifolds
embedded in $S^{n+1}$ with writhe divisible by $w$, under framed
cobordism also with writhe divisible by $w$.

\prf Let $U_w$ denote the unknot in $S^3$ with framing $w$. Then
$U_w$ is irreducible, non-trivial and the theorem applies. But
writhe divisible by $w$ is the same as having a representation
in the cyclic rack of order $w$ which is $\Gamma(U_w)$.\qed

Since we now know that $\pi_{n+1}(S^3)$ is isomorphic to the set of
equivalence classes of framed $(n-1)$--manifolds embedded in
$S^{n+1}$ with writhe divisible by $w$ we can consider the forgetful
map which ignores the writhe condition on framings.  

\proc{Corollary} The forgetful map is multiplication by $w$.

\prf This follows from observing that the map from $\pi_3(S^3)$ to
$\pi_3(S^2)$ which is given by applying the Thom--Ponrjagin
construction to $U_w$ is $w$ times the Hopf map. But the forgetful
map is effectively composition with this map. \qed

By contrast in the case $w=0$, the cobordism groups are all zero.
This is essentially what theorem 5.13 below says.

Wiest [\Wiest] has extended these results in a number of ways.  He has
shown that the augmented rack space $B_GX$ (where $G$ is $\pi_1(M-L)$)
has the same homotopy type as $\Omega(M^3)$. This implies that $BX$
has the homotopy type of $\Omega(M^3)$ factored by an action of
$G$.  Further, he has shown that, for irreducible links in homotopy
three spheres, the homotopy type of the rack space is determined by
the fundamental group of the link.  Thus for example if $R$ is the
reef knot in $S^3$ (square knot in American English) and $G$ is the
granny knot (both taken with framing zero for definiteness) then
$B\Gamma(R)$ and $B\Gamma(G)$ have the same homotopy type.  In this
case the fundamental rack is a classifying invariant and the two racks
differ.  Thus the homotopy type of the rack space (as distinct from
the combinatorial structure) contains strictly less information that
the rack itself.

\sh{The second homotopy group for links in $S^3$}

Now let $L$ be a framed link in $S^3$. We say that $L$ is {\sl
non-split} or {\sl irreducible} if no embedded 2-sphere in $S^3-L$
divides the components of $L$ into two non-empty subsets. (This is
consistent with the usage of irreducible for links in $M^3$ given
earlier.)

In general a link can be written as a union $L=L_1\cup\ldots\cup L_k$
where each $L_i$ is a maximal irreducible sublink. We call the
sublinks $L_i$ the  {\sl blocks} of $L$. A block is said to be {\sl
trivial} if it is equivalent to the unknot with zero framing.

\proc{Theorem}Let $L$ be a framed link in $S^3$. Then
$\pi_2(B(\Gamma(L)))\cong\Z^p$ where $p$ is the number of non-trivial
blocks of $L$.  Furthermore a basis of $\pi_2(B(\Gamma(L))$ is given
by diagrams representing these blocks.\rm

\rm Paraphrased, the theorem says that any link in $S^3$ with
representation in $\Gamma(L)$ is cobordant respecting the
representation to a unique {\sl standard} link comprising a number of
separate copies of the blocks of $L$.

\prf {\bf Case 1: $L$ is irreducible and non-trivial}

This is a special case of theorem 5.4. 

\rk{Case 2: $L$ is irreducible and trivial}

In this case the homomorphism $\phi\co\pi_3(S^3)\to\F(n+1,\Gamma(L))$
defined in case 1 is surjective by exactly the same argument. But
$\phi[{\rm id}]$ is represented by the trivial link with identity
representation which is null cobordant.

\ppar This completes case 2 and we now turn to the general case
$L=L_1\cup\ldots\cup L_k$ where each $L_i$ is maximal irreducible. We
need the following lemmas.

\proc{Lemma} Let $M$ be an irreducible 3-manifold. Let $M_0$ be $M$
with the interior of $k$ balls $B_1,\ldots,B_k$ removed. Then
$\pi_2(M_0)$ is generated as a $\Z[\pi_1(M_0)]$ module by the spheres
$\d B_i$.

\prf Let ${\widetilde M}$ and ${\widetilde M_0}$ be the universal
covers of $M$ and $M_0$ respectively. Then ${\widetilde M}$ can be
obtained from ${\widetilde M_0}$ by filling in holes with copies
$gB_i$ of $B_i$  one for each element $g$ in $\pi_1(M_0)$, $i=1\ldots
n$. Since $M$ is  irreducible $\pi_2(M)\cong  H_2(\widetilde M)\cong
0$. Then using a Mayer--Vietoris sequence we see that as an abelian
group $\pi_2(M_0)\cong H_2(\widetilde M_0)$ has one generator for
each pair $(g,B_i)$ where $g\in\pi_1(M_0)$.\qed 

\proc{Lemma}Let $M$ be the connected sum $M=M_1\sharp\ldots\sharp
M_k$ of $k$ irreducible 3-manifolds each with non trivial fundamental
group.  Then as a $\Z[\pi_1(M)]$ module $\pi_2(M)$ is generated by
the separating spheres $S_1,\ldots,S_{k-1}$.

\prf Let an element of $\pi_2(M)$ be represented by a map  $f\co
S^2\to M$ of the 2-sphere into $M$ which we may assume is transverse
to the separating spheres. Consider an innermost disc $D$ in $S^2$
which has boundary in the intersection of $f(S^2)$ and $S_i$ say. Let
$D'$ be a (singular) disc in $S_i$ which bounds  $\d D$. Then by the
previous lemma the homotopy class of the sphere  $D\cup D'$ is in the
subgroup generated as a $\pi_1(M)$ module by the separating spheres
$S_1,\ldots,S_{k-1}$. By subtracting this element, we may perform a
homotopy to remove this intersection curve. We can now argue by
induction on the number of intersections.\endprf

Returning now to the proof of the main theorem we look at $\pi_3$ and
attempt to construct as before a map $f\co\re^3\to S^3$. This time
the complement of $L$ may not be an Eilenberg-MacLane space. The
construction of the mapping on the tubular neighbourhood of $\Lambda$
runs as before and there is no obstruction to extending to the
2-skeleton of the complement but obstructions to mapping in the
3-cells may be non zero. By the above lemma it will be sufficient to
consider the case of a single 3-cell, which may be taken to be far
away from $\Lambda$ and where the map on the bounding 2-sphere is as
follows. The map is constant on the equator and the upper hemi-disc
is wrapped around a separating sphere with degree $\pm1$. The
separating sphere contains just one component $L_i$ of $L$. Maps on
great arcs running from the south pole to the equator give the same 
path in $S^3$. We can now map the lower hemi 3-ball to the image of
this path so that the map is constant on the equatorial 2-disc.  At
this point, by considering the upper hemisphere together with the
equatorial 2-disc, we see that the problem is reduced to the  case
where the bounding 2-sphere is mapped to the separating 2-sphere. But
now if we add to $\Lambda$ a copy of $L_i$ in the 3-ball the map
extends in the obvious way. Notice that the representation on the
copy of $\Gamma(L_i)$ which this determines is conjugated by a fixed
element of $\pi_1( L^c)$, 
Call such a link \it conjugated. \rm This
means that the map $f$ can be defined for the extended link. Indeed
we may further extend the link by adding copies of $L$ so that $f$
has degree zero. If we now make the resulting null homotopy
transverse to $L$ we construct a bordism between $[\Lambda,F,\rho]$
and and a number of disjoint copies of conjugated links $L_i$. By
proposition 5.2 and the result in case 1 we can assume that these
added links are labelled by identities and not conjugated. Thus the 
$L_i$ with identity labelling form a generating set for $\pi_2$. Note
that any trivial sublinks can be eliminated by a bordism as in case
2.

Now suppose some non-trivial linear combination of the
$[L_i,F_i,id_i]$ is bordant to zero. We need the following lemma.

\proc{Lemma}Let $X$ be the free product of racks $X_i$ and let $Y_j$
be the union of the orbits in $X$ determined by $X_j$. Then there is
a retraction of $Y_j$ onto $X_j$.

\prf This follows from the definition of the free product of racks
[\FeRo; page 357]. The retraction is defined by setting the action of
other orbits equal to the trivial action. \endprf

The rack homomorphism commutes with the action of the associated
groups and components correspond to orbits of the action. Thus the
null bordism cannot mix components. Now observe that $\Gamma(L)$ is
the free product of the $\Gamma(L_i)$ and therefore there is a
retraction of the orbit determined by $L_i$ onto $\Gamma(L_i)$ by the
lemma. Applying this retraction to the appropriate pieces of the null
bordism we see that $[L_i,F_i,id_i]=0$. This contradicts case 1.

It follows that the $[L_i,F_i,id_i]$ form a basis for $\pi_2$ and
theorem 5.9 is proved.\qed

\sh{Homotopy type of the space of a trivial rack}

Let $X(n)=\{1,\ldots,n\}$ be the trivial rack with $n$ elements; so that
$x^y=x$ for all $x$ and $y$.

\proc{Theorem} The classifying space $B X(n)$ has the homotopy
type of $\Omega(S^2\vee\cdots\vee S^2)$, the loop space on the wedge
of $n$ copies of $S^2$.\endproc

\prf First observe the simple form of the faces in $BX(n)$;
$$\d_i^\ep\row x 1 n =(x_1,\ldots,{\hat x}_i,\ldots,x_n)
\quad\hbox{for}\quad 1\leq i\leq n,\quad\ep\in\{0,1\}$$ Now
$\Omega(S^2\vee\dots\vee S^2)\simeq\Omega S(S^1\vee\dots\vee S^1)$,
but by the James construction $\Omega S(S^1\vee\dots\vee S^1)$ has
the homotopy type of $(S^1\vee\dots\vee S^1)_\infty$ the free monoid
on $(S^1\vee\dots\vee S^1)\setminus\{*\}$. Identify $S^1$ with $I/\d
I$ and let $(k,t)$ denote the point $t$ in the k-th copy of $S^1$ in
$S^1\vee\dots\vee S^1$. Then there is a homeomorphism $B
X(n)\to(S^1\vee\dots\vee S^1)_\infty$ given by $$[\row i1n,\row
t1n]\to(i_1,t_1)\cdot\dots\cdot(i_n,t_n).$$ It follows that
$\pi_k(B(X(n)))$ can be viewed geometrically as bordism classes of
framed codimension two manifolds in $S^{k+1}$ with components
labelled in $\{1,\ldots,n\}$. \qed

\rk{Remark}For $k=2$, by the Hilton-Milnor theorem [\Hilt],
$\pi_2(BX(n))\cong\Z^{n+{n \choose 2}}$ where the second set of
generators are Whitehead products. Geometrically, $\pi_2(BX(n))$ is
interpreted as cobordism classes of links with components labelled by
$n$ distinct labels, or equivalently as a link divided into $n$
disjoint sublinks, and the first $n$ integer invariants are total
writhes of the sublinks and the remaining ${n \choose 2}$ are the
mutual linking numbers. For more details, see [\SBord].

\sh{Homotopy type of the space of a free rack}

\proc{Theorem}Let $FR_n$ denote the free rack on $n$ elements.
Then $BFR_n$ has the homotopy type of $S^1\vee\cdots\vee S^1$, the
wedge of $n$ copies of $S^1$.

\prf Recall from theorem \ref{link-class} that  $\pi_n(BFR_n)\cong
\F(n+1,FR_n)$.

\font\sspec=cmsy10 scaled \magstep 2
\def\s{\hbox{\sspec \char'003}} 

We observe that $FR_n$ is the fundamental rack of $D_n$, which is 
the link comprising $n$
framed points in the 2--disc $D^2$. The proof
of theorem 5.4 shows that if $(\Lambda,FR_n,\rho)$ represents an
element of $\pi_n(BFR_n)$ (for $n\ge2$) then $(\Lambda,F,\rho)$
pulls back from a transverse map of $\re^{n+1}$ to $D_n$. But since
$D^2$ is contractible this map is null homotopic and applying
relative transversality we obtain a null cobordism of 
$(\Lambda,FR_n,\rho)$. By 5.1
$\pi_1(BFR_n)\cong \s_n\Z$ (the free product of $n$ copies of
$\Z$). The result follows. \qed

\rk{Remark}Recall from [\FeRo; page 376] that there is the
concept of an {\sl extended free rack}, which has free operator
generators in addition to the usual free rack generators.  This
can be identified with the fundamental rack of a number of framed
points in an orientable surface.  A similar proof (using the fact
that the higher homotopy groups of a surface vanish) 
then shows that if $F$ is an extended
free rack then $BF$ has the homotopy type of a wedge of circles.  
Moreover both proofs extend with a little care to
arbitrary free (or extended free) racks (in other words there is
no need for the generating sets to be finite). Thus $BF$ has
the homotopy type of a 1--complex where $F$ is any free (or 
extended free) rack.

\sh{A remark on homology groups}

Let $X$ be any rack. Recall that there is the concept of the {\sl
extended} rack space $B_XX$ [\Trunks; example 3.1.1] which is a covering
space of $BX$ [\Trunks; theorem 3.7]. Now there is a chain equivalence
between $C^*(BX,*)$ and $C^{*-1}(B_XX)$ where $*\in BX_0$ is the
unique vertex.  There are corresponding isomorphisms
of homology and cohomology groups (with a shift of one
dimension). This is defined by mapping $(x,x_1,\ldots,x_n)\in
B_XX^{(n-1)}$ to $(x,x_1,\ldots,x_n)\in BX^{(n)}$ using the
description given in [\Trunks; examples 3.4.1 and 2]. Note that
$\d^\ep_{i-1}$ in $B_XX$ coincides with $\d^\ep_i$ in $BX$ for
$i=2,3,\ldots,n$ whilst $\d_1^1=\d_1^0$ in $BX$ (both are given by
$(x_1,x_2,\ldots,x_n)\mapsto(x_2,\ldots,x_n)$) and that these two
cancel out as a pair in the boundary formula.

Thus we have:

\proc{Theorem}The rack space $BX$ of any rack admits a covering space
with the same homology groups but in dimensions all shifted one
lower. \qed\endproc

The chain equivalence can be realised using a map $\mo B_XX\mo\times
S^1\to \mo BX\mo$ defined as follows. Embed the $(n-1)$--cube
$(x,x_1,\ldots,x_n)$ of $B_XX$ as the central $(n-1)$--cube
perpendicular to the first direction in the $n$--cube
$(x,x_1,\ldots,x_n)$ of $BX$. Then use the remaining coordinate to
extend to $B_XX\times[-1,1]$. Since $\d_1^1=\d_1^0$ in $BX$ this map
factors via $B_XX\times[-1,1]/-1\sim1$ that is $B_XX\times S^1$. Then
the above chain equivalence is given by crossing with the fundamental
class of $S^1$ and using this map.

It is an interesting question to characterise spaces which have the
property of admitting a covering space with the same homology groups
shifted one dimension.

\sh{Permutation and Dihedral racks}

Let $\rho\co P\to P$ be a fixed permutation of the set $P$. Then the
{\sl permutation rack} $P_\rho$ is $P$ with $i^j=\rho(i)$ for all
$i,j$.

Combining results of Flower [\Flower] and Greene [\Greene], which were proved
using the cobordism by moves technique of section 4, we have:

\proc{Theorem} For a permutation rack $P_\rho$ \items \item{\rm(1)}
$\pi_2(BP_\rho)$ is freely generated  by one twisted unknot for each
finite orbit, with the number of twists equal to the length of the
orbit, together with a pair of linked unknots (each unknot having a
single twist) for each unordered pair of unequal orbits. \item{\rm(2)}
$H_2(BP\rho)$ is as in $(1)$ except that there is a generator for
each ordered pair of unequal orbits. \qed\enditems\endproc

Let $D_n$ denote the {\sl dihedral rack} on $n$ elements:
$D_n=\{0,1,\ldots,n-1\}$, and $i^j=2j-i \>\hbox{mod}\>n$ for all
$i,j$.

\proc{Theorem}{\rm[\Greene]} $$H_2(BD_n)=\cases{\Z &for $n$ odd \cr
H_2(BD_n)=\Z^4 &for $n=2 \>mod\> 4$ \cr H_2(BD_n)\ge\Z^4 &otherwise
\cr}$$ \qed\endproc

\rk{Remark} A calculation using maple shows $H_2(BD_4)=\Z^4+{\Z_2}^2$, 
so the inequality of the theorem can be strict. 
For further details see [\Flower,\Greene].

\references

\bye